\documentclass[]{interact}

\usepackage{hyperref}
\usepackage{cleveref}
\usepackage{graphicx, graphics,color}
\usepackage{multimedia}
\usepackage{float}
\usepackage{placeins}
\usepackage{mathrsfs}
\usepackage{algorithm}
\usepackage{mathtools}

\usepackage{tikz}
\usepackage{tikz-3dplot}
\usetikzlibrary{arrows}
\usetikzlibrary{calc}
\usetikzlibrary{fit}

\newcommand \udotini {\dot{u}_0}

\newcommand \p {\partial}

\newcommand \y {\mathrm{y}}

\newcommand \R {\mathbb{R}}
\newcommand \N {\mathbb{N}}
\renewcommand \L {\mathrm{L}}
\newcommand \W {\mathrm{W}}
\newcommand \WW {\mathbf{W}}
\newcommand \WWW {\mathbb{W}}

\newcommand \BB {\mathbf{B}}

\newcommand \LL {\mathbf{L}}
\newcommand \LLL {\mathbb{L}}

\renewcommand \H {\mathrm{H}}
\newcommand \HH {\mathbf{H}}

\newcommand \I {\mathrm{I}}

\newcommand \III {\mathbb{I}}
\newcommand \Id {\mathrm{Id}}

\renewcommand \d {\mathrm{d}}

\renewcommand \det {\mathrm{det}}
\newcommand \trace {\mathrm{tr}}

\newcommand \cof {\mathrm{cof}}

\DeclareMathOperator{\divg}{div}

\begingroup\makeatletter\ifx\SetFigFont\undefined%
\gdef\SetFigFont#1#2#3#4#5{%
  \reset@font\fontsize{#1}{#2pt}%
  \fontfamily{#3}\fontseries{#4}\fontshape{#5}%
  \selectfont}%
\fi\endgroup%

\usepackage{epstopdf}
\usepackage[caption=false]{subfig}

\usepackage[numbers,sort&compress]{natbib}
\bibpunct[, ]{[}{]}{,}{n}{,}{,}

\theoremstyle{plain}
\newtheorem{theorem}{Theorem}[section]
\newtheorem{lemma}[theorem]{Lemma}
\newtheorem{corollary}[theorem]{Corollary}
\newtheorem{proposition}[theorem]{Proposition}

\theoremstyle{definition}
\newtheorem{definition}[theorem]{Definition}

\theoremstyle{remark}
\newtheorem{remark}{Remark}

\addtocounter{equation}{-1}
\begin{document}


\title{\textcolor{black}{A hybrid optimal control problem constrained with hyperelasticity and the global injectivity condition}}

\author{
\name{S.~Court\textsuperscript{a,b}\thanks{S.~Court. Email: sebastien.court@uibk.ac.at}}
\affil{\textsuperscript{a}Department of Mathematics, University of Innsbruck, Technikerstrasse 13, 6020 Innsbruck, Austria; \textsuperscript{b}Digital Science Center, University of Innsbruck, Innrain 15, 6020 Innsbruck, Austria.}
}

\maketitle

\begin{abstract}
The purpose of this paper \textcolor{black}{is to address a class of hybrid optimal control problems constrained with hyperelasticity and constant global volume. This type of problems can intervene for example in the mechanical aspects of cardiac activity.} The time deformation of the heart tissue is modeled with the elastodynamics equations dealing with the displacement field as main unknown. These equations are coupled with a pressure whose \textcolor{black}{time variations are aimed to be maximized}. This pressure variable corresponds to a Lagrange multiplier associated with the so-called global injectivity condition, \textcolor{black}{translating the fact that the total volume of the domain remains constant}. We develop an optimal control approach in a general framework that covers in particular the maximization of the variations of this pressure, and also the time the maximum is reached, \textcolor{black}{defining what we call a {\it hybrid} optimal control problem}. Mathematical analysis based on the $\L^p$-parabolic maximal regularity is provided for the state equations and the rigorous derivation of optimality conditions. Numerical simulations for a toy-model illustrate the capacity of the approach.
\end{abstract}

\begin{keywords}
Nonlinear elastodynamics, Global injectivity condition, Hybrid optimal control problems, PDE-constrained optimization problems, $L^p$-maximal parabolic regularity, Finite Element Method.
\end{keywords}

\begin{amscode} (2020): 49K20, 74B20, 74P99, 35K20, 35K55, 35K61, 74F99, 74H20, 74H30, 74S05, 74-10, 65N30.
\end{amscode}

\newpage
\tableofcontents

\section{Introduction} \label{sec-intro}

\textcolor{black}{While elasticity problems often utilizes the hypothesis of small displacements, leading to linear steady equations, some problems require to take into account deformations that cannot be assumed close to the identity, and therefore the corresponding system of partial differential equations governing the time evolution of the displacement is purely nonlinear and unsteady. For example, the heart can be described as an elastic tissue subject to non-trivial deformations. Further, optimal control problems constrained with hyperelasticity can model important tasks like defibrillation of resuscitation\cite{Chernysh, Page2011}. More specifically,} defibrillation of the heart in case of arrhytmias is realized with the use of electric shocks acting on muscular tissues, leading to a reset of its electrical activity and thus the re-oxygenation of the ill area. 
\textcolor{black}{One way to model mathematically such a problem is to consider the monodomain model coupled with the elastodynamics system, as it is proposed in~\cite{Gerach2021}: An electric potential $v_m$ is coupled with two other variables $w_m$ and $s_a$, and the displacement of the elastic material denoted by~$u$ defined on a reference configuration~$\Omega$. The displacement defines the deformed domain at time~$t$ via $\Omega(t) := (\Id +u)(\Omega,t)$. The corresponding coupled system writes then
\begin{subequations} \label{sysmaingenesis}
\begin{eqnarray}
\displaystyle c_m\frac{\p v_m}{\p t} - \frac{1}{\beta}\divg\left(
\sigma_m \nabla v_m \right) + I_{\mathrm{ion}}(v_m,w_m) = \xi & & \text{in }
\Omega(t), \ t\in (0,T), \label{cpr1}\\
\frac{\p w_m}{\p t} = aw_m -bv_m  & & \text{in }
\Omega(t), \ t\in (0,T),\\
\frac{\p s_a}{\p t} = \theta(v_m)\left(k_a(v_m-v_r)-s_a\right)  & & \text{in }
\Omega(t), \ t\in (0,T),\\
\rho \ddot{u} - \divg(\sigma(\nabla u)) = 
\divg(s_a \mathfrak{f} \otimes \mathfrak{f})  & & \text{in }
\Omega, \ t\in (0,T), \label{cpr4}
\end{eqnarray}
\end{subequations}
where $\sigma_m$ is a given conductivity tensor, $\beta$, $a$, $b$, $k_a$ are given coefficients, $I_{\mathrm{ion}}$ is a given nonlinearity (see for example the FitzHugh-Nagumo model~\cite{FitzHugh}), $v_r$ is a residual potential, $\theta$ is a smooth cut-off function, $\rho$ is the density field and $\mathfrak{f}$ is the fiber direction of the elastic tissues. The mapping $u\mapsto \sigma(\nabla u)$ is the second Piola-Kirchhoff stress tensor. Similar models coupling the monodomain model with elasticity were studied in~\cite{Nash2004, Whiteley2007, Goektepe2010, Fritz2014}. In this context visco-elastic models were also considered in \cite{Kaliske2014, Kaliske2015, Propp2020, Zhang2023} for example. From a control perspective, the actuator is exerted as an electric stimulus via the function $\xi$ in the right-hand-side of equation~\eqref{cpr1}. The variable $s_a$ defines the so-called active stress $s_a \mathfrak{f} \otimes \mathfrak{f}$ which corresponds to a volume force in the right-hand-side of~\eqref{cpr4}. This system is strongly coupled, as the displacement $u$ determines the domain~$\Omega(t)$ where the other variables are defined. Mathematical contributions that study the defibrillation problem~\cite{Trayanova2011, CKP11, CKP13, Bendahmane2015, CKP16} or more generally optimal control problems in the context of cardiac electrophysiology~\cite{Breiten1, Breiten2} often omit the coupling terms due to the elastic deformation, in order to treat the monodomain (or the bidomain model) only, set in a reference domain. In the same fashion, but by taking into account the domain deformation, in~\cite{Ambrosini2011} the authors assumed that the variable~$v_m$ dictates the active strain of the muscle. They also consider a simplified version of~\eqref{cpr4}. These simplifications allow to derive mathematical results for this type of problems, as addressing the full problem would require multiphysics coupling, and hence would imply too much technicalities and complexities. We propose a different direction, where the electrical contributions in system~\eqref{sysmaingenesis} are omitted, in order to address equation~\eqref{cpr4} only, by prescribing a control function in its right-hand-side. This would amount assuming that we could prescribe the values of variable~$s_a$, in order to influence the deformation $\Id+u$ of the elastic material over time. Therefore the present article does not have the ambition of modeling defibrillation, but rather aims at developing a rigorous mathematical approach dealing with a class of optimal control problems constrained by a second-order nonlinear partial differential equation like~\eqref{cpr4}.}

\textcolor{black}{We consider in this article an hyperelastic material occupying a domain~$\Omega$, submitted to volume forces. The latter are modeled via a distributed control. The boundary of~$\Omega$ is divided into two parts $\Gamma_D$ and $\Gamma_N$ (see Figure~\ref{fig0}). We assume that the exterior boundary~$\Gamma_D$ remains immobile, and so we impose homogeneous Dirichlet boundary conditions on it. The other boundary~$\Gamma_N$ is not necessarily connected, and is let free. This corresponds to Neumann boundary conditions imposed for the displacement on~$\Gamma_N$. Such a problem could model for example a heart tissue crossed by blood, in contact with the boundary~$\Gamma_N$. The exeterior boundary~$\Gamma_D$ would represent the pericardium, that we can assume to be relatively immobile~\cite{Fritz2014}. Since blood can be considered as an incompressible flow, we then assume additionally that the total volume of the deformed material $(\Id+u)(\Omega,t)$ remains constant over time.}
\hfill \\
\vspace*{30pt}
\begin{center}
\scalebox{0.4}{
\begin{picture}(500,0)
\includegraphics{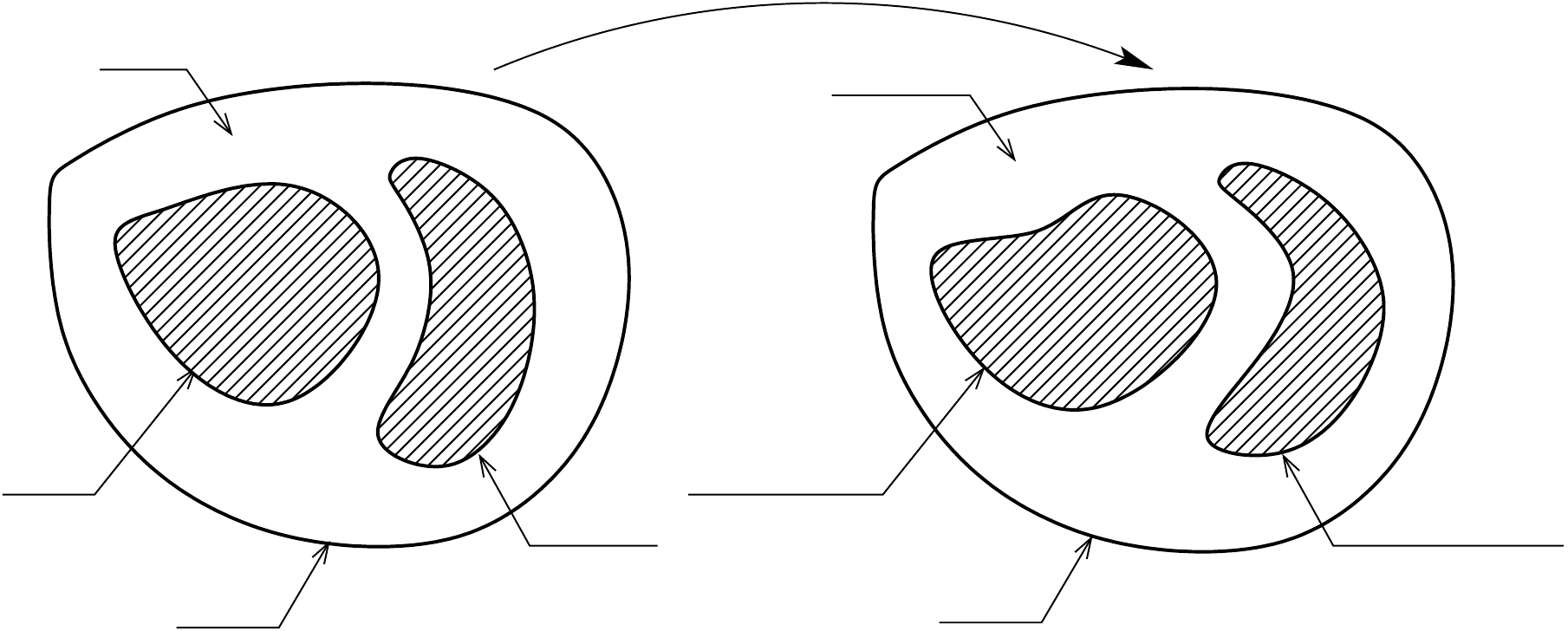}%
\end{picture}%
\setlength{\unitlength}{4144sp}%
\begin{picture}(5566,4712)(2139,-6028)
\put(500,-200){\makebox(0,0)[lb]{\smash{{\SetFigFont{24}{28.4}{\rmdefault}{\mddefault}{\updefault}{\color[rgb]{0,0,0}$\Id+u(\cdot,t)$}%
}}}}
\put(-4600,-1050){\makebox(0,0)[lb]{\smash{{\SetFigFont{26}{28.4}{\rmdefault}{\mddefault}{\updefault}{\color[rgb]{0,0,0}$\Omega$}%
}}}}
\put(1400,-1200){\makebox(0,0)[lb]{\smash{{\SetFigFont{22}{24.4}{\rmdefault}{\mddefault}{\updefault}{\color[rgb]{0,0,0}$(\Id+u)(\Omega,t)$}%
}}}}
\put(-100,-4700){\makebox(0,0)[lb]{\smash{{\SetFigFont{22}{24.4}{\rmdefault}{\mddefault}{\updefault}{\color[rgb]{0,0,0}$(\Id+u)(\Gamma_{N,1},t)$}%
}}}}
\put(6200,-5150){\makebox(0,0)[lb]{\smash{{\SetFigFont{22}{24.4}{\rmdefault}{\mddefault}{\updefault}{\color[rgb]{0,0,0}$(\Id+u)(\Gamma_{N,2},t)$}%
}}}}
\put(2600,-5850){\makebox(0,0)[lb]{\smash{{\SetFigFont{22}{24.4}{\rmdefault}{\mddefault}{\updefault}{\color[rgb]{0,0,0}$\Gamma_D$}%
}}}}
\put(-4100,-5900){\makebox(0,0)[lb]{\smash{{\SetFigFont{22}{24.4}{\rmdefault}{\mddefault}{\updefault}{\color[rgb]{0,0,0}$\Gamma_D$}%
}}}}
\put(-5600,-4700){\makebox(0,0)[lb]{\smash{{\SetFigFont{22}{24.4}{\rmdefault}{\mddefault}{\updefault}{\color[rgb]{0,0,0}$\Gamma_{N,1}$}%
}}}}
\put(-5600,-4700){\makebox(0,0)[lb]{\smash{{\SetFigFont{22}{24.4}{\rmdefault}{\mddefault}{\updefault}{\color[rgb]{0,0,0}$\Gamma_{N,1}$}%
}}}}
\put(-950,-5150){\makebox(0,0)[lb]{\smash{{\SetFigFont{22}{24.4}{\rmdefault}{\mddefault}{\updefault}{\color[rgb]{0,0,0}$\Gamma_{N,2}$}%
}}}}
\end{picture}%
}
\nopagebreak
\begin{figure}[H]
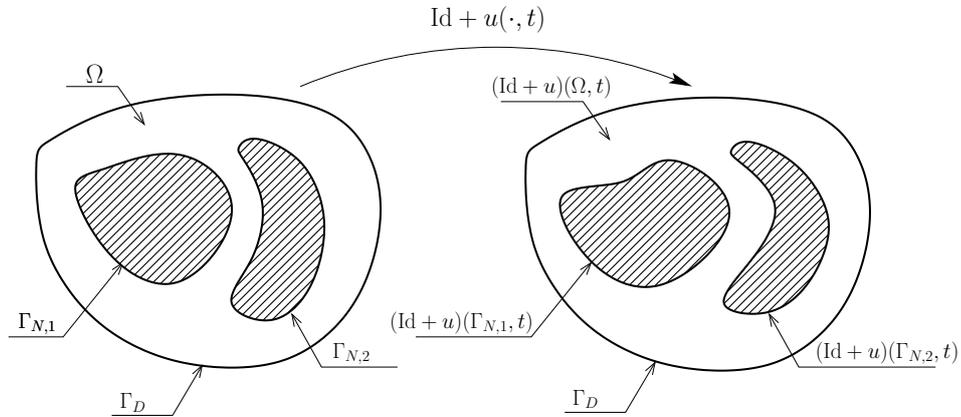

\centering
\caption{Slice representation of \textcolor{black}{the reference configuration $\Omega$ and its deformation $(\Id+u)(\Omega,t)$ at time~$t$}. We impose $u=0$ on the exterior boundary $\Gamma_D$, and the interior \textcolor{black}{boundary $\Gamma_N = \Gamma_{N,1} \cup \Gamma_{N,2}$} is subject to Neumann conditions.}
\label{fig0}
\end{figure}
\end{center}

\FloatBarrier

\subsection{Elastodynamics equations with damping and global injectivity constraint.} In a smooth bounded domain $\Omega$ of $\R^d$ ($d \leq 3$) we consider the elastodynamics system with damping and global volume preserving constraint. The unknown are a displacement field denoted by $u$, and a pressure variable denoted by~$\mathfrak{p}$. The pressure~$\mathfrak{p}$ is a Lagrange multiplier associated with the volume preserving constraint, and depends only on the time variable. \textcolor{black}{We shall not confound this pressure with the hydrostatic pressure associated with the incompressibility condition ($\det(\I+\nabla u) \equiv 1$) commonly considered in fluid mechanics, but also in~\cite{Ambrosini2011, Fritz2014, Bendahmane2019, Bendahmane2024} for elasticity models where cardiac cells are assumed to be composed mainly of water. In the present article we do not make this strong assumption, and prefer to take into account this global volume constraint.} The control operator is distributed, represented by a smooth mapping $\xi \mapsto f(\xi)$ \textcolor{black}{that play the role of a volume force}. The initial state is given by the couple $(u_0,\udotini)$. The couple $(u,\mathfrak{p})$ satisfies the following system: 
\begin{subequations} \label{sysmain}
\begin{eqnarray}
\displaystyle \ddot{u} -  \kappa\Delta \dot{u} -\divg \left(\sigma(\nabla u) \right) = f(\xi) & & \text{in } \Omega \times (0,T), \label{sysmain1}\\
\displaystyle\kappa \frac{\p \dot{u}}{\p n} + \sigma(\nabla u)n 
+ \mathfrak{p}\, \cof (\I +\nabla u)n = g & & 
\text{on } \Gamma_N \times (0,T), \label{sysmain2} \\
\displaystyle\int_{\Omega} \det(\I+\nabla u)\, \d \Omega = 
\int_{\Omega} \det(\I+\nabla u_0)\, \d \Omega & & \text{in } (0,T), \label{sysmain4}\\
u=0 & & \text{on } \Gamma_D \times (0,T),\label{sysmain3}\\
\displaystyle u(\cdot,0) = u_0, \quad \dot{u}(\cdot,0) = \udotini & & \text{in } \Omega.
\label{sysmain5}
\end{eqnarray}
\end{subequations}
\textcolor{black}{The physical derivation of this system is provided in section~\ref{sec-deriv-Lag}.} The tensor field $\sigma(\nabla u)$ is derived from the elasticity model that we adopt. The damping term $\kappa\Delta \dot{u}$ is added for the sake of mathematical convenience. Indeed, it enables us to use well-established results for parabolic equations, while the original system is hyperbolic and nonlinear. The constraint~\eqref{sysmain4} is the so-called Ciarlet-Ne\v{c}as global injectivity condition, studied in~\cite{Necas1987}. It represents the time preservation of the total volume of $\Omega$ under the deformation $\Id +u$. The right-hand-side~$g$ in~\eqref{sysmain2} represents possible given surface forces. Local-in-time wellposedness for system~\eqref{sysmain} has been established in~\cite{Court2023}. In the present paper, for the analysis we prefer to assume smallness of the data rather than relying on small time~$T>0$, as it offers more straightforward wellposedness results which rely on the application of the inverse function theorem. The derivation of other type of existence results (like global wellposedness), that would not require smallness assumptions for example, would go beyond the scope of the present article. Note that under smallness assumptions we can keep the displacement~$u$ small enough, so that the deformation~$\Id+u$ remains invertible, that will be assumed in the rest of the paper.\\
\textcolor{black}{Without damping, and without pressure variable, wellposedness questions for nonlinear elastic waves are well-known to be difficult to address: Global-in-time existence of solutions for equations similar to~\eqref{sysmain} was obtained in~\cite{Sideris1996, Agemi2000} by assuming the {\it null condition} and in~\cite{Sideris2000} by assuming that the stored energy of the elastic domain obeys a {\it nonresonance} condition. In the context incompressibility condition $\det(\I+\nabla u) \equiv 1$, existence results were obtained in~\cite{Ebin1993, Ebin1996, Becca2005, Becca2007, Lei2016} under smallness assumptions on the data.\\
Note that, beside the case of heart tissues, system~\eqref{sysmain} can also model any mechanical problem where a compressible elastic material encloses an incompressible fluid like a gas (pneumatic system) or a liquid (hydraulic system).}

\subsection{A hybrid optimal control problem.}
We assume that the control function $\xi$ acts on a subdomain $\omega \subset \Omega$. The goal is to maximize an objective function~$\Phi^{(1)}$ at time~$\tau$. Consider a cost functional $c$, the general optimal control problem that the present paper proposes to address is the following:
\begin{equation} \tag{$\mathcal{P}$}
\left\{ \begin{array}{l}\displaystyle \max_{\xi \in \mathcal{X}_{p,T}(\omega), \tau \in (0,T)} \left(
J(\xi,\tau):= \int_0^T c(u,\dot{u},\xi) \, \d t + \phi^{(1)}(u,\dot{u})(\tau) + \phi^{(2)}(u,\dot{u})(T)\right), \\[10pt]
\text{subject to~\eqref{sysmain}.}
\end{array} \right. \label{mainpb}
\end{equation}
The control space~$\mathcal{X}_{p,T}(\omega)$ will be defined later. The functional~$\phi^{(2)}$ represents the terminal cost, and~$\phi^{(1)}$ the objective functional that we want to maximize at time~$\tau$. Note that in the formulation of Problem~\eqref{mainpb}, the state variables are coupled with the time parameter $\tau$ to be optimized. The values of the state variables at time~$\tau$ depend on the control function $\xi$. Therefore the two control variables, $\xi$ and $\tau$ are coupled in a complex manner, defining why we call such a problem a {\it hybrid} optimal control problem. This type of problems were treated in~\cite{Maxmax1, Maxmax2} when dealing with a time-parameter, and in~\cite{Maxmax3} when dealing with a space parameter. \textcolor{black}{One of the novelties lies in the fact that the state equation is a second-order time evolution equation. We notice very few contributions dealing with optimal control problems constrained with nonlinear unsteady elastic models: The contributions~\cite{Askar1991, Sun2009, Barsegyan2012, Kimmerle2018, Kroener2011} deal with linear wave-type equations, while~\cite{Kroener2012, Kimmerle2018} deal with numerical realizations only. Let us mention that static nonlinear problems were addressed in~\cite{Givoli2001, Frutos2020, Ortigosa2021}, and in~\cite{Lubkoll2014} under the local orientation-preserving condition $\det(\I+\nabla u) >0$. Another novelty of the present article is the specific form of the functional~$\phi^{(1)}$, explained in section~\ref{subsec-choice-func}.}

\subsection{On the choice of the objective function.} \label{subsec-choice-func}
\textcolor{black}{Cardiopulmonary resuscitation (CPR) consists in applying chest compressions in order to re-ignite heart beats. The goal is to maximize the aortic diastolic pressure~\cite{Page2011} (see the paragraph entitled {\it CPR: Push hard, push fast, and minimize interruptions}): These compressions must be {\it rapid and forceful} in order to maximize {\it the return of blood to the heart}. Therefore maximizing the variations of the pressure is of interest. Further, optimizing the time at which these variations can be maximal is also of interest. Let us then choose to} maximize the variations of the pressure~$\mathfrak{p}$, on a short time interval $(\tau, \tau +\varepsilon)$, where the parameter $\tau$ is let free and also chosen optimally \textcolor{black}{over all the possible times of~$(0,T)$}, while $\varepsilon >0$ is fixed and chosen as small as possible. For several reasons, one would not make tend $\varepsilon$ to $0$, as the pressure variable may not necessarily be differentiable in time in view of the regularity of the data and the functional framework in which the unknown of system~\eqref{sysmain} will be considered. Besides, for practical realization it may not be possible in practice to choose technically~$\varepsilon$ as small as desired. Problem~\eqref{mainpb} would then involve the following objective function at time~$\tau$:
\begin{equation}
\phi_1(u,\dot{u})(\tau) = 
\frac{\mathfrak{p}(\tau+\varepsilon) -\mathfrak{p}(\tau)}{\varepsilon}. \label{sysmainpressure}
\end{equation}
From equation~\eqref{sysmain2} (with $g=0$) the variable $\mathfrak{p}$ is indeed a function of~$(u,\dot{u})$: 
\begin{equation*}
\mathfrak{p} = \displaystyle 
-\frac{1}{|(\Id+u)(\Gamma_N)|} \int_{\Gamma_N} \left(
\kappa \frac{\p \dot{u}}{\p n} + \sigma(\nabla u)n
\right)\d \Gamma_N, 
\end{equation*}
with
\begin{equation*}
|(\Id+u)(\Gamma_N)| = \int_{\Gamma_N} |\cof(\I+\nabla u)n|_{\R^d} \d \Gamma_N.
\end{equation*}
In~\eqref{sysmainpressure}, assuming that the parameter $\varepsilon$ could be supposed to tend towards zero would lead to maximize directly the time-derivative of the pressure. This would require more regularity to consider for the state variables, and would also involve to incorporate $u$, $\dot{u}$ and~$\ddot {u}$ in the objective function, leading to undesirable complexities.\\
\textcolor{black}{Taking into account this type of functionals, involving originally a dual variable, constitutes another novelty of our approach compared with contributions dealing with optimal control problems constrained with elasticity models.}

\subsection{Methodology.} 
The type of objective functional that we consider in Problem~\eqref{mainpb} requires that the state variables~$u$ and~$\dot{u}$ are continuous in time, with values in smooth trace spaces. Therefore a strong functional framework is adopted, corresponding to the so-called $L^p$-maximal parabolic regularity, leading us to assume that a solution $u$ of~\eqref{sysmain} satisfies
\begin{equation*}
\dot{u}\in \L^p(0,T; \WW^{2,p}(\Omega)) \cap \W^{1,p}(0,T;\LL^p(\Omega)),
\end{equation*}
with $p>3$. See sections~\ref{sec-notation-func} and~\ref{sec-well} for more details. With the help of results obtained in~\cite{Pruess2002, Arendt2007}, we first study system~\eqref{sysmain} linearized around $0$, in order to deduce local existence of solutions for~\eqref{sysmain}, while assuming that the data -- initial conditions and right-and-sides -- small enough. We also prove wellposedness for a non-autonomous linear system (namely system~\eqref{sysmain} linearized around a non-trivial state) that will be used for deriving rigorously optimality conditions. From there we are in position to address the question of optimality conditions for Problem~\eqref{mainpb}. Due to the (lack of) regularity of the given right-hand-sides, the time-derivative of the state variable, namely~$\ddot{u}$, is not continuous in time, and therefore it is not possible to determine the optimal parameter~$\tau$ directly by deriving in time the function~$\phi^{(1)}(u,\dot{u})$. In order to uncouple the time parameter and the state variables, we first introduce a change of variable, and reformulate Problem~\eqref{mainpb} in terms of the new variables, like in~\cite{Maxmax1, Maxmax2} (see Figure~\ref{fig-diagram}). Next we introduce an adjoint system whose wellposedness is obtained by transposition, via the non-autonomous linear system aforementioned. Necessary first-order optimality conditions are then derived and expressed in terms of the state, the control function, the time parameter~$\tau$ and the adjoint state. Further, the optimality conditions so obtained are re-expressed in terms of the original state variables, by reversing the change of variables previously used. This is convenient for numerical implementation, even if this was not the approach adopted in~\cite{Maxmax1, Maxmax2}. \textcolor{black}{These new expressions for the optimality conditions can be more suitable for performing numerical illustrations. The latter are performed on a 1D model with finite elements approximation combined with an augmented Lagrangian technique.}

\begin{figure} \label{fig-diagram}
\begin{equation*}
\boxed{
\begin{array} {rcclclcll}
& (u,\dot{u}) & \longrightarrow & (\mathcal{P}) & \longrightarrow & J & \text{\scalebox{1.2}{$\overset{\mathop{?}}{\dashrightarrow}$}} & \nabla J & \\
\text{\scriptsize{ }} \hspace*{-15pt} & \mu \big\downarrow & &
\hspace{0pt} \mu\big\downarrow & & \hspace{-5pt} \mu\big\downarrow & &
\hspace*{5pt}\big\uparrow\mu^{-1} \\ [5pt]
& (\tilde{u},\dot{\tilde{u}}) & \longrightarrow & (\tilde{\mathcal{P}}) & \longrightarrow & \tilde{J} &
\longrightarrow & \nabla \tilde{J} &
\end{array}}
\end{equation*}
\centering
\caption{Transformation of the variables~$(u,\dot{u})$ via a change of variable~$\mu$, leading to a new problem~$(\tilde{\mathcal{P}})$ with a new objective function~$\tilde{J}$, for which the derivation of optimality conditions is tractable.}
\end{figure}
\FloatBarrier

The paper is organized as follows: Notation, assumptions and functional framework are provided in section~\ref{sec-notation}. In particular, we show with classical examples in section~\ref{sec-app-energies} that the assumptions made on the strain energy in section~\ref{sec-notation-ass} are reasonable. \textcolor{black}{Derivation of system~\eqref{sysmain} from the principle of least action is provided in section~\ref{sec-deriv-Lag}.} Section~\ref{sec-well} is dedicated to wellposedness results in the context of the $L^p$-maximal parabolic regularity. In section~\ref{sec-well1} we study system~\eqref{sysmain} linearized around~$0$, leading in section~\ref{sec-well2} to local existence results for system~\eqref{sysmain} and also for its linearized version around some non-trivial states. Solutions for a general adjoint system are studied in section~\ref{sec-well3}. Next section~\ref{sec-optcond} is devoted to the derivation of optimality conditions: Problem~\eqref{mainpb} is transformed in section~\ref{sec-transformation}, the corresponding control-to-state mapping is studied in section~\ref{sec-cts}, and the main results are obtained in sections~\ref{sec-optcond_1st} and~\ref{sec-optcond-final}. Numerical illustrations are presented in section~\ref{sec-num}. Final comments are given in the conclusion (section~\ref{sec-conclusion}). In the Appendix, we present modeling aspects of the problem, in particular the formal derivation of the optimality conditions from a Lagrangian mechanics perspective (section~\ref{sec-app-Lag}), the investigation of classical elasticity models with regards to the assumptions made in this article (section~\ref{sec-app-energies}) and the physical description of the control operator (section~\ref{sec-app-control}).

\section{Notation and assumptions} \label{sec-notation}
Let us introduce the notation for functional spaces and assumptions on the data and the elasticity model. The reader is invited to refer to the present section when reading the rest of the paper.

\subsection{General notation} \label{sec-notation-general}
We denote by $u\cdot v$ the inner product between two vectors $u$, $v \in \R^d$, and the corresponding Euclidean norm by $|u|_{\R^d}$. We define the tensor product $u\otimes v \in \R^{d\times d}$, such that $(u\otimes v)_{ij} := u_i v_j$. The inner product between two matrices $A$, $B\in \R^{d\times d}$ is given by $A:B = \mathrm{trace}(A^TB)$, and we recall that the associated Euclidean norm satisfies $| AB |_{\R^{d\times d}} \leq | A|_{\R^{d\times d}} |B|_{\R^{d\times d}}$. The tensor product between matrices is denoted by $A \otimes B \in \R^{d\times d \times d \times d}$, such that for all matrix $C \in \R^{d\times d}$ we have $(A\otimes B)C := (B:C)A \in \R^{d\times d}$.\\

\noindent\textbf{On the cofactor matrix.} We denote by $\cof (A)$ the cofactor matrix of any matrix field $A$. Recall that this is a polynomial function of the coefficients of $A$. When $A$ is invertible, the following formula holds
\begin{equation*}
\cof (A) = (\det (A))A^{-T} .
\end{equation*}
Recall that $H \mapsto (\cof A):H$ is the differential of $A \mapsto \det(A)$ at point $A$. Further, recall the differential of $A \mapsto \cof(A)$, given by the following formula
\begin{equation*}
D_A(\cof).H = 
\frac{1}{\det A} 
\left(\left(\cof(A) \otimes \cof(A) \right)H 
- \cof (A) H^T \cof (A)\right),
\end{equation*}
for all matrix $H \in \R^{d\times d}$.

\subsection{Geometric assumptions and the global injectivity condition} \label{subsec-geo}
The domain $\Omega \subset \R^d$ (with~$d=2$ or~$3$) is assumed to be smooth and bounded. Its boundary~$\p \Omega$ is made of two smooth parts~$\Gamma_D$ and~$\Gamma_N$ such that $\Gamma_D \cap \Gamma_N = \emptyset$ (see Figure~\ref{fig0}), and their respective surface Lebesgue measures are positive. We assume that~$\Gamma_D$ and~$\Gamma_N$ are smooth in the sense that the surfaces $\Gamma_D$ and $\Gamma_N$ are {\it regular}, meaning that at any point of $\Gamma_D$ and $\Gamma_N$ we can define a tangent plane. Therefore the outward unit normal of~$\p \Omega$ is well-defined. On~$\Gamma_N$, we will assume that~$n\in \WW^{2-1/p,p}(\Gamma_N)$.\\
The deformation gradient tensor associated with the displacement field~$u$ is denoted by
\begin{equation*}
\Phi(u) = \nabla (\Id+u) = \I+\nabla u.
\end{equation*}
Equation~\eqref{sysmain4} translates the fact that the total volume of~$\Omega$ must remain constant over time. Differentiating this equality in the direction~$v$ yields
\begin{equation*}
\int_{\Omega} \det(\Phi(u)) \d \Omega = 
\int_{\Omega} \det(\Phi(u_0)) \d \Omega
\quad \Rightarrow \quad
\int_{\Omega} \cof(\Phi(u)): \nabla v \, \d \Omega = 0.
\end{equation*}
Further, using the Piola's identity, we have $\divg(\cof(\Phi(u))^T v) = \cof(\Phi(u)): \nabla v$, and then the divergence formula enables us to rewrite the quantity above as an integral over~$\Gamma_N$ only, as follows
\begin{equation}\label{myPiola}
\int_{\Omega} \cof(\Phi(u)): \nabla v \, \d \Omega = 
\int_{\Gamma_N} v\cdot \cof(\Phi(u))n \, \d \Gamma_N = 0,
\end{equation}
if we assume $v_{|\Gamma_D} = 0$.
In particular, equation~\eqref{sysmain4} can be equivalently replaced by its time-derivative, namely
\begin{equation}
\int_{\Gamma_N} \dot{u}\cdot \cof(\Phi(u))n \, \d \Gamma_N = 0. \label{sysmain4bis}
\end{equation}
dealing with the boundary~$\Gamma_N$ only.

\subsection{Functional spaces} \label{sec-notation-func}
Throughout we consider the exponent $p>3$. In order to distinguish scalar fields, vector fields and matrix fields, we use the following notation
\begin{equation*}
\L^{p}(\Omega) = \left\{\varphi :\Omega \rightarrow \R \mid \int_{\Omega} |\varphi|_{\R}^p\d \Omega <\infty \right\}, \quad
\LL^{p}(\Omega) = [\L^p(\Omega)]^d, \quad
\LLL^{p}(\Omega) = [\L^p(\Omega)]^{d\times d},
\end{equation*}
that we transpose by analogy to other types of Lebesgue and Sobolev spaces. Denote
\begin{equation*}
\WW^{1,p}_{D,0}(\Omega) := \left\{ 
\varphi \in \WW^{1,p}(\Omega) \mid \quad  \varphi_{|\Gamma_D} = 0
\right\}.
\end{equation*}
The displacement $u$ and its time-derivative $\dot{u}$ are considered in the spaces given below:
\begin{equation*}
\begin{array} {l}
 u \in \mathcal{U}_{p,T}(\Omega)  :=  \W^{1,p}(0,T; \WW^{2,p}(\Omega)\cap \WW^{1,p}_{0,D}(\Omega)) 
\cap \W^{2,p}(0,T; \LL^p(\Omega )), \\
 \dot{u} \in \dot{\mathcal{U}}_{p,T}(\Omega)  :=  
\L^p(0,T; \WW^{2,p}(\Omega)\cap \WW^{1,p}_{0,D}(\Omega)) 
\cap \W^{1,p}(0,T;\LL^p(\Omega)).
\end{array}
\end{equation*}
Given $r\in (1,\infty)$, we denote by $r'$ its dual exponent satisfying $1/r+1/{r'} = 1$. The trace space for $\dot{u} \in \dot{\mathcal{U}}_{p,T}(\Omega)$ involves the Besov spaces obtained by real interpolation as $\displaystyle \left(\LL^p(\Omega);\WW^{2,p}(\Omega)\right)_{1/{p'},p} =: \BB^{2/{p'}}_{pp}(\Omega)$ and $\big(\LL^p(\Omega) ;\WW^{1,p}_{0,D}(\Omega)\big)_{1/{p'},p} =: \mathring{\BB}^{1/{p'}}_{pp}(\Omega)$, which coincide with $\WW^{2/{p'},p}(\Omega)$ and $\WW^{1/{p'},p}_{0,D}(\Omega)$, respectively. See for instance~\cite{Triebel}. The initial conditions are assumed to lie in the trace space $
\left\{		(u(0),\dot{u}(0)) \mid u \in \mathcal{U}_{p,T}(\Omega)\times \dot{\mathcal{U}}_{p,T}(\Omega)	\right\}$, namely:
\begin{equation*}
(u_0,\udotini) \in \mathcal{U}_p^{(0,1)}(\Omega)  := 
\left(\WW^{2,p}(\Omega)\cap \WW^{1,p}_{0,D}(\Omega)  \right)
\times
\left(\WW^{2/{p'},p}(\Omega) \cap \WW^{1/{p'},p}_{0,D}(\Omega)\right).
\end{equation*}
We refer to~\cite{Chill2005} and~\cite[section~6]{Arendt2007} for more details. The choice of such a strong functional frame is motivated by the fact that the trace space described above guarantees that the gradient of the displacement is continuous in space. Further, introduce the following spaces
\begin{equation*}
\begin{array} {l}
\mathcal{F}_{p,T}(\Omega) := \L^p(0,T;\LL^p(\Omega)),\\
\mathcal{G}_{p,T}(\Gamma_N) := \W^{1/2-1/{2p},p}(0,T; \LL^p(\Gamma_N)) \cap \L^p(0,T;\WW^{1-1/p,p}(\Gamma_N)),\\
\mathcal{H}_{p,T}(\Gamma_N) := \W^{1-1/{2p},p}(0,T; \LL^p(\Gamma_N)) \cap \L^p(0,T;\WW^{2-1/p,p}(\Gamma_N)), \\
\mathcal{H}_{p,T} := \W^{1-1/{2p},p}(0,T; \R).
\end{array}
\end{equation*}
Following~\cite{Pruess2002}, the Neumann boundary condition~\eqref{sysmain2} is considered in~$\mathcal{G}_{p,T}(\Gamma_N)$, and the trace of~$\dot{u}$ on~$\Gamma_N$ is considered in~$\mathcal{H}_{p,T}(\Gamma_N)$. More precisely, we recall the following boundary trace embedding (see for example~\cite[Lemma~3.5]{Denk2006}):
\begin{equation}
\|v_{|\Gamma_N}\|_{\mathcal{H}_{p,T}(\Gamma_N)}
\leq C\|v \|_{\dot{\mathcal{U}}_{p,T}(\Omega)},
\label{est-trace-emb}
\end{equation}
where the constant $C>0$ is independent of~$v$. The constraint~\eqref{sysmain4bis} involves the trace of~$\dot{u}$ on~$\Gamma_N$, and the space~$\mathcal{H}_{p,T}$ is where this constraint~\eqref{sysmain4bis} is considered (namely constraint~\eqref{sysmain} derived in time). 
Further, we also recall this other trace embedding:
\begin{equation}
\displaystyle
\left\|\frac{\p v}{\p n}\right\|_{\mathcal{G}_{p,T}(\Gamma_N)}
\leq C\|v \|_{\dot{\mathcal{U}}_{p,T}(\Omega)}.
\label{est-trace-emb2}
\end{equation}
In the Hilbert case, we will need the following estimate:
\begin{equation}
\left\|\frac{\p v}{\p n}\right\|_{\W^{1/(2p')}(0,T;\mathbf{H}^{1/2-1/p}(\Gamma_N))'}
\leq C\| v\|_{\L^p(0,T;\mathbf{H}^2(\Omega))\cap \W^{1,p}(0,T;\mathbf{L}^2(\Omega))}.
\label{est-trace-emb3}
\end{equation}
Note that the trace space of~$\L^p(0,T;\mathbf{H}^2(\Omega))\cap \W^{1,p}(0,T;\mathbf{L}^2(\Omega))$ coincides with~$\mathbf{H}^{2/(p')}(\Omega)$.\\
Finally, the pressure variable~$\mathfrak{p}$ that appears in the Neumann condition will be considered 
such that
\begin{equation*}
\mathfrak{p} \in \mathcal{P}_{p,T} := \W^{1/2-1/{2p},p}(0,T; \R) =
\W^{1/(2p'),p}(0,T; \R).
\end{equation*}

\subsection{Assumptions on the strain energy and operator notation} \label{sec-notation-ass}
In th rest of the paper we will assume that the functionals which appear in Problem~\eqref{mainpb}, namely~$c$, $\phi^{(1)}$ and~$\phi^{(2)}$, are Fr\'echet-differentiable on~$\mathcal{U}_{p,T}(\Omega) \times \dot{\mathcal{U}}_{p,T}(\Omega)\times \mathcal{X}_{p,T}(\omega)$ for functional~$c$, and on $\mathcal{U}_{p,T}(\Omega) \times \dot{\mathcal{U}}_{p,T}(\Omega)$ for functionals~$\phi^{(1)}$ and~$\phi^{(2)}$. Let us give the assumptions that we make on the other operators of the problem.

Recall that for $p>d$, the space $\WWW^{1,p}(\Omega)$ is an algebra. In particular, there exists a positive constant $C$, depending only on $\Omega$ and $p$, such that for all $A, \, B \in \WWW^{1,p}(\Omega)$, we have 
\begin{equation}
\|AB \|_{\WWW^{1,p}(\Omega)}  \leq 
C \| A \|_{\WWW^{1,p}(\Omega)} \|B  \|_{\WWW^{1,p}(\Omega)}.
\label{W-algebra}
\end{equation}
See for instance~\cite[Lemma~A.1]{BB1974}. Therefore the different products between the elasticity-related tensors will be mainly understood in~$\WWW^{1,p}(\Omega)$. Recall the expression of the deformation gradient $\Phi(u) = \I + \nabla u$ introduced previously. The strain energy of the elastic material is denoted by $\mathcal{W}$, and is a function of the Green--Saint-venant strain tensor
\begin{equation*}
E(u) := \frac{1}{2}\left(\Phi(u)^T \Phi(u) - \I\right) = 
\frac{1}{2}\left((\I+\nabla u)^T(\I+\nabla u) - \I \right).
\end{equation*}
We denote classically~\cite{Ciarlet} by $\check{\Sigma}$ the differential of~$\mathcal{W}$:
\begin{equation*}
\check{\Sigma}(E) = \frac{\p \mathcal{W}}{\p E}(E),
\end{equation*}
and by $\Sigma$ its composition by~$E(u)$ as follows
\begin{equation*}
\Sigma(u) := \frac{\p \mathcal{W}}{\p E}(E(u)).
\end{equation*}
We further introduce \textcolor{black}{the second Piola-Kirchhoff stress tensor, that is}
\begin{equation*}
\sigma(\nabla u) = (\I+\nabla u )\Sigma(u) =
\Phi(u)\Sigma(u) = 
 (\I + \nabla u) \frac{\p \mathcal{W}}{\p E}(E(u)),
\end{equation*}
that is the operator which appears in system~\eqref{sysmain}. Note that $\sigma$ is a function of $\nabla u$ only, since the strain energy is chosen to be a function of the Green -- St-Venant strain tensor $E(u)$, which is itself a function of~$\nabla u$ only. The derivation of system~\eqref{sysmain} from $\mathcal{W}$ is given in section~\ref{sec-deriv-Lag}.

\textcolor{black}{
\subsubsection{Differentiation of nonlinear operators} \label{sec-linear2}}
The tensor $E$ linearized around $u$ in the direction $v$ is given by
\begin{equation*}
E'(u).v = \frac{1}{2}\left(\Phi(u)^T\nabla v + \nabla v^T\Phi(u) \right).
\end{equation*}
The linearized systems, around $0$ in section~\ref{sec-well1} and around some unsteady state $u$ in sections~\ref{sec-well2} and~\ref{sec-well3}, involve the differentials of $\sigma(\nabla u)$ and $\cof(\Phi(u))$ (with respect to~$\nabla u$), denoted respectively by~$\sigma_L(\nabla u)$ and~$\sigma_N(\nabla u)$, and given as follow
\begin{subequations} \label{tensor-linear}
\begin{eqnarray}
\sigma_L(\nabla u).\nabla v & = & \nabla v \Sigma(u) + (\I+\nabla u)\frac{\p^2 \mathcal{W}}{\p E^2}(E(u)).(E'(u).v), \label{tensor-lin1}\\
\sigma_N(\nabla u).\nabla v 
& = & \frac{1}{\det \Phi(u)} 
\left(\left(\cof \Phi(u) \otimes \cof \Phi(u) \right)\nabla v 
- \cof \Phi(u) \nabla v^T \cof \Phi(u)\right). 
 \label{tensor-lin2}
\end{eqnarray}
\end{subequations}
\textcolor{black}{Given the assumptions~$\mathbf{A2}$ and~$\mathbf{A3}$ of section~\ref{subsec-assx} below, we note that~$\sigma_L(\nabla u)$ is symmetric}. A variational formulation of its expressions gives, for all vector field $w$
\begin{equation*}
\begin{array} {rcl}
(\sigma_L(\nabla u).\nabla v): \nabla w & = & (\nabla v \Sigma(u)):\nabla w + 
\displaystyle
\left(\frac{\p^2 \mathcal{W}}{\p E^2}(E(u)).(E'(u).v)\right): (E'(u).w). 
\end{array}
\end{equation*}
The operator~$ v \mapsto (\sigma_N(\nabla u).\nabla v)n$ is symmetric too. Indeed, for all smooth test function~$\zeta$ such that $\zeta_{|\Gamma_D} = 0$, we first express, using the divergence formula and the Piola's identity 
\begin{equation*}
\begin{array} {rcl} 
\displaystyle
\int_{\Gamma_N} \zeta \cdot \cof(\Phi(u))n \d \Gamma_N =  
\displaystyle
\int_{\Gamma_N} (\cof(\Phi(u))^T\zeta) \cdot n \d \Gamma_N & = &
\displaystyle
\int_{\Omega} \divg \left(\cof(\Phi(u))^T \zeta\right) \d \Omega \\[10pt]
& = & \displaystyle
\int_{\Omega} \cof(\Phi(u)): \nabla \zeta \d \Omega.
\end{array}
\end{equation*}
Next, differentiating the left- and right-hand-sides of this equality yields\begin{small}
\begin{equation*}
\begin{array} {rcl} 
& & \displaystyle\int_{\Gamma_N} \zeta \cdot (\sigma_N(\nabla u).\nabla v)n\, \d \Gamma_N  
 =  \displaystyle
\int_{\Omega} (\sigma_N(\nabla u).\nabla v): \nabla \zeta \, \d \Omega \\[10pt]
& = & \displaystyle \int_{\Omega} \frac{1}{\det \Phi(u)} \left(
(\cof(\Phi(u)) : \nabla v)(\cof(\Phi(u)):\nabla \zeta) - 
(\cof(\Phi(u))\nabla v^T)(\cof(\Phi(u))\nabla \zeta^T)
\right) \d \Omega.
\end{array}
\end{equation*}\end{small}\noindent This symmetric form shows that the operator $v\mapsto (\sigma_N(\nabla u).\nabla v)n$ is symmetric, and in particular we have
\begin{equation*}
\displaystyle\int_{\Gamma_N} \zeta \cdot (\sigma_N(\nabla u).\nabla v)n \, \d \Gamma_N =
\displaystyle\int_{\Gamma_N} (\sigma_N(\nabla u).\nabla \zeta)n \cdot  v\, \d \Gamma_N.
\end{equation*}
When dealing with the adjoint system (from section~\ref{sec-well3}) we will still use~$\sigma_N(\nabla u)^{\ast}.\nabla \zeta$, for the sake of consistency.

\subsubsection{Assumptions on the strain energy} \label{subsec-assx}
First, we define what we call an {\it admissible} operator for a second-order linear parabolic system.

\begin{definition} \label{def-ass-op}
Introduce the following Hilbert spaces
\begin{equation*}
\mathcal{V}_0(\Omega):= \displaystyle 
\left\{v\in \HH^1(\Omega) \mid v_{|\Gamma_D} = 0 \right\}, \quad
\mathcal{V}_0(\Gamma_N) := \HH^{1/2}.
\end{equation*}
Consider the following abstract system
\begin{equation}
\begin{array} {rcl}
\ddot{u} - \kappa \Delta \dot{u} -\divg(B.\nabla u) = f & & 
\text{in $\Omega \times (0,T)$},\\
\kappa \displaystyle \frac{\p \dot{u}}{\p n} + (B.\nabla u) n = g & & 
\text{on $\Gamma_N \times (0,T)$},\\
u = 0 & & \text{on $\Gamma_D \times (0,T)$},\\
u(\cdot,0) = u_0, \quad \dot{u}(\cdot,0) = \udotini & &
\text{in $\Omega$},\\
\end{array} \label{sys-beta}
\end{equation}
Given~$T>0$, we say that the operator $B$ is {\it admissible} if for all
\begin{equation*}
f \in \L^2(0,T;\mathcal{V}_0(\Omega)'),  \quad
g \in \L^2(0,T;\mathcal{V}_0(\Gamma_N)'), \quad
u_0 \in \LL^2(\Omega), \quad \dot{u}_0 \in \LL^2(\Omega)
\end{equation*}
there exists a unique solution~$\dot{u}$ to system~\eqref{sys-beta} such that $
\dot{u} \in \L^2(0,T;\mathcal{V}_0(\Omega)) \cap \H^1(0,T;\mathcal{V}_0(\Omega)')$.
\end{definition}
The functional framework in the definition above corresponds to the standard notion of weak solutions in Hilbert spaces for a second-order parabolic system.\\

We summarize the set of general assumptions we make on the strain energy~$\mathcal{W}$, and that are needed for the analysis. Using the notation introduced previously from~$\mathcal{W}$, these assumptions are listed below:

\begin{description}
\item[$\mathbf{A1}$] The Nemytskii operator $\mathcal{W}: \WWW^{1,p}(\Omega) \ni E \mapsto \mathcal{W} (E) \in \R$ is twice continuously Fr\'echet-differentiable.
	
\item[$\mathbf{A2}$] For all matrix $E \in \R^{d\times d}$ the tensor $\check{\Sigma}(E)$ defines a symmetric matrix field.

\item[$\mathbf{A3}$] The operator~$\sigma_L(0)$ is {\it admissible} in the sense of Definition~\ref{def-ass-op}.
\end{description}

Assumptions~$\mathrm{A1}-\mathrm{A2}$ are quite natural. About Assumption~$\mathbf{A3}$, we have from~\eqref{tensor-lin1} the following expression:
\begin{equation*}
\sigma_L(0).\nabla v = \nabla v \Sigma(0) + \displaystyle
\frac{\p^2\mathcal{W}}{\p E^2}(0).(E'(0).v) = 
\nabla v \Sigma(0) + \displaystyle \frac{1}{2}
\frac{\p^2\mathcal{W}}{\p E^2}(0).(\nabla v + \nabla v^T).
\end{equation*}
In section~\ref{sec-app-energies} we show that well-known examples of strain energies from the literature fulfill these assumptions.

\subsubsection{On the control operator} \label{sec-control-operator}

The control operator, denoted as follows 
\begin{equation*}
f: \mathcal{X}_{p,T}(\omega)\ni \xi \mapsto f(\xi) \in \mathcal{F}_{p,T}(\Omega),
\end{equation*}
is distributed on a subdomain~$\omega \subset\subset \Omega$, as it appears in equation~\eqref{sysmain1}. We assume that $f$ is Fr\'echet-differentiable on~$\mathcal{X}_{p,T}(\omega)$, with values in~$\mathcal{F}_{p,T}(\Omega)$. We refer to section~\ref{sec-app-control} for modeling related comments on the control operator. In particular, $f$ may possibly be linear.

\textcolor{black}{
\subsection{Derivation of system~\eqref{sysmain}} \label{sec-deriv-Lag}
For the sake of simplicity (but without loss of generality) we can assume that the density of the material in the reference configuration is constant equal to~$1$. The kinetic and potential energies of the system are given respectively by
\begin{equation*}
\mathcal{E}_c(\dot{u}) = \displaystyle \int_{\Omega} \frac{1}{2} |\dot{u}|^2\d \Omega, 
\quad
\mathcal{E}_p(u) = \int_{\Omega} \mathcal{W}(E(u)) \d \Omega 
+ \int_{\Omega} f\cdot u \, \d \Omega
+ \int_{\Gamma_N} g \cdot u \, \d \Gamma_N.
\end{equation*}
Further, imposing that the global volume of the domain remains constant in time, we introduce the constraint
\begin{equation*}
\int_{\Omega(t)} 1 \, \d \Omega(t) = \int_{\Omega(0)} 1 \, \d \Omega(0) 
\Leftrightarrow
\int_{\Omega} \Phi(u(\cdot,t)) \d \Omega = 
\int_{\Omega} \Phi(u(\cdot,0)) \d \Omega.
\end{equation*}
We associate to this constraint a Lagrange multiplier denoted by~$\mathfrak{p}$. Consider the following {\it action}
\begin{equation*}
\mathcal{A}[T](u,\dot{u},\mathfrak{p}) = \int_0^T \mathscr{L}(u,\dot{u},\mathfrak{p}) \d t
\end{equation*}
for vector fields $u$ such that $u_{|\Gamma_D} = 0$, where
\begin{equation*}
\begin{array} {rcl}
\mathscr{L}(u,\dot{u},\mathfrak{p}) & = & \displaystyle
\mathcal{E}_c(\dot{u}) - \mathcal{E}_p(u)
 + \mathfrak{p}\left( 
\int_{\Omega} \det(\Phi(u)) \d \Omega - \int_{\Omega} \det(\Phi(u_0)) \d \Omega
\right).
\end{array}
\end{equation*}
The principle of least action applied to$\mathcal{A}[T](u,\dot{u},\mathfrak{p})$ yields
\begin{equation*}
\displaystyle 
\frac{\d}{\d t} \frac{\p \mathscr{L}}{\p \dot{u}} - 
\frac{\p \mathscr{L}}{\p u} = 0,
\end{equation*}
and hence after integration by parts, and after using~\eqref{myPiola}, the following system of equations:
\begin{equation*}
\begin{array}{rcl}
\displaystyle \ddot{u}  -\divg \left(\sigma(\nabla u) \right) = f & & \text{in } \Omega \times (0,T), \\
\displaystyle \sigma(\nabla u)n 
+ \mathfrak{p}\, \cof (\I +\nabla u)n = g & & 
\text{on } \Gamma_N \times (0,T),  \\
\displaystyle\int_{\Omega} \det(\I+\nabla u)\, \d \Omega = 
\int_{\Omega} \det(\I+\nabla u_0)\, \d \Omega & & \text{in } (0,T),\\
u=0 & & \text{on } \Gamma_D \times (0,T),\\
\displaystyle u(\cdot,0) = u_0, \quad \dot{u}(\cdot,0) = \udotini & & \text{in } \Omega.
\end{array}
\end{equation*}
Adding the damping term $\kappa \Delta \dot{u}$ in $\Omega$, corresponding to $\displaystyle \kappa \frac{\p \dot{u}}{\p n}$ on~$\Gamma_N$, yields system~\eqref{sysmain}.
}

\section{Wellposedness results} \label{sec-well}
The goal of this section is to establish existence of solutions for system~\eqref{sysmain}, and also for its linearized version around a non-trivial state, which will be used in section~\ref{sec-optcond}. We first show in section~\ref{sec-well1} that the linearized system around~$0$ is well-posed in the context of the $L^p$-maximal regularity. Via the inverse function theorem, we deduce in section~\ref{sec-well2} that the same property holds for the state system~\eqref{sysmain} under smallness assumptions on the data, and also for the non-autonomous linear system, namely system~\eqref{sysmain} linearized around a non-trivial state~$(u,\mathfrak{p})$. We rely on these results for studying in section~\ref{sec-well3} the adjoint system.

\subsection{$L^p$-maximal regularity for the linear autonomous system} \label{sec-well1}
System~\eqref{sysmain} linearized around~$(u,\mathfrak{p}) = (0,0)$ writes formally as follows:
\begin{equation} \label{siszero}
\begin{array} {rcl}
\ddot{u} - \kappa \dot{u} -\divg(\sigma_L(0).\nabla u) = f & & 
\text{in $\Omega \times (0,T)$}, \\
\kappa \displaystyle \frac{\p \dot{u}}{\p n} + 
(\sigma_L(0).\nabla u)n + \mathfrak{p}\, n = g & & 
\text{on $\Gamma_N \times (0,T)$},\\
\displaystyle \int_{\Gamma_N} u\cdot n\, \d \Gamma_N = h & & 
\text{on $(0,T)$}, \\
u = 0 & & \text{on $\Gamma_D \times (0,T)$}, \\
u(\cdot,0) = u_0, \quad \dot{u}(\cdot,0) = \udotini & &
\text{in $\Omega$.}
\end{array}
\end{equation}
The goal of this subsection is to provide wellposedness for system~\eqref{siszero}, which is stated in Corollary~\ref{coro-well-auto-NH}.
Let us first omit the constraint~$\displaystyle \int_{\Gamma_N} u\cdot n\, \d \Gamma_N = h$, and the associated pressure~$\mathfrak{p}$. We state the following result dealing with a general second-order linear parabolic system with non-homogeneous right-hand-sides and initial value conditions.

\begin{proposition} \label{prop-well-auto0}
Let be $T \in (0,\infty)$. Assume that $
f\in \mathcal{F}_{p,T}(\Omega), \quad
g\in \mathcal{G}_{p,T}(\Gamma_N), \quad
(u_0,\udotini) \in \mathcal{U}^{(0,1)}_p(\Omega)$ with the compatibility condition~$\displaystyle \kappa\frac{\p \udotini}{\p n} + (\sigma_L(0).\nabla u_0)n = g(\cdot,0)$ on~$\Gamma_N$. Then the following system
\begin{equation}
\begin{array} {rcl}
\ddot{u} - \kappa \Delta \dot{u} -\divg(\sigma_L(0).\nabla u) = f & & 
\text{in $\Omega \times (0,T)$},\\
\kappa \displaystyle \frac{\p \dot{u}}{\p n} + (\sigma_L(0).\nabla u) n = g & & 
\text{on $\Gamma_N \times (0,T)$},\\
u = 0 & & \text{on $\Gamma_D \times (0,T)$},\\
u(\cdot,0) = u_0, \quad \dot{u}(\cdot,0) = \udotini & &
\text{in $\Omega$},\\
\end{array}
\label{sys-begin}
\end{equation}
admits a unique solution~$u\in \mathcal{U}_{p,T}(\Omega)$. Further, the following estimate holds
\begin{equation*}
\|u\|_{\mathcal{U}_{p,T}(\Omega)} \leq C(T)\left(
\|f\|_{\mathcal{F}_{p,T}(\Omega)} + 
\|g\|_{\mathcal{G}_{p,T}(\Gamma_N)} +
\|(u_0,\udotini)\|_{\mathcal{U}^{(0,1)}_p(\Omega)}
\right),
\end{equation*}
where the constant~$C(T)>0$ is non-decreasing with respect to~$T$.
\end{proposition}

\begin{proof}
Such a result falls into the frame of the so-called $L^p$-maximal parabolic regularity~\cite{AMS2003}. Given Assumption~$\mathbf{A3}$, the result stated in~\cite[Theorem~6.1]{Arendt2007} addresses the question of existence of solutions for this type of second-order autonomous equations with homogeneous Dirichlet conditions on~$\p \Omega$. System~\eqref{sys-begin} can be rewritten in terms of~$\dot{u}$ as main unknown variable, and thus becomes a first-order parabolic system. Then the results provided in~\cite{Pruess2002} when considering mixed boundary conditions apply, and thus we obtain the $L^p$-maximal regularity property for system~\eqref{sys-begin}.
\end{proof}

Now we establish the same type of results for system~\eqref{siszero}, when its constraint is homogeneous:
\begin{proposition} \label{prop-well-auto}
Let be $T \in (0,\infty)$, and assume the hypotheses of Proposition~\ref{prop-well-auto0}, with additionally
\begin{equation*}
\displaystyle \int_{\Gamma_N} u_0 \cdot n \, \d \Gamma_N = 0.
\end{equation*}
Then the following system
\begin{equation} \label{syslih}
\begin{array} {rcl}
\ddot{u} - \kappa \Delta \dot{u} -\divg(\sigma_L(0).\nabla u) = f & & 
\text{in $\Omega \times (0,T)$}, \\
\kappa \displaystyle \frac{\p \dot{u}}{\p n} + (\sigma_L(0).\nabla u) n
+\mathfrak{p}\, n = g & &  \text{on $\Gamma_N \times (0,T)$}, \\
\displaystyle
\int_{\Gamma_N} u\cdot n\, \d \Gamma_N = 0 & & \text{in $(0,T)$} , \\
u = 0 & & \text{on $\Gamma_D \times (0,T)$},\\
u(\cdot,0) = u_0, \quad \dot{u}(\cdot,0) = \udotini & &
\text{in $\Omega$},
\end{array}
\end{equation}
admits a unique solution~$(u,\mathfrak{p})\in \mathcal{U}_{p,T}(\Omega)\times \mathcal{P}_{p,T}$. Further, the following estimate holds
\begin{equation}
\|u\|_{\mathcal{U}_{p,T}(\Omega)} + 
\|\mathfrak{p}\|_{\mathcal{P}_{p,T}}\leq C(T)\left(
\|f\|_{\mathcal{F}_{p,T}(\Omega)} + 
\|g\|_{\mathcal{G}_{p,T}(\Gamma_N)} +
\|(u_0,\udotini)\|_{\mathcal{U}^{(0,1)}_p(\Omega)}
\right), \label{est-prop32}
\end{equation}
where the constant~$C(T)>0$ is non-decreasing with respect to~$T$.
\end{proposition}

\noindent The proof of Proposition~\ref{prop-well-auto} is given in the section~\ref{sec-app-tek1}. We deduce the same result when the constraint of system~\eqref{syslih} is non-homogeneous. More precisely, we consider constraint of system~\eqref{siszero} with $h\in \mathcal{H}_{p,T} = \W^{1-1/(2p),p}(0,T;\R)$.

\begin{corollary} \label{coro-well-auto-NH}
Let be $T \in (0,\infty)$, and assume $
f\in \mathcal{F}_{p,T}(\Omega), \quad
g\in \mathcal{G}_{p,T}(\Gamma_N), \quad
h\in \mathcal{H}_{p,T}, \quad 
(u_0,\udotini) \in \mathcal{U}^{(0,1)}_p(\Omega)$ satisfying the compatibility conditions
\begin{equation*}
\begin{array}{l}
\displaystyle g(\cdot,0) = \kappa\frac{\p \udotini}{\p n} +(\sigma_L(0).\nabla u_0)n 
\ \text{on } \Gamma_N, \\
h(0) = \displaystyle \int_{\Gamma_N} u_0 \cdot n\, \d \Gamma_N
\quad \text{and} \quad 
\dot{h}(0) = \displaystyle \int_{\Gamma_N} \dot{u}_0 \cdot n\, \d \Gamma_N.
\end{array}
\end{equation*} 
Then there exists a unique solution~$(u,\mathfrak{p})$ to the following system
\begin{equation} \label{siszeroNH}
\begin{array} {rcl}
\ddot{u} - \kappa \Delta \dot{u} -\divg(\sigma_L(0).\nabla u) = f & & 
\text{in $\Omega \times (0,T)$},\\
\kappa \displaystyle \frac{\p \dot{u}}{\p n} + (\sigma_L(0).\nabla u) n
+\mathfrak{p}\, n = g & &  \text{on $\Gamma_N \times (0,T)$},\\
\displaystyle
\int_{\Gamma_N} u\cdot n\, \d \Gamma_N = h & & \text{in $(0,T)$} , 
\\
u = 0 & & \text{on $\Gamma_D \times (0,T)$},\\
u(\cdot,0) = u_0, \quad \dot{u}(\cdot,0) = \udotini & &
\text{in $\Omega$},
\end{array}
\end{equation}
and it satisfies the following estimate
\begin{equation*}
\|u\|_{\mathcal{U}_{p,T}(\Omega)} + 
\|\mathfrak{p}\|_{\mathcal{P}_{p,T}}\leq C(T)\left(
\|f\|_{\mathcal{F}_{p,T}(\Omega)} + 
\|g\|_{\mathcal{G}_{p,T}(\Gamma_N)} +
\|h\|_{\mathcal{H}_{p,T}} +
\|(u_0,\udotini)\|_{\mathcal{U}^{(0,1)}_p(\Omega)}
\right),
\end{equation*}
where the constant~$C(T)>0$ is non-decreasing with respect to~$T$.
\end{corollary}
\noindent The technical proof of Corollary~\ref{coro-well-auto-NH} is given in section~\ref{sec-app-tek2}.

\subsection{Local existence result for the state system} \label{sec-well2}

Define the mapping
\begin{equation*}
\begin{array} {rrcc}
\mathcal{K}:  & 
\mathcal{U}_{p,T}(\Omega) \times \mathcal{P}_{p,T}
& \rightarrow & 
\mathcal{F}_{p,T}(\Omega) \times \mathcal{G}_{p,T}(\Gamma_N) \times 
\mathcal{H}_{p,T} \times \mathcal{U}_p^{(0,1)}(\Omega)
\\[5pt]
&  (u,\mathfrak{p}) & \mapsto & 
\left( \begin{array}{c}
\ddot{u} - \kappa \Delta\dot{u} - \divg(\sigma(\nabla u))\\
\kappa \displaystyle \frac{\p \dot{u}}{\p n} + \sigma(\nabla u)n 
+\mathfrak{p}\, \cof(\Phi(u))n \\
\displaystyle \int_{\Omega} \det(\Phi(u)) \d \Omega\\
	(u(\cdot,0), \dot{u}(\cdot,0))
	\end{array} \right).
\end{array}
\end{equation*}
We have $\mathcal{K}(0,0) = (-\divg(\sigma(0)),\sigma(0)n,|\Omega|,0)^T$. From Corollary~\ref{coro-well-auto-NH}, the differential of mapping~$\mathcal{K}$ at~$(u,\mathfrak{p}) = (0,0)$ is an isomorphism. Therefore, from the inverse function theorem, system~\eqref{sysmain} is locally wellposed. More precisely, we state the following result:

\begin{proposition} \label{prop-well-sysmain}
Let be $T\in(0,\infty)$. There exists~$\eta>0$ such that if
\begin{equation*}
\|f+ \divg(\sigma(0))\|_{\mathcal{F}_{p,T}(\Omega)} + 
\|g-\sigma(0)n\|_{\mathcal{G}_{p,T}(\Gamma_N)} +
\|(u_0,\udotini)\|_{\mathcal{U}^{(0,1)}_p(\Omega)} \leq \eta
\end{equation*}
with the compatibility condition $\displaystyle \kappa \frac{\p \dot{u}_0}{\p n} + \sigma(\nabla u_0)n = g(\cdot,0)$ on~$\Gamma_N$, then system~\eqref{sysmain} admits a unique solution~$(u,\mathfrak{p}) \in \mathcal{U}_{p,T}(\Omega) \times \mathcal{P}_{p,T}$.
\end{proposition}

\begin{proof}
Note that the smallness assumption on~$u_0 \in \WW^{2,p}(\Omega)$ implies\begin{small}
\begin{equation*}
\begin{array} {rcl}
\displaystyle
\left||\Omega|  - \int_{\Omega}\det (\Phi(u_0))\d \Omega \right| & = &
\displaystyle \left| \int_{\Omega}\det (\Phi(0))\d \Omega - \int_{\Omega}\det (\Phi(u_0))\d \Omega \right|\\[10pt]
& \leq & 
\|\nabla u_0\|_{\WWW^{1,p}(\Omega)} 
\displaystyle \sup_{\alpha \in [0,1]}
\int_{\Omega} \| \cof(\I + \alpha \nabla u)\|\d \Omega\\
& \leq & \displaystyle C\eta\sum_{i=0}^{d-1} \eta^i,
\end{array}
\end{equation*}\end{small}\noindent where we have used the mean-value theorem in the algebra~$\WWW^{1,p}(\Omega)$. Therefore, provided that~$\eta>0$ is chosen small enough, the assumptions of the inverse function theorem are satisfied, which yields the result.
\end{proof}

Further, the differential of mapping~$\mathcal{K}$ is also locally an isomorphism, that means that system~\eqref{sysmain} linearized around some state~$(u,\mathfrak{p})$ is well-posed, provided smallness assumptions on~$(u,\mathfrak{p})$. This is achieved by assuming the data small enough, in virtue of Proposition~\ref{prop-well-sysmain}. This non-autonomous linear system writes formally as follows, where~$(v,\mathfrak{q})$ denotes its unknown:
\begin{equation} \label{syslin-nonauto}
\begin{array} {rcl}
\ddot{v} - \kappa \Delta\dot{v} -\divg(\sigma_L(\nabla u).\nabla v) = f & & 
\text{in $\Omega \times (0,T)$},\\
\kappa \displaystyle \frac{\p \dot{v}}{\p n} + 
\Big(\big(\sigma_L(\nabla u)+\mathfrak{p}\, \sigma_N(\nabla u)\big).\nabla v \Big) n
+ \mathfrak{q}\, \cof(\Phi(u))n = g & & 
\text{on $\Gamma_N \times (0,T)$},\\
\displaystyle \int_{\Gamma_N} v\cdot \cof(\Phi(u))n\, \d \Gamma_N = 0 & & 
\text{on $(0,T)$}, \\
v = 0 & & \text{on $\Gamma_D \times (0,T)$},\\
v(\cdot,0) = 0,
\quad \dot{v}(\cdot,0) = 0, & & 
\text{in $\Omega$}.
\end{array}
\end{equation}
Recall that~$\sigma_N$ is introduced in~\eqref{tensor-lin2}. Note that in the writing of system~\eqref{syslin-nonauto}, only~$u$ appears, not~$\dot{u}$, as the terms involving~$\dot{u}$ are linear in system~\eqref{sysmain}. We state the following result:

\begin{proposition} \label{prop-well-syslin}
Let be $T\in(0,\infty)$, and assume that~$(u,\mathfrak{p})$ is small enough in~$\mathcal{U}_{p,T}(\Omega) \times \mathcal{P}_{p,T}$. Then for all $f\in \mathcal{F}_{p,T}(\Omega)$, $g \in \mathcal{G}_{p,T}(\Gamma_N)$
satisfying the compatibility condition~$g(\cdot,0) = 0$,
system~\eqref{syslin-nonauto} admits a unique solution~$(v,\mathfrak{q}) \in \mathcal{U}_{p,T}(\Omega) \times \mathcal{P}_{p,T}$, and it satisfies
\begin{equation}
\|v\|_{\mathcal{U}_{p,T}(\Omega)} +
\|\mathfrak{q} \|_{\mathcal{P}_{p,T}} \leq
C(u,\mathfrak{p}) \left(
\|f\|_{\mathcal{F}_{p,T}(\Omega)} + 
\|g\|_{\mathcal{G}_{p,T}(\Gamma_N)} 
\right).
\label{est-prop-nonauto}
\end{equation}
The constant $C(u,\mathfrak{p})>0$ is independent of~$(v,\mathfrak{q})$.
\end{proposition}

\subsection{The adjoint system} \label{sec-well3}
Let us first rewrite the second-order linear system~\eqref{syslin-nonauto} in the form of a first-order parabolic system, by setting~$(y_0,y_1) = (u,\dot{u})$ and~$(z_0,z_1) = (v,\dot{v})$. More generally we consider the following system:
\begin{equation} \label{sysx}
\begin{array} {rcl}
\dot{z}_0 - z_1 = f_0
 & & \text{in $\Omega \times(0,T)$}, \\
\dot{z}_1 - \kappa \Delta z_1
-  \divg (\sigma_L(\nabla y_0).\nabla z_0) = 
 f_1 & & \text{in $\Omega \times(0,T)$},  \\
\kappa \displaystyle \frac{\p z_1}{\p n} +  
\Big((\sigma_L + \mathfrak{p}\, \sigma_N)(\nabla y_0).\nabla z_0\Big)n +   \mathfrak{q}\, \cof\left(\Phi(y_0) \right) n = g & & \text{on $\Gamma_N\times (0,T)$}, \\
\displaystyle \int_{\Gamma_N} z_0 \cdot \cof(\Phi(y_0))n\, \d \Gamma_N
= 0 & & \text{in } (0,T),\\
z_1 = 0 & & \text{on $\Gamma_D\times (0,T)$},  \\
z_0(\cdot, 0) = 0, \displaystyle \quad z_1(\cdot, 0) = 0 & & \text{in $\Omega$}.
\end{array}
\end{equation}
Note that only~$y_0 = u$ appears in the writing of system~\eqref{sysx}, not~$y_1$. Using Proposition~\ref{prop-well-syslin}, we state the following proposition:

\begin{proposition} \label{prop-well-syslin-init}
Let be $T\in(0,\infty)$, and assume that~$(y_0,\mathfrak{p})$ is small enough in~$\mathcal{U}_{p,T}(\Omega) \times \mathcal{P}_{p,T}$. Then for all $
f_0 \in \dot{\mathcal{U}}_{p,T}(\Omega), \quad
f_1\in \mathcal{F}_{p,T}(\Omega), \quad
g \in \mathcal{G}_{p,T}(\Gamma_N)$ satisfying the compatibility condition~$g(\cdot,0) = 0$, system~\eqref{syslin-nonauto} admits a unique solution~$(z_0,z_1,\mathfrak{q}) \in \mathcal{U}_{p,T}(\Omega) \times \dot{\mathcal{U}}_{p,T}(\Omega) \times \mathcal{P}_{p,T}$, and it satisfies
\begin{equation*}
\|z_0\|_{\mathcal{U}_{p,T}(\Omega)} +
\|z_1\|_{\dot{\mathcal{U}}_{p,T}(\Omega)} +
\|\mathfrak{q} \|_{\mathcal{P}_{p,T}} \leq
C(y_0,\mathfrak{p}) \left(
\|f_0\|_{\dot{\mathcal{U}}_{p,T}(\Omega)} + 
\|f_1\|_{\mathcal{F}_{p,T}(\Omega)} +
\|g\|_{\mathcal{G}_{p,T}(\Gamma_N)}
\right).
\end{equation*}
The constant $C(y_0,\mathfrak{p})>0$ is independent of~$(z_0,z_1,\mathfrak{q})$.
\end{proposition}

\begin{proof}
The first equation of system~\eqref{sysx} derived in time, and combined with the second equation, shows that variable~$z_0$ satisfies system~\eqref{syslin-nonauto}, with respectively $z_0$ in the role of~$v$, $f_1+\dot{f}_0 -\kappa \Delta f_0 \in \mathcal{F}_{p,T}(\Omega)$ in the role of~$f$, and $g +\kappa \displaystyle \frac{\p f_0}{\p n} \in \mathcal{G}_{p,T}(\Gamma_N)$ in the role of~$g$. Proposition~\ref{prop-well-syslin} states the existence and uniqueness of~$z_0 \in \mathcal{U}_{p,T}(\Omega)$, which satisfies\begin{small}
\begin{eqnarray}
\|z_0\|_{\mathcal{U}_{p,T}(\Omega)} + \|\mathfrak{q}\|_{\mathcal{P}_{p,T}} & \leq & 
C(y_0,\mathfrak{p}) \left(
\|f_1\|_{\mathcal{F}_{p,T}(\Omega)} +
\|\dot{f}_0\|_{\mathcal{F}_{p,T}(\Omega)} + 
\|\Delta f_0\|_{\mathcal{F}_{p,T}(\Omega)}\right. \nonumber \\ & & \left.
\displaystyle +
\|g\|_{\mathcal{G}_{p,T}(\Gamma_N)}+ 
\left\|\frac{\p f_0}{\p n}\right\|_{\mathcal{G}_{p,T}(\Gamma_N)}
\right) \nonumber \\
& \leq &
CC(y_0,\mathfrak{p}) \left(
\|f_0\|_{\dot{\mathcal{U}}_{p,T}(\Omega)}+
\|f_1\|_{\mathcal{F}_{p,T}(\Omega)} +
\|g\|_{\mathcal{G}_{p,T}(\Gamma_N)} 
\right), \label{est-z0}
\end{eqnarray}\end{small}\noindent where the constant~$C(y_0,\mathfrak{p})$ is the one of estimate~\eqref{est-prop-nonauto}. We still denote $C(y_0,\mathfrak{p})$ constants of type $CC(y_0,\mathfrak{p})$. Further, the first equation of system~\eqref{sysx}, namely~$z_1 = \dot{z}_0 - f_0$, yields
\begin{equation*}
\|z_1\|_{\dot{\mathcal{U}}_{p,T}(\Omega)} \leq 
\|z_0\|_{\mathcal{U}_{p,T}(\Omega)} + 
\|f_0\|_{\dot{\mathcal{U}}_{p,T}(\Omega)}.
\end{equation*}
Combined with~\eqref{est-z0}, we deduce the announced estimate and complete the proof.
\end{proof}

We stress that solutions~$((z_0,z_1),\mathfrak{q})$ of system~\eqref{sysx} are continuous on~$[0,T]$ -- with values in $\mathcal{U}_p^{(0,1)}(\Omega) \times \R$. Recall that our goal is to address Problem~\eqref{mainpb}, involving the functionals~$c: \mathcal{U}_{p,T}(\Omega) \times \dot{\mathcal{U}}_{p,T}(\Omega) \times \mathcal{X}_{p,T}(\omega)\rightarrow \R$, $\phi^{(1)}: \mathcal{U}^{(0,1)}_{p}(\Omega) \rightarrow \R$ and $\phi^{(2)}: \mathcal{U}^{(0,1)}_{p}(\Omega) \rightarrow \R$. Given $\tau \in (0,T)$, $\xi \in \mathcal{X}_{p,T}(\omega)$ and $(y_0,y_1, \mathfrak{p}) \in \mathcal{U}_{p,T}(\Omega) \times \dot{\mathcal{U}}_{p,T}(\Omega)\times \mathcal{P}_{p,T}$, we introduce the adjoint system, namely\begin{small}
\begin{equation} \label{sysadjoint-init}
\begin{array} {rcl}
-\dot{\zeta}_0 -  \divg(\sigma_L(\nabla y_0)^{\ast}.\nabla \zeta_1) = 
-c'_{y_0}(y_0, y_1, \xi) 
& & \text{in $\Omega \times\left((0,\tau) \cup (\tau,T)\right)$}, \\
-\dot{\zeta}_1 -  \zeta_0 -\kappa  \Delta \zeta_1 = 
-c'_{y_1}(y_0,y_1, \xi) 
& & \text{in $\Omega \times(0,T)$}, \\ 
\Big((\sigma_L + \mathfrak{p}\, \sigma_N)(\nabla y_0)^{\ast}.\nabla \zeta_1\Big)n  
+ \pi\, \cof(\Phi(y_0)) n = 0 & & 
\text{on $\Gamma_N\times (0,T)$},  \\
\kappa \displaystyle \frac{\p \zeta_1}{\p n}  =0
 & &  \text{on $\Gamma_N\times (0,T)$},\\[5pt]
\zeta_1 = 0 & & \text{on $\Gamma_D\times(0,T)$},  \\
\displaystyle \left\langle \zeta_1 \, ; \cof(\Phi(y_0))n\right\rangle_{\WW^{1/(p'),p}(\Gamma_N)',\WW^{1/(p'),p}(\Gamma_N)} 
= 0 & & \text{in } (0,T)\\
\left[ \zeta_0 \right]_{\tau} = \phi^{(1)'}_{y_0}(y_0,y_1) (\tau)
\quad \text{and} \quad 
\left[ \zeta_1 \right]_{\tau} = \phi^{(1)'}_{y_1}(y_0,y_1) (\tau)
& & \text{in $\Omega$},
\\
\zeta_0(\cdot, T) = -\phi^{(2)'}_{y_0}(y_0,y_1)(T), \displaystyle 
\quad \zeta_1(\cdot, T) = -\phi^{(2)'}_{y_1}(y_0,y_1)(T) & & \text{in $\Omega$}.
\end{array}
\end{equation}\end{small}\noindent We have introduced the notation~$\left[ \zeta \right]_{\tau} := 
\displaystyle\lim_{t\rightarrow \tau^+} \zeta(t) - \lim_{t\rightarrow \tau^-} \zeta(t)$ which describes the jump of a variable~$\zeta$ at time~$t=\tau$. We define solutions of the adjoint system~\eqref{sysadjoint-init} by~{\it transposition}.

\begin{definition} \label{def-trans-init}
Let be $\tau \in (0,T)$ and $(y_0, y_1,\mathfrak{p},\xi) \in \mathcal{U}_{p,T}(\Omega) \times \dot{\mathcal{U}}_{p,T}(\Omega)\times \mathcal{P}_{p,T}\times \mathcal{X}_{p,T}(\omega)$. We say that~$(\zeta_0,\zeta_1,\pi)$ is a solution of~\eqref{sysadjoint-init} associated with $(y_0,y_1,\mathfrak{p})$, if for all $(f_0, f_1,g) \in \dot{\mathcal{U}}_{p,T}(\Omega) \times \mathcal{F}_{p,T}(\Omega) \times \mathcal{G}_{p,T}(\Gamma_N)$ we have
\begin{equation}
\begin{array} {l}
 \displaystyle \left\langle \zeta_0 ;f_0 \right\rangle_{\dot{\mathcal{U}}_{p,T}(\Omega)',\dot{\mathcal{U}}_{p,T}(\Omega)}
 +\left\langle \zeta_1 ; f_1 \right\rangle_{\mathcal{F}_{p,T}(\Omega)',\mathcal{F}_{p,T}(\Omega)}
+ \left\langle \zeta_1 ; g \right\rangle_{\mathcal{G}_{p,T}(\Gamma_N)',\mathcal{G}_{p,T}(\Gamma_N)} 
\\
 \displaystyle =
-\left\langle c'_{y_0}(y_0,y_1, \xi)\, ; z_0 \right\rangle_{\mathcal{U}_{p,T}(\Omega)',\mathcal{U}_{p,T}(\Omega)}
- \left\langle  c'_{y_1}(y_0,y_1, \xi)\, ; z_1 \right\rangle_{\dot{\mathcal{U}}_{p,T}(\Omega)',\dot{\mathcal{U}}_{p,T}(\Omega)}
 \\
 - \left\langle(\phi^{(1)'}_{y_0}(y_0,y_1), \phi^{(1)'}_{y_1}(y_0,y_1) )\, ; (z_0(\cdot,\tau),z_1(\cdot,\tau))
\right\rangle_{\mathcal{U}_p^{(0,1)}(\Omega)',\mathcal{U}_p^{(0,1)}(\Omega)}
 \\
 -\left\langle 
(\phi^{(2)'}_{y_0}(y_0,y_1), \phi^{(2)'}_{y_1}(y_0,y_1) )\, ; (z_0(\cdot,T),z_1(\cdot,T))
\right\rangle_{\mathcal{U}_p^{(0,1)}(\Omega)',\mathcal{U}_p^{(0,1)}(\Omega)} ,
\label{id-def-trans-init}
\end{array}
\end{equation}
where~$(z_0,z_1,\mathfrak{q})$ is the solution of system~\eqref{sysx} with $(y_0,y_1,\mathfrak{p})$ and $(f_0,f_1,g)$ as data.
\end{definition}

\begin{remark}
Solutions~$(\zeta_0,\zeta_1)$ in the sense of Definition~\ref{def-trans-init} lie in $\dot{\mathcal{U}}_{p,T}(\Omega)' \times \mathcal{F}_{p,T}(\Omega)'$, and therefore satisfy
\begin{equation*} 
\begin{array}{l}
\zeta_0 \in \L^{p'}(0,T;\WW^{2,p}(\Omega)'), \quad
\zeta_1 \in \L^{p'}(0,T;\LL^{p'}(\Omega)), \quad
{\zeta_1}_{|\Gamma_N} \in \L^{p'}(0,T; \WW^{1/(p'),p}(\Gamma_N)') \\
\nabla \zeta_1 \in \L^{p'}(0,T; \WWW^{1,p}(\Omega)'), \quad
\displaystyle \frac{\p \zeta_1}{\p n} \in \L^{p'}(0,T; \WW^{2-1/p,p}(\Gamma_N)').
\end{array}
\end{equation*}
It is unnecessary to comment on the regularity of the variable~$\pi$, playing the role of Lagrange multiplier for the constraint imposed on~$\zeta_1$ (sixth equation of~system~\eqref{sysadjoint-init}).
\end{remark}
We prove the existence and uniqueness of a {\it very weak} solution for system~\eqref{sysadjoint-init}, in the sense of Definition~\ref{def-trans-init}.

\begin{proposition} \label{propadj-init}
Let be $(y_0,y_1,\mathfrak{p},\xi) \in \mathcal{U}_{p,T}(\Omega) \times \dot{\mathcal{U}}_{p,T}(\Omega)\times \mathcal{P}_{p,T}\times \mathcal{X}_{p,T}(\omega) $. If~$(y_0,\mathfrak{p})$ is small enough in~$\mathcal{U}_{p,T}(\Omega) \times \mathcal{P}_{p,T}$, then system~\eqref{sysadjoint-init} admits a unique solution $(\zeta_0,\zeta_1,\pi) \in \dot{\mathcal{U}}_{p,T}(\Omega)' \times \mathcal{F}_{p,T}(\Omega)' \times \mathcal{G}_{p,T}(\Gamma_N)'$, in the sense of Definition~\ref{def-trans-init}. Moreover, there exists a constant $C(y_0, \mathfrak{p}) >0$ depending only on $(y_0,\mathfrak{p})$ such that\begin{small}
\begin{equation*}
\begin{array} {rcl}
\|(\zeta_0,\zeta_1) \|_{\dot{\mathcal{U}}_{p,T}(\Omega)' \times \mathcal{F}_{p,T}(\Omega)'}
& \leq & C(y_0, \mathfrak{p}) \left(
\| c'_{y_0}(y_0,y_1,\xi) \|_{\mathcal{U}_{p,T}(\Omega)'} +
\| c'_{y_1}(y_0,y_1,\xi) \|_{\dot{\mathcal{U}}_{p,T}(\Omega)'} \right.\\
& & \left. + \| (\phi^{(1)'}_{y_0}(y_0,y_1)(\tau), \phi^{(1)'}_{y_1}(y_0,y_1)(\tau)) \|_{\mathcal{U}^{(0,1)}_p(\Omega)'} \right. \\
&  & \left. +\| (\phi^{(2)'}_{y_0}(y_0,y_1)(T), \phi^{(2)'}_{y_1}(y_0,y_1)(T)) \|_{\mathcal{U}^{(0,1)}_p(\Omega)'}
\right).
\end{array}
\end{equation*}\end{small}\noindent In particular, $C(y_0,\mathfrak{p})$ is independent of $c$, $\phi^{(1)}$ and $\phi^{(2)}$.
\end{proposition}

\begin{proof}
Define the operator
\begin{equation*}
\Lambda(y_0,y_1,\mathfrak{p}): (f_0,f_1,g) \mapsto \big(z_0,z_1, (z_0(\cdot,\tau),z_1(\cdot,\tau)), (z_0(\cdot,T),z_1(\cdot,T))\big),
\end{equation*}
where~$(z_0,z_1)$ is the solution of system~\eqref{sysx}. From Proposition~\ref{prop-well-syslin-init}, the linear operator~$\Lambda(y_0,y_1,\mathfrak{p})$ is bounded from~$\dot{\mathcal{U}}_{p,T}(\Omega) \times \mathcal{F}_{p,T}(\Omega) \times \mathcal{G}_{p,T}(\Gamma_N)$ into~$\mathcal{U}_{p,T}(\Omega) \times \dot{\mathcal{U}}_{p,T}(\Omega) \times \mathcal{U}_p^{(0,1)}(\Omega) \times \mathcal{U}_p^{(0,1)}(\Omega)$. Therefore~$\Lambda(y_0,y_1,\mathfrak{p})^{\ast}$ is bounded from~$\mathcal{U}_{p,T}(\Omega)' \times \dot{\mathcal{U}}_{p,T}(\Omega)' \times \mathcal{U}_p^{(0,1)}(\Omega)' \times \mathcal{U}_p^{(0,1)}(\Omega)'$ into~$\dot{\mathcal{U}}_{p,T}(\Omega)' \times \mathcal{F}_{p,T}(\Omega)' \times \mathcal{G}_{p,T}(\Gamma_N)'$. Defining\begin{scriptsize}
\begin{equation*}
(\zeta_0,\zeta_1,\pi) = -\Lambda(y_0,y_1,p)^{\ast}\left(
c_{y_0}'(y_0,y_1,\xi), c_{y_1}'(y_0,y_1,\xi), \phi^{(1)'}_{y_0}(y_0,y_1),
 \phi^{(1)'}_{y_1}(y_0,y_1), \phi^{(2)'}_{y_0}(y_0,y_1),
 \phi^{(2)'}_{y_1}(y_0,y_1)\right),
\end{equation*}\end{scriptsize}\noindent we can verify that~$(\zeta_0, \zeta_1,\pi)$ is solution of system~\eqref{sysadjoint-init} in the sense of Definition~\ref{def-trans-init}. For that we use the Green formula twice, and proceed by integration by parts in $(0,\tau) \cup (\tau,T)$. Uniqueness is due to the linearity of system~\eqref{sysadjoint-init}. Indeed, if\begin{small}
\begin{equation*}
\left(c_{y_0}'(y_0,y_1,\xi), c_{y_1}'(y_0,y_1,\xi), \phi^{(1)'}_{y_0}(y_0,y_1)(\tau),
 \phi^{(1)'}_{y_1}(y_0,y_1)(\tau), \phi^{(2)'}_{y_0}(y_0,y_1)(T),
 \phi^{(2)'}_{y_1}(y_0,y_1)(T)\right) 
\end{equation*}\end{small}\noindent is equal to~$(0,0,0,0,0,0)$, then from~\eqref{id-def-trans-init} we deduce that $(\zeta_0, \zeta_1) = (0,0)$ in $\dot{\mathcal{U}}_{p,T}(\Omega)' \times  \mathcal{F}_{p,T}(\Omega)'$, which also implies in the third equation of~\eqref{sysadjoint-init} that~$\pi = 0$, completing the proof. 
\end{proof}

\section{Optimal control formulation and optimality conditions} \label{sec-optcond}
The purpose of this section is the derivation of necessary optimality conditions for the original optimal control problem~\eqref{mainpb}:
\begin{equation} \tag{$\mathcal{P}$}
\left\{ \begin{array} {l}
	\displaystyle \max_{(\xi,\tau) \in \mathcal{X}_{p,T}(\omega) \times (0,T)}
\left( J(\xi,\tau) =
\int_0^T c(u,\dot{u}, \xi)\, \d t + \phi^{(1)}(u,\dot{u})(\tau) + 
\phi^{(2)}(u,\dot{u})(T) \right)
\\[10pt]
	 \text{where $(u, \dot{u})$ satisfies~\eqref{sysmain}.} 
\end{array} \right. \label{optcontprob}
\end{equation}

\begin{remark}\label{rk-box}
In order to guarantee the existence of solutions~$(u,\dot{u})$ to system~\eqref{sysmain} independently of the set of control/parameters~$(\xi,\tau)\in \mathcal{X}_{p,T}(\omega)\times(0,T)$ that is considered in Problem~\eqref{optcontprob}, we may need to add norm constraints on the control function~$\xi$. More precisely, in virtue of Proposition~\ref{prop-well-sysmain}, we could consider in addition the following constraint
\begin{equation*}
\|f(\xi)+ \divg(\sigma(0))\|_{\mathcal{F}_{p,T}(\Omega)} + 
\|g-\sigma(0)n\|_{\mathcal{G}_{p,T}(\Gamma_N)} +
\|(u_0,\udotini)\|_{\mathcal{U}^{(0,1)}_p(\Omega)} \leq \eta,
\end{equation*}
for some~$\eta>0$ chosen small enough. We would then proceed like in~\cite{Maxmax3} for deriving the corresponding optimality conditions, incorporating a Lagrange multiplier for taking into account such a constraint. In order to not make the complexity heavier, we choose to omit this point in what follows.
\end{remark}

For example, like in the illustrations presented in section~\ref{sec-num}, one can choose $c(u,\dot{u}, \xi)= - \frac{1}{2} \| \xi\|_{\L^2(\omega)}^2$ as cost functional.
Due to a lack of smoothness of the state at a (possible) optimal time~$\tau$, we need to make a change of variable, in order to uncouple the state variable~$(u,\dot{u})$ and the time parameter~$\tau$.

\subsection{Transformation of the problem and new formulation} \label{sec-transformation}
\begin{figure}[H]
\begin{center}
\begin{tikzpicture}[line/.style={>=latex}] 
\coordinate (V1) at (0,0);
\coordinate (V2) at (3,2);
\coordinate (V3) at (4,2.5);
\coordinate (V4) at (6,4);
\draw[step=14.25pt, color=black!10] (-0.5, -0.5) grid (6.5, 4.5);
\draw[->, line] (-0.5, 0) -- node [below, at end] {$s$} (6.5, 0);
\draw[->, line] (0, -0.5) -- node [left, at end] {$t$} (0, 4.5);
\foreach \x in {0}
   \draw (\x cm,1pt) -- (\x cm,-1pt) node[anchor=north east] {$\x$};
\foreach \x in {3}
   \draw (\x cm,1pt) -- (\x cm,-1pt) node[anchor=north] {$1$};
\foreach \x in {4}
   \draw (\x cm,1pt) -- (\x cm,-1pt) node[anchor=north] {$1+\tilde{\varepsilon}$};
\foreach \x in {6}
   \draw (\x cm,1pt) -- (\x cm,-1pt) node[anchor=north] {$2$};
\foreach \y in {2}
    \draw (1pt,\y cm) -- (-1pt,\y cm) node[anchor=east] {$\tau$};
\foreach \y in {2.5}
    \draw (1pt,\y cm) -- (-1pt,\y cm) node[anchor=east] {$\tau+\varepsilon$};
\foreach \y in {4}
    \draw (1pt,\y cm) -- (-1pt,\y cm) node[anchor=east] {$T$};
\draw[-, line, color=red!70!black, thick, densely dotted] (V1) -- (V2);
\draw[-, line, color=red!70!black, thick, densely dotted] (V2) -- (V3);
\draw[-, line, color=red!70!black, thick, densely dotted] (V3) -- (V4);
\draw[-, line, color=black, thick, dashed] (0,4) -- (6,4);
\draw[-, line, color=black, thick, dashed] (0,2.5) -- (4,2.5);
\draw[-, line, color=black, thick, dashed] (0,2) -- (3,2);
\draw[-, line, color=black, thick, dashed] (3,0) -- (3,2);
\draw[-, line, color=black, thick, dashed] (4,0) -- (4,2.5);
\draw[-, line, color=black, thick, dashed] (6,0) -- (6,4);
\draw (V1) node[circle,color=green!50!black,fill,inner sep=1pt,label=above:\textcolor{green!50!black}{$ $}] {};
\draw (V2) node[circle,color=green!50!black,fill,inner sep=1pt,label=above:\textcolor{green!50!black}{$ $}] {};
\draw (V3) node[circle,color=green!50!black,fill,inner sep=1pt,label=above:\textcolor{green!50!black}{$ $}] {};
\draw (V4) node[circle,color=green!50!black,fill,inner sep=1pt,label=above:\textcolor{green!50!black}{$ $}] {};
\end{tikzpicture}
\caption{Graph of the change of variables $s \mapsto \mu(s,\tau)$.}\label{fig-graph-mu}
\end{center}
\end{figure}
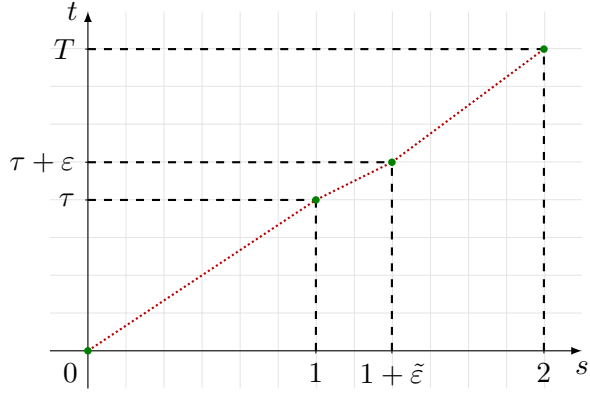

Let $\tilde{\varepsilon} \in (0,1)$ be a fix parameter chosen small enough, typically of the range of the time step when doing numerical simulations in section~\ref{sec-num}. When considering functionals~$\phi^{(1)}$ like~\eqref{sysmainpressure} for example, we introduce the change of variables $\mu :[0,2] \rightarrow [0,T]$ given as follows (see Figure~\ref{fig-graph-mu}):
\begin{equation*} 
	\mu(s,\tau) =  
	\left\{ \begin{array} {ll}
		\tau s & \text{if } s\in [0,1], \\[5pt]
		\displaystyle \frac{\varepsilon}{\tilde{\varepsilon}}(s-1) + \tau & \text{if } s\in [1,1+\tilde{\varepsilon}], \\[5pt]
		\displaystyle T- (2-s) \frac{T-(\tau+\varepsilon)}{1-\tilde{\varepsilon}} & \text{if } s\in [1+\tilde{\varepsilon}, 2].
	\end{array} \right.
\end{equation*}
This change of variable is designed such that~$\mu(\cdot,\tau)$ is bijective from $[0,2]$ to $[0,T]$, and
\begin{equation*}
\mu(0,\tau) = 0, \quad 
\mu(2,\tau) = T, \quad
\mu(1,\tau)= \tau, \quad
\mu(1+\tilde{\varepsilon},\tau) = \tau+ \varepsilon.
\end{equation*} 
\FloatBarrier

\begin{remark}
This kind of change of variables corresponds to functionals~$\phi^{(1)}$ that involve evaluations of the state variables at times~$\tau$ and~$\tau+\varepsilon$. Of course these changes of variables must be adapted when considering evaluations at other times (but still in function of~$\tau$). When the functional~$\phi^{(1)}$ involves only evaluations at time~$\tau$, then the one given above is still valid by choosing $\varepsilon = \tilde{\varepsilon} = 0$.
\end{remark}

The time-derivative $\dot{\mu}$ of $\mu$ (with respect to $s$), as well its the partial derivative $\dot{\mu}_\tau$ with respect to $\tau$ are given as follows:
\begin{equation*}
\begin{array}{rcl}
\dot{\mu}(s,\tau) & = & 
\left\{ \begin{array} {ll}
\tau  & \text{if } s\in [0,1), \\[5pt]
\displaystyle  \varepsilon/\tilde{\varepsilon} & \text{if } s \in (1, 1+\tilde{\varepsilon}), \\[5pt]
\displaystyle \frac{T-(\tau+\varepsilon)}{1-\tilde{\varepsilon}} & \text{if } s\in (1+\tilde{\varepsilon}, 2],
\end{array} \right.
\\[30pt]
\dot{\mu}_\tau(s,\tau) & = &  
\left\{ \begin{array} {ll}
1  & \text{if } s\in [0,1), \\
0 & \text{if } s \in (1, 1+\tilde{\varepsilon}), \\
-1/(1-\tilde{\varepsilon}) & \text{if } s\in (1+\tilde{\varepsilon}, 2].
\end{array} \right.
\end{array}
\end{equation*}
Note that $\dot{\mu}_\tau$ is actually independent of $\tau$, and that $\ddot{\mu} = 0$. 
For a given switching time $\tau$, we introduce the following change of unknowns and variables
\begin{equation}
\tilde{u}:s \mapsto u (\cdot,\mu(s,\tau)), \quad \tilde{\mathfrak{p}}:s \mapsto \mathfrak{p}( \mu(s,\tau)), \quad \tilde{\xi}:s \mapsto \xi(\cdot,\mu(s)),
\quad s \in [0,2]. \label{eqChangeOfVar}
\end{equation}
We then transform~\eqref{mainpb} into the following one:
\begin{equation} \tag{$\tilde{\mathcal{P}}$}
\left\{ \begin{array}{l}\displaystyle \max_{\tilde{\xi} \in \mathcal{X}_{p,2}(\omega), \tau \in (0,T)} \int_0^2 \dot{\mu}(s,\tau)c(\tilde{u},\dot{\tilde{u}}/\dot{\mu}, \tilde{\xi})\, \d s + \phi^{(1)}(\tilde{u},\dot{\tilde{u}}/\dot{\mu})(1) + 
\phi^{(2)}(\tilde{u},\dot{\tilde{u}}/\dot{\mu})(2), \\[10pt]
\text{subject to~\eqref{sysmaintilde}},
\end{array} \right. \label{mainpbtilde}
\end{equation}
where~\eqref{sysmaintilde} is the system satisfied by $(\tilde{u}, \tilde{\mathfrak{p}})$, namely
\begin{equation} \label{sysmaintilde}
\begin{array} {rcl}
\ddot{\tilde{u}} - \kappa \dot{\mu}\Delta \dot{\tilde{u}}
- \dot{\mu}^2 \divg \sigma(\nabla \tilde{u}) = \dot{\mu}^2 f(\tilde{\xi}) & & \text{in $\Omega \times(0,2)$},  \\
\displaystyle \frac{\kappa}{\dot{\mu}} \frac{\p \dot{\tilde{u}}}{\p n} +  \sigma(\nabla \tilde{u})n  + \tilde{\mathfrak{p}}\, \cof\left(\Phi(\tilde{u}) \right) n = \tilde{g} & & \text{on $\Gamma_N\times (0,2)$}, \\
\displaystyle \int_{\Omega} \det(\Phi(\tilde{u})) \d \Omega = 
\int_{\Omega} \det(\Phi(u_0)) \d \Omega
& & \text{in } (0,2)  \\
\dot{\tilde{u}} = 0 & & \text{on $\Gamma_D\times(0,2)$},  \\
\tilde{u}(\cdot, 0) = u_0, \displaystyle \quad \dot{\tilde{u}}(\cdot, 0) = \dot{\mu}(0)\udotini & & \text{in $\Omega$}, 
\end{array}
\end{equation}
with $\tilde{g}(\cdot,s) := g(\cdot \mu(s,\tau))$. The interest of this change of unknowns lies in the fact that in the new optimal control problem~\eqref{mainpbtilde}, the two variables to be optimized, namely the time parameter $\tau$ and the control function $\tilde{\xi}$, are no longer coupled. Let us rewrite~\eqref{sysmaintilde} as a first-order system of evolution equations, by introducing
\begin{equation}
(\tilde{y}_0, \tilde{y}_1) = \displaystyle \big(\tilde{u}, \dot{\tilde{u}}/\dot{\mu}\big),
\label{introduce-y}
\end{equation}
so that $(\tilde{y}_0,\tilde{y}_1) =( u\circ \mu, \dot{u} \circ \mu)$. Then $(\tilde{y}_0,\tilde{y}_1,\tilde{\mathfrak{p}})$ satisfies the following system: 
\begin{equation} \label{sysmaintilde2x}
\begin{array} {rcl}
\dot{\tilde{y}}_0 - \dot{\mu} \tilde{y}_1 = 0 
& & \text{in $\Omega \times(0,2)$},  \\
\dot{\tilde{y}}_1 - \kappa \dot{\mu}\Delta \tilde{y}_1
- \dot{\mu} \divg \sigma(\nabla \tilde{y}_0) = \dot{\mu} f(\tilde{\xi}) & & \text{in $\Omega \times(0,2)$}, \\
\kappa \displaystyle \frac{\p \tilde{y}_1}{\p n} +  \sigma(\nabla \tilde{y}_0)n  
+ \tilde{\mathfrak{p}}\, \cof\left(\Phi(\tilde{y}_0) \right) n = \tilde{g} & & \text{on $\Gamma_N\times (0,2)$},  \\
\displaystyle \int_{\Omega} \det(\Phi(\tilde{y}_0)) \d \Omega = 
\int_{\Omega} \det(\Phi(u_0)) \d \Omega
& & \text{in } (0,2) \\
\tilde{y}_1 = 0 & & \text{on $\Gamma_D\times(0,2)$},  \\
\tilde{y}_0(\cdot, 0) = u_0, \displaystyle \quad \dot{\tilde{y}}_1(\cdot, 0) = \udotini & & \text{in $\Omega$}.
\end{array}
\end{equation}

\begin{remark}\label{remark-id-boundary}
Note that deriving in time the constraint~$\displaystyle \int_{\Omega} \det(\Phi(\tilde{y}_0)) \d \Omega = 
\int_{\Omega} \det(\Phi(u_0)) \d \Omega$ combined with the first equation of~\eqref{sysmaintilde2x} yields
\begin{equation*}
\displaystyle
\dot{\mu}\int_{\Omega} \cof(\Phi(\tilde{y}_0)): \nabla \tilde{y}_1 \d \Omega = 0.
\end{equation*} Further, in the same way as we obtained~\eqref{sysmain4bis}, we deduce that
\begin{equation*}
\int_{\Gamma_N} \tilde{y}_1 \cdot \cof(\Phi(\tilde{y}_0))n\, \d \Gamma_N = 0.
\end{equation*}
\end{remark}
Problem~\eqref{mainpbtilde} is equivalent to the following one, for which we use the same notation:\begin{small}
\begin{equation} \tag{$\tilde{\mathcal{P}}$}
\left\{ \begin{array}{l}\displaystyle \max_{\tilde{\xi} \in \mathcal{X}_{p,2}(\omega), \tau \in (0,T)} \left(
J(\tilde{\xi},\tau) =
\int_0^2 \dot{\mu}(s,\tau)c(\tilde{y}_0,\tilde{y}_1, \tilde{\xi})\, \d s + \phi^{(1)}(\tilde{y}_0,\tilde{y}_1)(1) + 
\phi^{(2)}(\tilde{y}_0,\tilde{y}_1)(2)\right), \\[10pt]
\text{subject to~\eqref{sysmaintilde2x}}.
\end{array} \right. \label{mainpbtilde2}
\end{equation}\end{small}

\subsection{The control-to-state mapping} \label{sec-cts}

We first state a result for a general linear system that will be used several times in the rest.
\begin{proposition} \label{prop-superlinear}
Let be $(\tilde{y}_0,\tilde{y}_1,\tilde{\mathfrak{p}}) \in \mathcal{U}_{p,2}(\Omega) \times \dot{\mathcal{U}}_{p,2}(\Omega) \times \mathcal{P}_{p,2}$ and $\tau \in (0,T)$. Assume that $
 \tilde{f}_0 \in \mathcal{U}_{p,2}(\Omega)$, $\tilde{f}_1 \in \mathcal{F}_{p,2}(\Omega)$, and $\tilde{g} \in \mathcal{G}_{p,2}(\Gamma_N)$ with the compatibility condition~$\tilde{g}(\cdot,0) = 0$. Recall that the tensor fields $\sigma_L$ and $\sigma_N$ have been introduced in~\eqref{tensor-linear}. Then, if~$(\tilde{y_0},\tilde{\mathfrak{p}})$ is small enough in~$\mathcal{U}_{p,2}(\Omega) \times \mathcal{P}_{p,2}$, the following system
\begin{equation} \label{syssuperlinear}
\begin{array} {rcl}
\dot{\tilde{z}}_0 - \dot{\mu} \tilde{z}_1 = \dot{\mu}\tilde{f}_0
 & & \text{in $\Omega \times(0,2)$},  \\
\dot{\tilde{z}}_1 - \kappa \dot{\mu}\Delta \tilde{z}_1
- \dot{\mu} \divg (\sigma_L(\nabla \tilde{y}_0).\nabla \tilde{z}_0) = 
\dot{\mu} \tilde{f}_1 & & \text{in $\Omega \times(0,2)$}, \\
\kappa \frac{\p \tilde{z}_1}{\p n} +  
\Big((\sigma_L + \tilde{\mathfrak{p}}\, \sigma_N)(\nabla \tilde{y}_0).\nabla \tilde{z}_0\Big)n 
+   \tilde{\mathfrak{q}}\, \cof\left(\Phi(\tilde{y}_0) \right) n = \tilde{g} & & 
\text{on $\Gamma_N\times (0,2)$}, \\
\tilde{z}_1 = 0 & & \text{on $\Gamma_D\times(0,2)$},  \\
\displaystyle \int_{\Gamma_N} \tilde{z}_0 \cdot \cof(\Phi(\tilde{y}_0))n\, \d \Gamma
= 0 & & \text{in } (0,2)\\
\tilde{z}_0(\cdot, 0) = 0, \displaystyle \quad \tilde{z}_1(\cdot, 0) = 0 & & \text{in $\Omega$}.
\end{array}
\end{equation}
admits a unique solution~$(\tilde{z}_0, \tilde{z}_1,\tilde{\mathfrak{q}}) \in \mathcal{U}_{p,2}(\Omega)\times \dot{\mathcal{U}}_{p,2}(\Omega) \times  \mathcal{P}_{p,2}$. Moreover, it satisfies
\begin{equation*}
\|(\tilde{z}_0,\tilde{z}_1)\|_{\mathcal{U}_{p,2}(\Omega)\times \dot{\mathcal{U}}_{p,2}(\Omega)} \leq C(\tilde{y}_0,\tilde{\mathfrak{p}})
\left( \| \tilde{f}_0\|_{\mathcal{U}_{p,2}(\Omega)} + 
\| \tilde{f}_1\|_{\mathcal{F}_{p,2}(\Omega)}
+ \| \tilde{g}\|_{\mathcal{G}_{p,2}(\Gamma_N)}
\right).
\end{equation*}
where the constant $C(\tilde{y}_0,\tilde{\mathfrak{p}})>0$ depends only on $(\tilde{y}_0,\tilde{y}_1,\tilde{\mathfrak{p}})$.
\end{proposition}

Note here again that solutions~$(\tilde{z}_0,\tilde{z}_1,\tilde{\mathfrak{q}})$ of system~\eqref{syssuperlinear} are continuous on~$[0,2]$ -- with values in $\mathcal{U}_p^{(0,1)}(\Omega) \times \R$.

\begin{proof}
Notice that $(\tilde{z}_0,\tilde{z_1},\tilde{\mathfrak{q}})$ is solution of system~\eqref{syssuperlinear} with~$(\tilde{y}_0,\tilde{y}_1,\tilde{\mathfrak{p}})$ and $(\tilde{f}_0,\tilde{f}_1,\tilde{g})$ as data if and only if~$(z_0,z_1,\mathfrak{q})$ is solution of system~\eqref{sysx} with~$(y_0,y_1,\mathfrak{p})$ and $(f_0,f_1,g)$ as data, where we have
\begin{equation*}
\begin{array} {l}
y_0 = \tilde{y}_0(\cdot, \mu^{-1}(\cdot,\tau)), \quad 
y_1 = \tilde{y}_1 (\cdot, \mu^{-1}(\cdot,\tau)), \quad
\mathfrak{p} = \tilde{\mathfrak{p}}(\mu^{-1}(\cdot,\tau)), \\
f_0 = \tilde{f}_0 (\cdot, \mu^{-1}(\cdot,\tau)), \quad
f_1 = \tilde{f}_1 (\cdot, \mu^{-1}(\cdot,\tau)), \quad
g = \tilde{g} (\cdot, \mu^{-1}(\cdot,\tau)).
\end{array}
\end{equation*}
This is due to the regularity of the change of variables~$\mu(\cdot,\tau)$, in particular the fact that its derivatives and those of~$\mu^{-1}(\cdot,\tau)$ are in~$\L^{\infty}(0,2;\R)$ and~$\L^{\infty}(0,T;\R)$, respectively. Then the result follows from Proposition~\ref{prop-well-syslin-init}, and the announced estimate too, which concludes the proof.
\end{proof}

\subsubsection{Regularity}
In this subsection we study the regularity of the control-to-state mapping, given as
\begin{equation}
\begin{array} {rrcl}
\mathbb{S}: & \mathcal{X}_{p,2}(\omega) \times (0,T) & \rightarrow & \mathcal{U}_{p,2}(\Omega)\times \dot{\mathcal{U}}_{p,2}(\Omega) \times \mathcal{P}_{p,2} \\
& (\tilde{\xi} , \tau) & \mapsto & (\tilde{y}_0, \tilde{y}_1,\tilde{\mathfrak{p}}),
\end{array}
\end{equation}
where $(\tilde{y}_0,\tilde{y}_1,\tilde{\mathfrak{p}})$ is the solution of system~\eqref{sysmaintilde2x} corresponding to~$\tilde{\xi}$ and $\dot{\mu} = \dot{\mu}(\cdot,\tau)$. Let us show that $\mathbb{S}$ is {\it locally} well-defined. More precisely, we state:

\begin{proposition} \label{prop-genesis}
Let be $T\in(0,\infty)$. There exists~$\eta>0$ such that if
\begin{equation*}
\|\dot{\mu}f(\tilde{\xi})+ \divg(\sigma(0))\|_{\mathcal{F}_{p,2}(\Omega)} + 
\|\tilde{g}-\sigma(0)n\|_{\mathcal{G}_{p,2}(\Gamma_N)} +
\|(u_0,\udotini)\|_{\mathcal{U}^{(0,1)}_p(\Omega)} \leq \eta
\end{equation*}
with the compatibility condition $\displaystyle \kappa\frac{\p \dot{u}_0}{\p n} + \sigma(\nabla u_0)n = \tilde{g}(\cdot,0)$ on~$\Gamma_N$, then system~\eqref{sysmaintilde2x} admits a unique solution~$(\tilde{y}_0,\tilde{y}_1,\tilde{\mathfrak{p}}) \in \mathcal{U}_{p,2}(\Omega) \times \dot{\mathcal{U}}_{p,2}(\Omega) \times \mathcal{P}_{p,2}$.
\end{proposition}

\begin{proof}
From~\eqref{eqChangeOfVar} and~\eqref{introduce-y} we have 
\begin{equation*}
\begin{array} {l}
\tilde{y}_0(\cdot,s) = u(\cdot, \mu(s,\tau)), \quad 
\tilde{y}_1(\cdot,s) = \dot{u}(\cdot, \mu(s,\tau)),
\quad \tilde{\mathfrak{p}}(s) = \mathfrak{p}(\mu(s,\tau)),\\[5pt]
\tilde{\xi}(\cdot,s) = \xi(\cdot,\mu(s,\tau)),
\end{array}
\end{equation*}
and then it is clear that $(\tilde{y}_0,\tilde{y}_1,\tilde{\mathfrak{p}})$ is solution of~\eqref{sysmaintilde2x} with $(\tilde{\xi},\tau)$ as control if and only if $(u, \dot{u}, \mathfrak{p})$ is solution of~\eqref{sysmain} with $\xi$ as control , provided that~$\tilde{g}(\cdot,s) = g(\cdot, \mu(s,\tau))$. Therefore we conclude by invoking Proposition~\ref{prop-well-sysmain}.
\end{proof}

We are now in position to prove regularity for the control-to-state mapping.

\begin{theorem} \label{th-cts}
The control-to-state mapping $\mathbb{S}$ is locally of class $\mathcal{C}^1$ from $\mathcal{X}_{p,2}(\omega) \times (0,T)$ onto $\mathcal{U}_{p,2}(\Omega)\times \dot{\mathcal{U}}_{p,2}(\Omega) \times \mathcal{P}_{p,2}$.
\end{theorem}

\begin{proof}
The result is an application of the implicit function theorem. Define on the space $\mathcal{U}_{p,2}(\Omega)\times \dot{\mathcal{U}}_{p,2}(\Omega) \times \mathcal{P}_{p,2}
 \times  \mathcal{X}_{p,2}(\omega) \times  (0,T)$ the mapping
\begin{equation*}
	e:  (\tilde{y}_0, \tilde{y}_1,\tilde{\mathfrak{p}},\tilde{\xi},\tau) \mapsto 
	\left( \begin{array}{c}
	\dot{\tilde{y}}_0 - \dot{\mu} \tilde{y}_1 \\
	\dot{\tilde{y}}_1 - \dot{\mu}\left(\kappa \Delta \tilde{y}_1 
- \divg(\sigma(\nabla \tilde{y}_0)) - f(\tilde{\xi})\right) \\
\kappa \displaystyle \frac{\p \tilde{y}_1}{\p n} + \sigma(\nabla \tilde{y}_0)n 
+\mathfrak{p}\, \cof(\Phi(\tilde{y}_0))n - g \\[5pt]
\displaystyle \int_{\Omega} \det(\Phi(\tilde{y}_0)) \d \Omega
-\int_{\Omega} \det(\Phi(u_0))\d \Omega\\
	(\tilde{y}_0(\cdot,0), \tilde{y}_1(\cdot,0)) - (u_0,\udotini) 
	\end{array} \right),
\end{equation*}\noindent with values in $\dot{\mathcal{U}}_{p,2}(\Omega) \times	
\mathcal{F}_{p,2}(\Omega) \times \mathcal{G}_{p,2}(\Gamma_N) \times 
\mathcal{H}_{p,2} \times\mathcal{U}_p^{(0,1)}(\Omega)$, where the dependence on~$\tau$ lies in $\dot{\mu} = \dot{\mu}(\cdot,\tau)$. From Proposition~\ref{prop-genesis}, the mapping $e$ is locally well-defined, and the equality $e(\mathbb{S}(\tilde{\xi},\tau), \tilde{\xi},\tau) = 0$ holds for all $(\tilde{\xi}, \tau) \in \mathcal{X}_{p,2}(\omega) \times (0,T)$.
Furthermore, from assumption~$\mathbf{A1}$, the mapping $e$ is of class $\mathcal{C}^1$. Proposition~\ref{prop-superlinear} shows that the derivative of~$e$ with respect to $(\tilde{y}_0,\tilde{y}_1,\tilde{\mathfrak{p}})$ is invertible. Then the implicit function theorem provides us the existence of a $\mathcal{C}^1$ mapping on $\mathcal{U}_{p,2}(\Omega) \times \dot{\mathcal{U}}_{p,2}(\Omega)\times \mathcal{P}_{p,2}$ that coincides with $\mathbb{S}$, which concludes the proof.
\end{proof}

We can now describe the partial derivatives of the control-to-state mapping.

\subsubsection{Partial derivatives}
Let us introduce the linear system satisfied by $(\tilde{v}_0,\tilde{v}_1,\tilde{\mathfrak{q}}):=\mathbb{S}'_{\tilde{\xi}}(\tilde{\xi},\tau).\hat{\xi}$, denoting the sensitivity of $\mathbb{S}$ with respect to variable $\tilde{\xi}$ in the direction~$\hat{\xi}$, at point $(\tilde{y}_0,\tilde{y}_1,\tilde{\mathfrak{p}}) = \mathbb{S}(\tilde{\xi}, \tau)$:
\begin{equation} \label{sysprimexi}
\begin{array} {rcl}
\dot{\tilde{v}}_0 - \dot{\mu} \tilde{v}_1 = 0
 & & \text{in $\Omega \times(0,2)$},  \\
\dot{\tilde{v}}_1 - \kappa \dot{\mu}\Delta \tilde{v}_1
- \dot{\mu} \divg (\sigma_L(\nabla \tilde{y}_0).\nabla \tilde{v}_0) = 
\dot{\mu} f'(\tilde{\xi}).\hat{\xi} & & \text{in $\Omega \times(0,2)$},  \\
\kappa \displaystyle \frac{\p \tilde{v}_1}{\p n} +  
\Big((\sigma_L + \tilde{\mathfrak{p}}\, \sigma_N)(\nabla \tilde{y}_0).\nabla \tilde{v}_0\Big)n 
 + \tilde{\mathfrak{q}}\, \cof(\Phi(\tilde{y}_0)) n = 0 & & \text{on $\Gamma_N\times (0,2)$}, 
\\
\tilde{v}_1 = 0 & & \text{on $\Gamma_D\times(0,2)$},  \\
\displaystyle \int_{\Gamma_N} \tilde{v}_0 \cdot \cof(\Phi(\tilde{y}_0))n\, \d \Gamma_N
= 0 & & \text{in } (0,2)\\
\tilde{v}_0(\cdot, 0) = 0, \displaystyle \quad \tilde{v}_1(\cdot, 0) = 0 & & \text{in $\Omega$}.
\end{array}
\end{equation}

\noindent Note that when $\xi \mapsto f(\xi)$ is linear (see section~\ref{sec-app-control}), obviously we have $f'(\tilde{\xi}).\hat{\xi} = f(\hat{\xi})$ for all $\tilde{\xi} \in \mathcal{X}_{p,2}(\omega)$. 
Let us state that system~\eqref{sysprimexi} is well-posed.

\begin{proposition} \label{proplinxi}
Assume that~$(\tilde{\xi},\tau) \in \mathcal{X}_{p,2}(\omega) \times (0,T)$, and denote $(\tilde{y}_0, \tilde{y}_1,\tilde{\mathfrak{p}}) = \mathbb{S}(\tilde{\xi},\tau) \in \mathcal{U}_{p,2}(\Omega) \times \dot{\mathcal{U}}_{p,2}(\Omega)\times \mathcal{P}_{p,2}$. Then, if~$(\tilde{y_0},\tilde{\mathfrak{p}})$ is small enough in~$\mathcal{U}_{p,2}(\Omega) \times \mathcal{P}_{p,2}$, system~\eqref{sysprimexi} admits a unique solution~$(\tilde{v}_0, \tilde{v}_1, \tilde{\mathfrak{q}}) \in \mathcal{U}_{p,2}(\Omega)\times \dot{\mathcal{U}}_{p,2}(\Omega) \times \mathcal{P}_{p,2}$ for all~$\hat{\xi} \in \mathcal{X}_{p,2}(\omega)$. Moreover, there exists a constant $C(\tilde{y}_0,\tilde{\mathfrak{p}})$ depending only on $(\tilde{y}_0,\tilde{\mathfrak{p}})$ such that
\begin{equation*}
\| \tilde{v}_0\|_{\mathcal{U}_{p,2}(\Omega)}+
\| \tilde{v}_1\|_{\dot{\mathcal{U}}_{p,2}(\Omega)} + 
\| \tilde{\mathfrak{q}} \|_{\mathcal{P}_{p,2}} 
\leq  C(\tilde{y}_0,\tilde{\mathfrak{p}}) 
\|f'(\tilde{\xi}).\hat{\xi}\|_{\mathcal{F}_{p,2}(\Omega)}.
\end{equation*}
\end{proposition}

\begin{proof}
This is a consequence of Proposition~\ref{prop-superlinear} with $\tilde{f}_0 = 0$, $\tilde{f}_1 = f'(\hat{\xi}).\tilde{\xi}$ and $\tilde{g}=0$.
\end{proof}

We also introduce the linear system satisfied by $(\tilde{w}_0,\tilde{w}_1,\tilde{\mathfrak{r}}):=\mathbb{S}'_{\tau}(\tilde{\xi},\tau)$, denoting the sensitivity of $\mathbb{S}$ with respect to~$\tau$ at point $(\tilde{y}_0,\tilde{y}_1,\tilde{\mathfrak{p}}) = \mathbb{S}(\tilde{\xi},\tau)$:
\begin{equation} \label{sysprimetau}
\begin{array} {rcl}
\dot{\tilde{w}}_0 - \dot{\mu} \tilde{w}_1 = \dot{\mu}_{\tau} \tilde{y}_1
& & \text{in $\Omega \times(0,2)$}, \\
\dot{\tilde{w}}_1 - \kappa \dot{\mu}\Delta \tilde{w}_1
- \dot{\mu} \divg (\sigma_L(\nabla \tilde{y}_0).\nabla \tilde{w}_0) = &  &\\
\dot{\mu}_{\tau}
\left( \kappa \Delta \tilde{y}_1 + \divg(\sigma(\nabla \tilde{y}_0)) + f(\tilde{\xi}) \right)
& & \text{in $\Omega \times(0,2)$},  \\
\kappa \displaystyle \frac{\p \tilde{w}_1}{\p n} +  
\Big((\sigma_L + \tilde{\mathfrak{p}}\, \sigma_N)(\nabla \tilde{y}_0).\nabla \tilde{w}_0\Big)n   + 
\tilde{\mathfrak{r}} \, \cof\left(\Phi(\tilde{y}_0) \right) n = 0 & & \text{on $\Gamma_N\times (0,2)$}, \\
\tilde{w}_1 = 0 & & \text{on $\Gamma_D\times(0,2)$},  \\
\displaystyle \int_{\Gamma_N} \tilde{w}_0 \cdot \cof(\phi(\tilde{y}_0))n\, \d \Gamma_N
 = 0 & & \text{in } (0,2)\\
\tilde{w}_0(\cdot, 0) = 0, \displaystyle \quad \tilde{w}_1(\cdot, 0) = 0 & & \text{in $\Omega$}.
\end{array}
\end{equation}

\noindent We show that system~\eqref{sysprimetau} is also well-posed.

\begin{proposition} \label{proplintau}
Assume that~$(\tilde{\xi},\tau) \in \mathcal{X}_{p,2}(\omega) \times (0,T)$, and denote $(\tilde{y}_0, \tilde{y}_1,\tilde{\mathfrak{p}}) = \mathbb{S}(\tilde{\xi},\tau) \in \mathcal{U}_{p,2}(\Omega) \times \dot{\mathcal{U}}_{p,2}(\Omega)\times \mathcal{P}_{p,2}$. Then, if~$(\tilde{y_0},\tilde{\mathfrak{p}})$ is small enough in~$\mathcal{U}_{p,2}(\Omega) \times \mathcal{P}_{p,2}$, system~\eqref{sysprimetau} admits a unique solution~$(\tilde{w}_0,\tilde{w}_1,\tilde{\mathfrak{r}}) \in \mathcal{U}_{p,2}(\Omega)\times \dot{\mathcal{U}}_{p,2}(\Omega) \times \mathcal{P}_{p,2}$. Moreover, there exists a constant $C(\tilde{y}_0,\tilde{\mathfrak{p}})$ depending only on $(\tilde{y}_0,\tilde{\mathfrak{p}})$ such that
\begin{equation*}
\| \tilde{w}_0\|_{\mathcal{U}_{p,2}(\Omega)} +
 \|\tilde{w}_1\|_{\dot{\mathcal{U}}_{p,2}(\Omega)} +
 \| \tilde{\mathfrak{r}} \|_{\mathcal{P}_{p,2}} 
\leq  C(\tilde{y}_0,\tilde{\mathfrak{p}}) 
\left(1+ \|f(\tilde{\xi})\|_{\mathcal{F}_{p,2}(\Omega)}\right).
\end{equation*}
\end{proposition}

\begin{proof}
This is a consequence of Proposition~\ref{prop-superlinear} with
\begin{equation*}
\tilde{f}_0 = \frac{\dot{\mu}_{\tau}}{\dot{\mu}}\tilde{y}_1 
\in \dot{\mathcal{U}}_{p,2}(\Omega), \quad
\tilde{f}_1 = \displaystyle \frac{\dot{\mu}_{\tau}}{\dot{\mu}}
\left( \kappa \Delta \tilde{y}_1 + \divg(\sigma(\nabla \tilde{y}_0)) + f(\tilde{\xi}) \right)
\in \mathcal{F}_{p,2}(\Omega),
\end{equation*}
and $\tilde{g}=0$. Since~$\dot{\mu}$ and~$\dot{\mu}_{\tau}$ are in~$\L^{\infty}(0,2;\R)$, the parameter~$\tau$ does not appear in the dependence of the constant~$C(\tilde{y}_0,\tilde{\mathfrak{p}})$.
\end{proof}

\subsection{The adjoint system} \label{sec-adjoint}
Recall the notation introduced in system~\eqref{sysadjoint-init}: For a function~$\varphi$ continuous on~$(0,1)\cup(1,2)$ we define the jump of~$\varphi$ at $s=1$ as folows:
\begin{equation*}
\left[\varphi\right]_1 :=
 \lim_{s\rightarrow 1^+} \varphi(s) - \lim_{s\rightarrow 1^-} \varphi(s).
\end{equation*}
Let be $(\tilde{y}_0, \tilde{y}_1,\tilde{\mathfrak{p}}) \in \mathcal{U}_{p,2}(\Omega)\times \dot{\mathcal{U}}_{p,2}(\Omega) \times \mathcal{P}_{p,2} $. We associate with $(\tilde{y}_0, \tilde{y}_1,\tilde{\mathfrak{p}})$ the adjoint state~$(\tilde{\zeta}_0, \tilde{\zeta}_1,\tilde{\pi})$, assumed to satisfy the following system\begin{small}
\begin{equation} \label{sysadjoint}
\begin{array} {rcl}
-\dot{\tilde{\zeta}}_0 - \dot{\mu} \divg(\sigma_L(\nabla \tilde{y}_0)^{\ast}.\nabla \tilde{\zeta}_1) = 
-\dot{\mu}c'_{y_0}(\tilde{y}_0, \tilde{y}_1, \tilde{\xi}) 
& & \text{in $\Omega \times\left((0,1) \cup (1,2)\right)$},  \\
-\dot{\tilde{\zeta}}_1 - \dot{\mu} \tilde{\zeta}_0 -\kappa \dot{\mu} \Delta \tilde{\zeta}_1 = 
-\dot{\mu}c'_{y_1}(\tilde{y}_0,\tilde{y}_1, \tilde{\xi}) 
& & \text{in $\Omega \times\left((0,1) \cup (1,2)\right)$},  \\ 
\Big((\sigma_L + \tilde{\mathfrak{p}}\, \sigma_N)(\nabla \tilde{y}_0)^{\ast}.\nabla \tilde{\zeta}_1\Big)n  
+ \tilde{\pi}\, \cof(\Phi(\tilde{y}_0)) n = 0 & & 
\text{on $\Gamma_N\times \left((0,1) \cup (1,2)\right)$},  \\
\kappa \displaystyle \frac{\p \tilde{\zeta}_1}{\p n}  =0
 & & 
\text{on $\Gamma_N\times \left((0,1) \cup (1,2)\right)$},\\[5pt]
\tilde{\zeta}_1 = 0 & & \text{on $\Gamma_D\times\left((0,1) \cup (1,2)\right)$},  \\
\left\langle \tilde{\zeta}_1 \, ; \cof(\Phi(\tilde{y}_0))n\right\rangle_{\WW^{1/(p'),p}(\Gamma_N)',\WW^{1/(p'),p}(\Gamma_N)} = 0
& & \text{in } \left((0,1) \cup (1,2)\right) \\
\left[ \tilde{\zeta}_0 \right]_1 = \phi^{(1)'}_{y_0}(\tilde{y}_0,\tilde{y}_1)(1)
\quad \text{and} \quad 
\left[ \tilde{\zeta}_1 \right]_1 = \phi^{(1)'}_{y_1}(\tilde{y}_0,\tilde{y}_1)(1)
& & \text{in $\Omega$},
\\
\tilde{\zeta}_0(\cdot, 2) = -\phi^{(2)'}_{y_0}(\tilde{y}_0,\tilde{y}_1)(2), 
\displaystyle 
\quad \tilde{\zeta}_1(\cdot, 2) = -\phi^{(2)'}_{y_1}(\tilde{y}_0,\tilde{y}_1)(2) & & \text{in $\Omega$}.
\end{array}
\end{equation}\end{small}\noindent 
Note that the derivatives of mapping $\phi^{(1)}$ induce jumps for variables $\tilde{\zeta}_0$ and $\tilde{\zeta}_1$ at time $s=1$. Similarly to Definition~\ref{def-trans-init} that deals with solutions of system~\eqref{sysadjoint-init}, we define solutions of system~\eqref{sysadjoint} by {\it transposition} as follows:

\begin{definition} \label{def-trans}
Let be $(\tilde{y}_0, \tilde{y}_1,\tilde{\mathfrak{p}},\tilde{\xi}, \tau) \in \mathcal{U}_{p,2}(\Omega) \times \dot{\mathcal{U}}_{p,2}(\Omega)\times \mathcal{P}_{p,2}\times \mathcal{X}_{p,2}(\omega) \times (0,T)$. We say that~$(\tilde{\zeta}_0,\tilde{\zeta}_1,\tilde{\pi})$ is a solution of~\eqref{sysadjoint} associated with $(\tilde{y}_0,\tilde{y}_1,\tilde{\mathfrak{p}})$, if for all $(\tilde{f}_0, \tilde{f}_1,\tilde{g}) \in \dot{\mathcal{U}}_{p,2}(\Omega) \times \mathcal{F}_{p,2}(\Omega) \times \mathcal{G}_{p,2}(\Gamma_N)$ we have\begin{small}
\begin{equation}
\begin{array} {l}
 \displaystyle \left\langle \tilde{\zeta}_0 ; \dot{\mu}\tilde{f}_0 \right\rangle_{\dot{\mathcal{U}}_{p,2}(\Omega)',\dot{\mathcal{U}}_{p,2}(\Omega)}
 +\left\langle \tilde{\zeta}_1 ; \dot{\mu}\tilde{f}_1 \right\rangle_{\mathcal{F}_{p,2}(\Omega)',\mathcal{F}_{p,2}(\Omega)}
+ \left\langle \tilde{\zeta}_1 ; \dot{\mu}\tilde{g} \right\rangle_{\mathcal{G}_{p,2}(\Gamma_N)',\mathcal{G}_{p,2}(\Gamma_N)} 
\\
 \displaystyle =
-\int_0^2\dot{\mu}\left\langle c'_{y_0}(\tilde{y}_0,\tilde{y}_1, \tilde{\xi})\, ; \tilde{z}_0 \right\rangle_{\WW^{2,p}(\Omega)',\WW^{2,p}(\Omega)} \d s
- \int_0^2\dot{\mu}\left\langle  c'_{y_1}(\tilde{y}_0,\tilde{y}_1, \tilde{\xi})\, ; \tilde{z}_1 \right\rangle_{\WW^{2,p}(\Omega)',\WW^{2,p}(\Omega)} \d s
 \\
 - \left\langle(\phi^{(1)'}_{y_0}(\tilde{y}_0,\tilde{y}_1), \phi^{(1)'}_{y_1}(\tilde{y}_0,\tilde{y}_1) )\, ; (\tilde{z}_0(\cdot,1),\tilde{z}_1(\cdot,1))
\right\rangle_{\mathcal{U}_p^{(0,1)}(\Omega)',\mathcal{U}_p^{(0,1)}(\Omega)}
 \\
 -\left\langle 
(\phi^{(2)'}_{y_0}(\tilde{y}_0,\tilde{y}_1), \phi^{(2)'}_{y_1}(\tilde{y}_0,\tilde{y}_1) )\, ; (\tilde{z}_0(\cdot,2),\tilde{z}_1(\cdot,2))
\right\rangle_{\mathcal{U}_p^{(0,1)}(\Omega)',\mathcal{U}_p^{(0,1)}(\Omega)} ,
\label{id-def-trans}
\end{array}
\end{equation}\end{small}\noindent where~$(\tilde{z}_0,\tilde{z}_1,\tilde{\mathfrak{q}})$ is the solution of system~\eqref{syssuperlinear} with $(\tilde{y}_0,\tilde{y}_1,\tilde{\mathfrak{p}})$ and $(\tilde{f}_0,\tilde{f}_1,\tilde{g})$ as data.
\end{definition}

\noindent It is clear that $(\tilde{\zeta}_0,\tilde{\zeta}_1,\tilde{\pi})$ is solution of system~\eqref{sysadjoint} associated with~$(\tilde{y}_0,\tilde{y}_1,\tilde{\mathfrak{p}})$ in the sense of Definition~\ref{def-trans} if and only if $(\zeta_0,\zeta_1,\pi)$ is solution of system~\eqref{sysadjoint-init} associated with~$(y_0,y_1,\mathfrak{p})$ in the sense of Definition~\ref{def-trans-init}, provided that
\begin{equation*}\begin{array} {l}
\tilde{y}_0(\cdot,s) = y_0(\cdot, \mu(s,\tau)), \quad 
\tilde{y}_1(\cdot,s) = y_1(\cdot, \mu(s,\tau)),
\quad \tilde{\mathfrak{p}}(s) = \mathfrak{p}(\mu(s,\tau)),\\ 
\tilde{\xi}(\cdot,s) = \xi(\cdot,\mu(s,\tau)), \\
\tilde{f}_0(\cdot,s) = f_0(\cdot,\mu(s,\tau)), \quad
\tilde{f}_1(\cdot,s) = f_1(\cdot,\mu(s,\tau)), \quad
\tilde{g}(\cdot,s) = g(\cdot,\mu(s,\tau)).
\end{array}
\end{equation*}
The solutions of systems~\eqref{sysadjoint-init} and~\eqref{sysadjoint} then satisfy the relations
\begin{equation*}
\tilde{\zeta}_0(\cdot,s) = \zeta_0(\cdot, \mu(s,\tau)), \quad 
\tilde{\zeta}_1(\cdot,s) = \zeta_1(\cdot, \mu(s,\tau)),
\quad \tilde{\pi}(s) = \pi(\mu(s,\tau)).
\end{equation*}
Therefore we rely on Proposition~\ref{propadj-init} for stating the following result:

\begin{proposition} \label{propadj}
Let be $(\tilde{y}_0,\tilde{y}_1,\tilde{\mathfrak{p}},\tilde{\xi},\tau) \in \mathcal{U}_{p,2}(\Omega) \times \dot{\mathcal{U}}_{p,2}(\Omega)\times \mathcal{P}_{p,2}\times \mathcal{X}_{p,2}(\omega) \times (0,T)$. System~\eqref{sysadjoint} admits a unique solution $(\tilde{\zeta}_0,\tilde{\zeta}_1,\tilde{\pi}) \in \dot{\mathcal{U}}_{p,2}(\Omega)' \times \mathcal{F}_{p,2}(\Omega)' \times \mathcal{G}_{p,2}(\Gamma_N)'$, in the sense of Definition~\ref{def-trans}. Moreover, there exists a constant $C(\tilde{y}_0,\tilde{\mathfrak{p}}) >0$ depending only on $(\tilde{y}_0,\tilde{\mathfrak{p}},\tau)$ such that
\begin{equation*}
\begin{array} {rcl}
\|(\tilde{\zeta}_0,\tilde{\zeta}_1) \|_{\dot{\mathcal{U}}_{p,2}(\Omega)' \times \mathcal{F}_{p,2}(\Omega)'}
& \leq & C(\tilde{y}_0, \tilde{\mathfrak{p}}) \left(
\| c'_{y_0}(\tilde{y}_0,\tilde{y}_1,\tilde{\xi}) \|_{\mathcal{U}_{p,2}(\Omega)'} +
\| c'_{y_1}(\tilde{y}_0,\tilde{y}_1,\tilde{\xi}) \|_{\dot{\mathcal{U}}_{p,2}(\Omega)'} \right.\\
& &  + \| (\phi^{(1)'}_{y_0}(\tilde{y}_0,\tilde{y}_1), \phi^{(1)'}_{y_1}(\tilde{y}_0,\tilde{y}_1) \|_{\mathcal{U}^{(0,1)}_p(\Omega)'}\\
& &  \left.+
\| (\phi^{(2)'}_{y_0}(\tilde{y}_0,\tilde{y}_1), \phi^{(2)'}_{y_1}(\tilde{y}_0,\tilde{y}_1) \|_{\mathcal{U}^{(0,1)}_p(\Omega)'}
\right).
\end{array}
\end{equation*}
In particular, $C(\tilde{y}_0,\tilde{\mathfrak{p}})$ is independent of $c$, $\phi^{(1)}$ and $\phi^{(2)}$.
\end{proposition}

\subsection{First-order necessary optimality conditions} \label{sec-optcond_1st}

Introduce the functional of problem~\eqref{mainpbtilde2}:
\begin{equation}
\begin{array}{rccl}
\tilde{J}: & \mathcal{X}_{p,2}(\omega) \times (0,T) & \rightarrow & \R\\
& (\tilde{\xi},\tau) & \mapsto & 
\begin{array} {l}\displaystyle 
\int_0^2 \dot{\mu}(s,\tau)c(\mathbb{S}(\tilde{\xi}, \tau)(s),\tilde{\xi}(s))\d s\\
+ \phi_1(\mathbb{S}(\tilde{\xi}, \tau))(1) 
+ \phi_2(\mathbb{S}(\tilde{\xi}, \tau))(2). 
\end{array}
\end{array} \label{def-J}
\end{equation}
Define the Hamiltonian~$\mathcal{H}$ for problem~\eqref{mainpbtilde2}, formally given by
\begin{equation}
\begin{array} {rcl}
\mathcal{H} \big(y_0,y_1,\mathfrak{p},\xi, \zeta_0, \zeta_1, \pi\big) & := &
c(y_0,y_1, \xi)  -
\left\langle \zeta_1 ;  f(\xi)\right\rangle_{\LL^{p'}(\Omega), \LL^{p}(\Omega)}
\\[5pt]  &  &
- \left\langle \zeta_1; g \right\rangle_{\WW^{1/(p'),p}(\Gamma_N)',\WW^{1/(p'),p}(\Gamma_N)}\\[5pt]
& &- \langle \zeta_0 \, ; y_1 \rangle_{\WW^{2,p}(\Omega)', \WW^{2,p}(\Omega)}\\[5pt]
& &  +\left\langle \nabla\zeta_1 ; \kappa\nabla y_1 + 
\sigma(\nabla y_0)\right\rangle_{\WWW^{1,p}(\Omega)', \WWW^{1,p}(\Omega)} \\[5pt]
& & + \pi \displaystyle
\int_{\Omega} \det(\Phi(y_0)) \d \Omega \\[5pt]
& & + \mathfrak{p} \left\langle \zeta_1 \, ; \cof(\Phi(y_0))n\right\rangle_{\WW^{1/(p'),p}(\Gamma_N)',\WW^{1/(p'),p}(\Gamma_N)}.
\end{array}
\label{eqHami}
\end{equation}

\noindent We use the results of sections~\ref{sec-cts} and~\ref{sec-adjoint} to calculate the first-order derivatives of $\tilde{J}$.

\begin{proposition} \label{prop-optcond}
The functional~$\tilde{J}$ is of class~$\mathcal{C}^1$ and its first-order derivatives write as follows
\begin{subequations} \label{id-deriv}
\begin{eqnarray}
\tilde{J}'_{\tilde{\xi}}(\tilde{\xi},\tau).\hat{\xi} & = & \dot{\mu} \mathcal{H}_{\xi}\big(\tilde{y}_0,\tilde{y}_1,\tilde{\mathfrak{p}}, \tilde{\xi}, \tilde{\zeta}_0, \tilde{\zeta}_1,\tilde{\pi}\big).\hat{\xi}, \label{id-deriv1} \\
\tilde{J}'_{\tau}(\tilde{\xi},\tau) & = & \int_0^2\dot{\mu}_{\tau} (s,\tau)
\mathcal{H}\big(\tilde{y}_0,\tilde{y}_1,\tilde{\mathfrak{p}}, \tilde{\xi}, \tilde{\zeta}_0, \tilde{\zeta}_1,\tilde{\pi}\big)(s) \d s, \label{id-deriv2}
\end{eqnarray}
\end{subequations}
for all $\hat{\xi} \in \mathcal{X}_{p,2}(\omega)$, where~$(\tilde{y}_0,\tilde{y}_1,\tilde{\mathfrak{p}})$ satisfies system~\eqref{sysmaintilde2x} corresponding to~$(\tilde{\xi},\tau)$, and~$(\tilde{\zeta}_0, \tilde{\zeta}_1,\tilde{\pi})$ satisfies system~\eqref{sysadjoint} associated with~$(\tilde{y}_0, \tilde{y}_1,\tilde{\mathfrak{p}})$ and~$(\tilde{\xi},\tau)$.
\end{proposition}

\begin{proof}
Denote~$(\tilde{y}_0,\tilde{y}_1) = \mathbb{S}(\tilde{\xi},\tau)$. For the sake of concision we will denote~$\varphi(s)$ when dealing with~$\varphi(\cdot,s)$, for any $s\in[0,2]$. Differentiating functional~$\tilde{J}$ with respect to the variable~$\tilde{\xi}$ gives
\begin{equation*}
\begin{array} {rcl}
J'_{\tilde{\xi}}(\tilde{\xi},\tau).\hat{\xi} & = & 
\displaystyle \int_0^2  \dot{\mu}(s,\tau) 
c'_{\xi}\big(\tilde{y}_0(s),\tilde{y}_1(s),\tilde{\xi}(s)\big).\hat{\xi}(s)\, \d s \\
& & +\displaystyle \int_0^2 \dot{\mu}(s,\tau) c'_{y_0}\big(\tilde{y}_0(s),\tilde{y}_1(s),\tilde{\xi}(s)\big).\tilde{v}_0(s)\, \d s
\\
& & \displaystyle + \int_0^2 \dot{\mu}(s,\tau) c'_{y_1}\big(\tilde{y}_0(s),\dot{\tilde{y}}_1(s),\tilde{\xi}(s)\big).\tilde{v}_1(s) \, \d s 
\\[10pt]
& & +  
\left\langle\left(\tilde{\phi}^{(1)'}_{y_0}(\tilde{y}_0,\tilde{y}_1)(1),  \tilde{\phi}^{(1)'}_{y_1}(\tilde{y}_0,\tilde{y}_1)(1)\right);(\tilde{v}_0(1),\tilde{v}_1(1))
\right\rangle_{\mathcal{U}_p^{(0,1)}(\Omega)',\mathcal{U}_p^{(0,1)}(\Omega)}\\
& & +
 \left\langle\left(\tilde{\phi}^{(2)'}_{y_0}(\tilde{y}_0,\tilde{y}_1)(2),
  \tilde{\phi}^{(2)'}_{y_1}(\tilde{y}_0,\tilde{y}_1)(2)\right);(\tilde{v}_0(2),\tilde{v}_1(2))
\right\rangle_{\mathcal{U}_p^{(0,1)}(\Omega)',\mathcal{U}_p^{(0,1)}(\Omega)},
\end{array}
\end{equation*}
where~$(\tilde{v}_0,\tilde{v}_1) := \mathbb{S}_{\tilde{\xi}}(\tilde{\xi}, \tau).\hat{\xi}$ satisfies system~\eqref{sysprimexi}. Taking the duality product in~$\LL^{p'}(\Omega)\times \LL^p(\Omega)$ of the second equation of~\eqref{sysprimexi} by~$\tilde{\zeta}_1$, integrating by parts on $(0,1) \cup (1,2)$, and using the Green formula two times, leads us to
\begin{equation*}
\begin{array} {rcl}
J'_{\tilde{\xi}}(\tilde{\xi},\tau).\hat{\xi} & = & \displaystyle\int_0^2 \dot{\mu}(s,\tau)c'_{\xi}\big(\tilde{y}_0(s),\tilde{y}_1(s),\tilde{\xi}(s)\big).\hat{\xi}(s)\, \d s \\
& &  \displaystyle -
\int_0^2 \langle \tilde{\zeta}_1(s)\, ; \dot{\mu}(s,\tau) f'(\tilde{\xi}(s)).\hat{\xi} \rangle_{\LL^{p'}(\Omega), \LL^p(\Omega)} \d s.
\end{array}
\end{equation*}
Noticing that
\begin{equation*}
\begin{array} {rcl}
\mathcal{H}_{\xi}\big(\tilde{y}_0,\tilde{y}_1, \tilde{\mathfrak{p}},\tilde{\xi}, \tilde{\zeta}_0, \tilde{\zeta}_1,\tilde{\pi}\big).\hat{\xi} & = & c'_{\xi}(\tilde{y}_0,\tilde{y}_1,\tilde{\xi}).\hat{\xi} - \langle \tilde{\zeta}_1; f'(\tilde{\xi}).\hat{\xi} \rangle_{\LL^{p'}(\Omega),\LL^p(\Omega)}\\[5pt]
& = & \left\langle c'_{\xi}(\tilde{y}_0,\tilde{y}_1,\tilde{\xi}) - f'(\tilde{\xi})^{\ast}.\tilde{\zeta}_1; \hat{\xi} \right\rangle_{\mathcal{X}_{p,2}(\omega)',\mathcal{X}_{p,2}(\omega)}, 
\end{array}
\end{equation*}
identity~\eqref{id-deriv1} follows. Differentiating~$\tilde{J}$ with respect to~$\tau$ gives
\begin{eqnarray*}
\begin{array} {rcl}
J'_{\tau}(\tilde{\xi},\tau) & = & \displaystyle
 \int_0^2 \dot{\mu}_{\tau}(s,\tau) c\big(\tilde{y}_0(s),\tilde{y}_1(s),\tilde{\xi}(s)(s)\big)\, \d s\\
& & + \displaystyle \int_0^2  \dot{\mu}(s,\tau) 
c'_{y_0}\big(\tilde{y}_0(s),\tilde{y}_1(s),\tilde{\xi}(s)\big).\tilde{w}_0(s)\, \d s \\
& & + \displaystyle \int_0^2  \dot{\mu}_{\tau}(s,\tau)
+ c'_{y_1}\big(\tilde{y}_0(s),\tilde{y}_1(s),\tilde{\xi}(s)\big).\tilde{w}_1(s)\,  \d s 
\\[10pt]
& & 
+  
\left\langle\left(\tilde{\phi}^{(1)'}_{y_0}(\tilde{y}_0,\tilde{y}_1)(1),  \tilde{\phi}^{(1)'}_{y_1}(\tilde{y}_0,\tilde{y}_1)(1)\right);(\tilde{w}_0(1),\tilde{w}_1(1))
\right\rangle_{\mathcal{U}_p^{(0,1)}(\Omega)',\mathcal{U}_p^{(0,1)}(\Omega)}\\
& & +
 \left\langle\left(\tilde{\phi}^{(2)'}_{y_0}(\tilde{y}_0,\tilde{y}_1)(2),
  \tilde{\phi}^{(2)'}_{y_1}(\tilde{y}_0,\tilde{y}_1)(2)\right);(\tilde{w}_0(2),\tilde{w}_1(2))
\right\rangle_{\mathcal{U}_p^{(0,1)}(\Omega)',\mathcal{U}_p^{(0,1)}(\Omega)},
\end{array}
\end{eqnarray*}
where~$(\tilde{w}_0,\tilde{w}_1) := \mathbb{S}_{\tau}(\tilde{\xi}, \tau)$ satisfies system~\eqref{sysprimetau}. Taking the inner product of the second equation of~\eqref{sysprimetau} by~$\tilde{\zeta}_1$, integrating by parts and using the Green formula, two times, leads us to\begin{small}
\begin{equation*}
\begin{array} {rcl}
J'_{\tau}(\tilde{\xi},\tau) & = & \displaystyle
\int_0^2 \dot{\mu}_{\tau}(s,\tau)c(\tilde{y}_0(s),\tilde{y}_1(s),\tilde{\xi}(s)) \d s -
\int_0^2 \langle \tilde{\zeta}_0(s) ; \dot{\mu}_{\tau}(s,\tau) \tilde{y}_1(s) \rangle_{\WW^{2,p}(\Omega)', \WW^{2,p}(\Omega)} \d s \\
& & \displaystyle
- \int_0^2 \left\langle \tilde{\zeta}_1(s) ; \dot{\mu}_{\tau}(s,\tau)
\left(\kappa \Delta \tilde{y}_1(s) + \divg(\sigma(\nabla \tilde{y}_0(s))) + f(\tilde{\xi}(s)) \right) \right\rangle_{\LL^{p'}(\Omega), \LL^{p}(\Omega)} \d s,
\end{array}
\end{equation*}\end{small}\noindent where in particular we have used the identity $\displaystyle \int_{\Gamma_N} \tilde{y}_1 \cdot \cof(\Phi(\tilde{y}_0))n\, \d \Gamma_N = 0$ (see Remark~\ref{remark-id-boundary}). We have also used $\displaystyle \left\langle \tilde{\zeta}_1 \, ; \cof(\Phi(\tilde{y}_0))n\right\rangle_{\WW^{1/(p'),p}(\Gamma_N)',\WW^{1/(p'),p}(\Gamma_N)} = 0$, imposed by the constraint of system~\eqref{sysadjoint}. Further, by using the Green's formula we obtain
\begin{equation*}
\begin{array} {rcl}
& & \dot{\mu}_{\tau}c(\tilde{y}_0,\tilde{y}_1,\tilde{\xi}) - 
\langle \tilde{\zeta}_0 ; \dot{\mu}_{\tau} \tilde{y}_1 \rangle_{\WW^{2,p}(\Omega)', \WW^{2,p}(\Omega)}  \\[5pt]
& &- \left\langle \tilde{\zeta}_1 ; \dot{\mu}_{\tau}
\left(\kappa \Delta \tilde{y}_1 + \divg(\sigma(\nabla \tilde{y}_0)) + f(\tilde{\xi}) \right) \right\rangle_{\LL^{p'}(\Omega), \LL^{p}(\Omega)} \\[5pt]
& = & \dot{\mu}_{\tau}\left(c(\tilde{y}_0,\tilde{y}_1,\tilde{\xi}) - 
\langle \tilde{\zeta}_0 ; \tilde{y}_1 \rangle_{\WW^{2,p}(\Omega)', \WW^{2,p}(\Omega)} 
- \langle \tilde{\zeta}_1 ; f(\tilde{\xi}) \rangle_{\LL^{p'}(\Omega), \LL^p(\Omega)} 
\right.
\\[5pt]
& & + \left. \left\langle \nabla \tilde{\zeta}_1 ; 
\kappa \nabla \tilde{y}_1 + \sigma(\nabla \tilde{y}_0)  \right\rangle_{\WWW^{1,p}(\Omega)', \WWW^{1,p}(\Omega)} \right. \\[5pt]
&  & \left.- \displaystyle \left\langle \tilde{\zeta}_1 ; \frac{\p \tilde{y}_1}{\p n} + \sigma(\nabla \tilde{y}_0) \right\rangle_{\WW^{1/(p'),p}(\Gamma_N)',\WW^{1/(p'),p}(\Gamma_N)} \right) \\[15pt]
& = & \dot{\mu}_{\tau} \mathcal{H}(\tilde{y}_0,\tilde{y}_1, \tilde{\mathfrak{p}},\tilde{\xi},\tilde{\zeta}_0,\tilde{\zeta}_1, \tilde{\pi}).
\end{array} 
\end{equation*}
Thus we obtain~\eqref{id-deriv2}, which completes the proof.
\end{proof}

\begin{remark}
Since the chosen control operator appears in~\eqref{eqHami} in the specific form $\tilde{\xi} \mapsto f(\tilde{\xi})$, the derivative~\eqref{id-deriv1} reduces to
\begin{equation*}
\tilde{J}'_{\tilde{\xi}}(\tilde{\xi},\tau) = 
-\dot{\mu}f'(\tilde{\xi})^{\ast}.\tilde{\zeta}_1.
\end{equation*}
We keep the general form~\eqref{id-deriv1} because this formula applies in a more case (see~\cite{Maxmax2}). 
\end{remark}

Then the first main result follows, namely the first-order optimality conditions for problem~\eqref{mainpbtilde2}:

\begin{theorem} \label{th-optcond}
Let be~$(\tilde{\xi}, \tau) \in \mathcal{X}_{p,2}(\omega) \times (0,T)$ an optimal solution of problem~\eqref{mainpbtilde2}. Then we have
\begin{subequations} \label{ende}
\begin{eqnarray}
c'_{\xi}(\tilde{y}_0,\tilde{y}_1,\tilde{\xi}) - f'(\tilde{\xi})^{\ast}.\tilde{\zeta}_1 & = & 0, \label{ende1}\\
\displaystyle\int_0^2\dot{\mu}_{\tau} (s,\tau)
\mathcal{H}\big(\tilde{y}_0,\tilde{y}_1,\tilde{\mathfrak{p}}, \tilde{\xi}, \tilde{\zeta}_0, \tilde{\zeta}_1,\tilde{\pi}\big)(s)\, \d s & = & 0, \label{ende2}
\end{eqnarray}
\end{subequations}
where $(\tilde{y}_0,\tilde{y}_1,\tilde{\mathfrak{p}}) = \mathbb{S}(\tilde{\xi}, \tau)$ is the solution of~\eqref{sysmaintilde2x} and $(\tilde{\zeta}_0,\tilde{\zeta}_1,\tilde{\pi})$ is the solution of~\eqref{sysadjoint} associated with $(\tilde{y}_0,\tilde{y}_1,\tilde{\mathfrak{p}})$.
\end{theorem}

\begin{proof}
Problem~\eqref{mainpbtilde2} consists in minimizing the functional $\tilde{J}$ defined in~\eqref{def-J}. From Theorem~\ref{th-cts}, the functional $\tilde{J}$ is $C^1$. Its derivatives with respect to $\tilde{\xi}$ and $\tau$ are given in Proposition~\ref{prop-optcond}. As mentioned in the prrof of the latter, we have $\mathcal{H}_{\xi}\big(\tilde{y}_0,\tilde{y}_1, \tilde{\mathfrak{p}},\tilde{\xi}, \tilde{\zeta}_0, \tilde{\zeta}_1,\tilde{\pi}\big).\hat{\xi} = c'_{\xi}(\tilde{y}_0,\tilde{y}_1,\tilde{\xi}).\hat{\xi} - f'(\hat{\xi}).\hat{\xi}$. We conclude by utilizing the Karush-Kuhn-Tucker conditions.
\end{proof}

\subsection{Other formulation of the optimality conditions} \label{sec-optcond-final}
The optimality conditions stated in Theorem~\ref{th-optcond} deal with transformed systems whose time variables is~$s\in(0,2)$. For practical purposes, like implementation, it might be more convenient to rewrite them in terms of variables satisfying systems whose time variable stands for~$t\in (0,T)$.

We state the second main result of the paper, as a corollary to Theorem~\ref{th-optcond}.
\begin{corollary} \label{cor-optcond}
Let be~$(\tilde{\xi}, \tau) \in \mathcal{X}_{p,2}(\omega) \times (0,T)$ an optimal solution of problem~\eqref{mainpbtilde2}. Then we have
\begin{subequations} \label{endebis}
\begin{eqnarray}
c'_{\xi}(u,\dot{u},\xi) - f'(\xi) & = & 0, \label{ende3}\\
\displaystyle \int_0^T \left(\dot{\mu}_{\tau}\circ \mu^{-1} (t,\tau)\right)
\mathcal{H}\big(u,\dot{u},p, \xi, \zeta_0, \zeta_1,\pi)(t)\, \d t & = & 0, \label{ende4}
\end{eqnarray}
\end{subequations}
where $(u,\dot{u},\mathfrak{p}, \xi) = (\mathbb{S}(\tilde{\xi}, \tau),\tilde{\xi})\circ \mu^{-1}(\cdot,\tau)$ satisfies~\eqref{sysmain}, and $(\zeta_0,\zeta_1,\pi)$ is solution of system~\eqref{sysadjoint-init} associated with $(y_0,y_1,\mathfrak{p}) = (u,\dot{u},\mathfrak{p})$ in the sense of Definition~\ref{def-trans-init}.
\end{corollary}

\begin{proof}
Recall from section~\ref{sec-transformation} that we have
\begin{equation}
\begin{array} {l}
\tilde{y}_0(\cdot,s) = u(\cdot,\mu(s,\tau)), \quad
\tilde{y}_1(\cdot,s) = \dot{u}(\cdot,\mu(s,\tau)), \quad
\tilde{\mathfrak{p}}(s) = \mathfrak{p}(\mu(s,\tau)),\\[5pt]
\tilde{\xi}(\cdot,s) = \xi(\cdot,\mu(s,\tau)), \quad
s\in[0,2],
\end{array}\label{change-var-proof}
\end{equation}
so that $(u,\dot{u},\mathfrak{p})$ satisfies~\eqref{sysmain} if and only if $(\tilde{y}_0,\tilde{y}_1, \tilde{\mathfrak{p}})$ satisfies~\eqref{sysmaintilde2x}. Now introduce
\begin{equation*}
\zeta_0(\cdot,t) = \tilde{\zeta}_0(\cdot,\mu^{-1}(t,\tau)), \ 
\zeta_1(\cdot,t) = \tilde{\zeta}_1(\cdot,\mu^{-1}(t,\tau)), \ 
\pi(t) = \tilde{\pi}(\mu^{-1}(t,\tau)), \quad
t\in(0,T),
\end{equation*}
where the notation~$\mu^{-1}(t,\tau)$ refers to the inverse of~$s\mapsto \mu(s,\tau)$. Composing system~\eqref{sysadjoint} (satisfied by~$(\tilde{\zeta}_0,\tilde{\zeta}_1,\tilde{\pi})$) by~$\mu^{-1}(t,\tau)$ yields system~\eqref{sysadjoint-init} (satisfied by~$(\zeta_0,\zeta_1,\pi)$) corresponding to $(y_0,y_1,\mathfrak{p}) = (u,\dot{u},\mathfrak{p})$. Then we compose~\eqref{ende1} by~$\mu^{-1}(\cdot,\tau)$ in order to obtain~\eqref{ende3}, and use the change of variable formula in the integral of~\eqref{ende2} in order to obtain~\eqref{ende4}, which concludes the proof.
\end{proof}

\begin{remark}
The term $\dot{\mu}_{\tau}\circ \mu^{-1}(t,\tau)$ which appears in~\eqref{ende4} is an Eulerian representation of the sensitivity of $\dot{\mu}$ with respect to~$\tau$. Further, choosing $\mu$ as in section~\ref{sec-transformation}, this term is actually piecewise constant:
\begin{equation*}
\dot{\mu}_\tau(\mu^{-1}(t,\tau),\tau)  =   
\left\{ \begin{array} {ll}
1  & \text{if } t\in [0,\tau), \\
0 & \text{if } t \in (\tau, \tau+\varepsilon), \\
-1/(1-\tilde{\varepsilon}) & \text{if } t\in (\tau+\varepsilon, T].
\end{array} \right.
\end{equation*}
\end{remark}

\begin{remark}
The functional $J(\tilde{\xi},\tau)$ of Problem~\eqref{mainpbtilde2}, which deals with variables $(\tilde{y}_0,\tilde{y}_1,\tilde{\xi})$, is also expressed in terms of the original variables as follows
\begin{equation}
J(\xi,\tau) =
\int_0^T c(u,\dot{u}, \xi)\, \d t + \phi^{(1)}(u,\dot{u})(\tau) + 
\phi^{(2)}(u,\dot{u})(T), \label{myfunctional}
\end{equation}
using the change of variables~\eqref{change-var-proof}.
\end{remark}

\section{Numerical \textcolor{black}{simulations}} \label{sec-num}

We propose to illustrate the optimality conditions obtained in Corollary~\ref{cor-optcond} by realizing numerical simulations which rely on finite element formulations for the space discretization. While a strong functional framework has been considered for the theoretical study of Problem~\eqref{mainpbtilde2}, the variational formulation corresponding to the finite element discretization only requires weaker regularity.

\subsection{Variational formulations} \label{sec-variational}

Solving numerically the optimality conditions for Problem~\eqref{mainpbtilde2} -- for example those provided by Corollary~\ref{cor-optcond} -- requires solving the state system~\eqref{sysmain} and the adjoint system~\eqref{sysadjoint-init}. Let us write their respective variational formulations, whose space discretizations lead to their respective finite element formulations.\\

\noindent\textbf{Weak formulation of the state system.}\\
The variational formulation of the state system~\eqref{sysmain} is given, for all test function~$\varphi \in \WW^{1,p'}(\Omega)$ such that $\varphi_{|\Gamma_D} = 0$, and for all multiplier~$\mathfrak{q} \in \R$, as follows
\begin{subequations}
\begin{eqnarray}
& & \langle \ddot{u}; \varphi\rangle_{\LL^p(\Omega),\LL^{p'}(\Omega)} 
+ \kappa \langle \nabla \dot{u}; \nabla \varphi \rangle_{\LLL^p(\Omega),\LLL^{p'}(\Omega)}
+ \displaystyle\left\langle
\frac{\p \mathcal{W}}{\p E}(E(u));(E'(u).\varphi)
\right\rangle_{\LLL^p(\Omega),\LLL^{p'}(\Omega)} \nonumber \\
& & = 
\langle f; \varphi\rangle_{\LL^p(\Omega),\LL^{p'}(\Omega)} 
+ \langle g; \varphi\rangle_{\LL^p(\Gamma_N),\LL^{p'}(\Gamma_N)} , \label{varfor1}\\
& & \mathfrak{q} \int_{\Omega} \det(\Phi(u)) \d \Omega = 
\mathfrak{q} \int_{\Omega} \det(\Phi(u_0))\d \Omega,
\end{eqnarray}
\end{subequations}
almost everywhere in~$(0,T)$. We obtained the bilinear form associated with the strain energy by using the symmetry of tensor $\displaystyle \frac{\p \mathcal{W}}{\p E}$, provided by Assumption~$\mathbf{A2}$ (we refer to Remark~\ref{remark-sym} for more details). Note that the Neumann condition~\eqref{sysmain2} is implicitly contained in~\eqref{varfor1}: By using the Green formula on~\eqref{varfor1}, we deduce both~\eqref{sysmain1} and~\eqref{sysmain2}.\\

\noindent\textbf{Weak formulation of the adjoint system.}\\
Using the Hamiltonian functional introduced in~\eqref{eqHami}, we notice that the weak formulation of the adjoint system~\eqref{sysadjoint-init} writes
\begin{equation*}
\left\{\begin{array} {rcl}
-\dot{\zeta}_0 + \mathcal{H}_{y_0}(y_0,y_1, \mathfrak{p},\xi,\zeta_0, \zeta_1,\pi) = 0 & &
 \text{in $\Omega \times\left((0,\tau) \cup (\tau, T)\right)$},\\
-\dot{\zeta}_1 + \mathcal{H}_{y_1}(y_0,y_1,\mathfrak{p},\xi,\zeta_0, \zeta_1,\pi) = 0
& & \text{in $\Omega \times\left((0,\tau) \cup (\tau, T)\right)$}, \\
\mathcal{H}_{\mathfrak{p}}(y_0,y_1,\mathfrak{p},\xi,\zeta_0, \zeta_1,\pi) = 0 & &
\text{in $\left((0,\tau) \cup (\tau, T)\right)$},
\end{array} \right.
\end{equation*}
where, using the Green formula, and denoting~$X=(y_0,y_1, \mathfrak{p},\xi,\zeta_0, \zeta_1,\pi)$, we have\begin{small}
\begin{equation*}
\begin{array} {rcl}
\mathcal{H}_{y_0}(X) .\varphi_0 & = & c'_{y_0}(y_0,y_1,\xi).\varphi_0 + 
\left\langle \nabla\zeta_1 ; 
\sigma_L(\nabla y_0).\nabla \varphi_0\right\rangle_{\WWW^{1,p}(\Omega)', \WWW^{1,p}(\Omega)} \\[5pt]
& &  +  \left\langle \varphi_0 ; \pi\, \cof(\Phi(y_0))n \right\rangle_{\LL^2(\Gamma_N)} +
\left\langle \mathfrak{p}\, \zeta_1 \, ; (\sigma_N(\nabla y_0).\nabla \varphi_0)n\right\rangle_{\LL^2(\Gamma_N)} \\[5pt]
& = & c'_{y_0}(y_0,y_1,\xi).\varphi_0
- \left\langle \divg (\sigma_L(\nabla y_0)^{\ast}. \nabla \zeta_1) ;
 \varphi_0\right\rangle_{\LL^{p'}(\Omega), \LL^{p}(\Omega)}  \\[5pt]
& & +\left\langle \varphi_0 ;  ((\sigma_L+\mathfrak{p}\, \sigma_N)(\nabla y_0)^{\ast}. \nabla \zeta_1)n + \pi\, \cof(\Phi(y_0))n \right\rangle_{\LL^2(\Gamma_N)}  , \\[5pt]
\mathcal{H}_{y_1}(X) .\varphi_1 
& = & 
c'_{y_1}(y_0,y_1,\xi).\varphi_1 - 
\langle \zeta_0 \, ; \varphi_1 \rangle_{\WW^{2,p}(\Omega)', \WW^{2,p}(\Omega)}
+ \kappa\left\langle \nabla\zeta_1 ; \nabla \varphi_1 \right\rangle_{\WWW^{1,p}(\Omega)', \WWW^{1,p}(\Omega)}  \\[5pt]
& = & c'_{y_1}(y_0,y_1,\xi).\varphi_1 - 
\langle \zeta_0 + \kappa \Delta \zeta_1 \, ; \varphi_1 \rangle_{\WW^{2,p}(\Omega)', \WW^{2,p}(\Omega)}\\[5pt]
& &  + \left\langle  \displaystyle\kappa\frac{\p \zeta_1}{\p n}\, ; \varphi_1 \right\rangle_{\WW^{2-1/p,p}(\Gamma_N)',\WW^{2-1/p,p}(\Gamma_N)}, \\[10pt]
\mathcal{H}_{\mathfrak{p}}(X)
& = & \left\langle \zeta_1 \, ; \cof(\Phi(y_0))n\right\rangle_{\WW^{1/(p'),p}(\Gamma_N)',\WW^{1/(p'),p}(\Gamma_N)}.
\end{array}
\end{equation*}\end{small}
Actually it has been shown in~\cite{Maxmax1} that the variational formulation for the adjoint system~\eqref{sysadjoint-init} can be derived as from the Hamiltonian functional as above. Note that the adjoint system is solved backward in time. It remains to comment on the initial values of system~\eqref{sysadjoint-init}.\\

\noindent\textbf{Expressions of the initial value conditions for the adjoint system.}\\
As mentioned in the introduction, we aim to maximize the variations of the pressure, namely
\begin{equation*}
\phi^{(1)}(u,\dot{u})(\tau) = (\mathfrak{p}(\tau+\varepsilon) - \mathfrak{p}(\tau))/\varepsilon,
\end{equation*}
where~$\mathfrak{p}$ is a function of~$(u,\dot{u})$, as from~\eqref{sysmain2} we obtain
\begin{equation*}
\mathfrak{p} = -\frac{1}{|\Gamma_N|}\int_{\Gamma_N} 
(\det(\Phi(u)))^{-1} \Phi(u)^T \left(
\kappa \frac{\p \dot{u}}{\p n} + \sigma(\nabla u) n
\right)
\d \Gamma_N.
\end{equation*}
The sensitivity of~$\mathfrak{p}$ with respect to~$u$ and~$\dot{u}$ is given in the variational sense by\begin{small}
\begin{equation*}
\begin{array} {rcl}
\displaystyle 
\frac{\p \mathfrak{p}}{\p u}.v & = & 
\displaystyle
-\frac{1}{|\Gamma_N|}\int_{\Gamma_N} 
(\det(\Phi(u)))^{-1} 
\Big(\nabla v^T - \big(\Phi(u)^{-T}:\nabla v\big)\Phi(u)^T\Big)
\left(\kappa \frac{\p \dot{u}}{\p n} + \sigma(\nabla u)n\right) \d \Gamma_N \\[10pt]
& & \displaystyle 
-\frac{1}{|\Gamma_N|}\int_{\Gamma_N} 
(\det(\Phi(u)))^{-1} \Phi(u)^T\big(\sigma_L(\nabla u).\nabla v\big)n \d \Gamma_N,\\[10pt]
\displaystyle\frac{\p \mathfrak{p}}{\p \dot{u}}.\dot{v} & = &
\displaystyle
-\frac{\kappa}{|\Gamma_N|}\int_{\Gamma_N} 
(\det(\Phi(u)))^{-1} \Phi(u)^T 
 \frac{\p \dot{v}}{\p n} \d \Gamma_N,
\end{array}
\end{equation*}\end{small}
and consequently the first-order derivatives of functional~$\phi^{(1)}$ are expressed as
\begin{equation*}
\phi^{(1)'}_u(u,\dot{u}) = \frac{1}{\varepsilon}\left(
\frac{\p \mathfrak{p}}{\p u}(\tau+\varepsilon) 
- \frac{\p \mathfrak{p}}{\p u}(\tau)
\right), \quad
\phi^{(1)'}_{\dot{u}}(u,\dot{u}) = \frac{1}{\varepsilon}\left(
\frac{\p \mathfrak{p}}{\p \dot{u}}(\tau+\varepsilon) 
- \frac{\p \mathfrak{p}}{\p \dot{u}}(\tau)\right).
\end{equation*}
These expressions are needed when implementing the numerical solution of the adjoint system~\eqref{sysadjoint-init}.\\

\noindent\textbf{Implementation of the jump conditions for the adjoint state.}\\
In practice we consider a subdivision of the time interval~$(0,2)$ when addressing the optimality conditions given in Theorem~\ref{th-optcond}, or of the time interval~$(0,T)$ when addressing those given in Corollary~\ref{cor-optcond}. The adjoint system~\eqref{sysadjoint} deals with jump conditions at fixed time~$s=1$, and therefore if the time subdivision of the interval~$(0,2)$ meets~$s=1$, then implementing the jump conditions at~$s=1$ does not present any difficulties. On the other hand, the adjoint system~\eqref{sysadjoint-init} involves jumps at free time $t=\tau$, which in general does not coincide with a point of the time subdivision. Therefore we need to approximate the values of the right-and-sides of these jump conditions at time~$\tau$. For example, if $\tau \in [t_i,t_{i+1}]$, where the $\{t_i\}_{i\in I}$ define the time subdivision,  one can use the following linear approximation:
\begin{equation*}
\phi^{(1)}(y_0,y_1)(\tau) \approx \displaystyle
\frac{t_{i+1}-\tau}{t_{i+1}-t_i}\phi^{(1)}(y_0,y_1)(t_i) +
\frac{\tau-t_i}{t_{i+1}-t_i}\phi^{(1)}(y_0,y_1)(t_{i+1}).
\end{equation*}
Such an approximation introduces an error of order~$1$ in the time scheme. In the numerical realizations presented in section~\ref{sec-results}, we chose to deal with the transformed systems~\eqref{sysmaintilde} and~\eqref{sysadjoint}, and thus with Problem~\ref{mainpbtilde}.

\subsection{Algorithm} \label{sec-algorithm}
We adopt the so-called {\it optimize-discretize} approach, meaning that we discretize the optimality conditions initially obtained in Corollary~\ref{cor-optcond} for the continuous problem. The other approach would consist in first discretizing Problem~\eqref{mainpb}, and next deriving optimality conditions for the corresponding discretized problem (which would be the focus of another approach). Note that the optimality conditions so derived would be specific to the discretization chosen for the state system and the objective functional.

The optimality conditions obtained in Corollary~\ref{cor-optcond} provide a gradient to vanish, namely
\begin{equation}
\mathcal{G}(\xi, \tau) := \left(\begin{array} {c}
c'_{\xi}(u,\dot{u},\xi) - f'(\xi) \\
\displaystyle \int_0^T \left(\dot{\mu}_{\tau}\circ \mu^{-1} (t,\tau)\right)
\mathcal{H}\big(u,\dot{u},\mathfrak{p}, \xi, \zeta_0, \zeta_1,\pi)(t)\, \d t 
\end{array} \right), \label{mygradient}
\end{equation}
where $(u,\dot{u}, \mathfrak{p}) = (y_0,y_1,\mathfrak{p})$ satisfies~\eqref{sysmain} and $(\zeta_0,\zeta_1,\pi)$ satisfies~\eqref{sysadjoint-init}. We solve Problem~\eqref{mainpb} with a gradient rule, more specifically the Barzilai-Borwein algorithm~\cite{BarBor}. The corresponding method is given in Algorithm~\ref{algo-super} below.

\begin{algorithm}[htpb]
\hfill \begin{small}
	\begin{description}
		\item[Initialization:] \hfill
		\begin{description}
			\item[-] Initialize $(\xi_{0}, \tau_0)= (0,T/2)$.
		\end{description}

		\item[Initial gradient:] From $(\xi_0, \tau_0)$, compute the (initial) gradient as follows: \\[5pt]
		\begin{tabular} {lr}
		$\left.\begin{minipage}{10.5cm}
		\begin{description}
			\item[-] Compute the state $(u,\dot{u},\mathfrak{p})$ corresponding to $\xi_0$, by solving~\eqref{sysmain}.
			\item[-] Compute the adjoint state $(\zeta_0,\zeta_1,\pi)$ corresponding to $(\xi_0,\tau_0)$,\\by solving~\eqref{sysadjoint-init}.
			\item[-] Compute the gradient $\mathcal{G}_0:= \mathcal{G}(\xi_0,\tau_0)$, using the expression~\eqref{mygradient}.
		\end{description}
		\end{minipage} 
		\right\}$ & 
		\begin{minipage}{4cm}
		Compute \\the gradient \\ in 3 steps
		\end{minipage}
		\end{tabular}
		\hfill \\Store $(\xi_0, \tau_0)$ and $\mathcal{G}_0$.

		\item[Armijo rule:] Choose $\alpha = 0.5$. 
		\begin{description}
			\item[-] Find the smallest $n\in\N$ such that $J((\xi_0,\tau_0)-\alpha^n \mathcal{G}_0) < J(\xi_0,\tau_0)$, using expression~\eqref{myfunctional}.
			\item[-] Define $(\xi_1,\tau_1) = (\xi_0,\tau_0) - \alpha^n G_0$. 
			\item[-] Compute the gradient $\mathcal{G}_1$ as above, corresponding to $(\xi_1,\tau_1)$.
			\item[-] Store $(\xi_1,\tau_1)$ and $\mathcal{G}_1$.
		\end{description}

	\item[Barzilai-Borwein gradient steps:] Initialization with $((\xi_0,\tau_0),\mathcal{G}_0)$ and $((\xi_1,\tau_1),\mathcal{G}_1)$.\\[5pt]
	Compute iteratively $(\xi_n,\tau_n)$ ($n\geq 2$) with the Barzilai-Borwein steps.\\
	While $|| \mathcal{G}(\xi_n,\tau_n)||_{\LL^2(\omega) \times \R} > 1.e^{-10}$, do gradient steps.

	\item[End:]
	Obtain $(\xi,\tau)$, approximated solution of~\eqref{mainpbtilde2} with.
	\end{description}
	\caption{Solving Problem~\eqref{mainpb} via the first-order necessary optimality conditions of Corollary~\ref{cor-optcond}.}\label{algo-super}\end{small}
\end{algorithm}
\FloatBarrier 

\vspace*{-30pt}
\textcolor{black}{
\subsection{Implementation}} \label{sec-results}

Consider a 1D toy-model for which $\Omega = (0,1)$, $\Gamma_D = \{0\}$ and $\Gamma_N = \{1\}$. \textcolor{black}{The goal of this subsection is to provide an illustration of the implementation of a solution for Problem~\eqref{mainpb} in order to illustrate the qualitative aspects of the problem. For this purpose, since we are mainly interested in the time evolution of the scalar-valued variable~$\mathfrak{p}$, it is unnecessary to provide illustrations in higher dimensions, as the latter would require non-trivial numerical methods that we do not aim to address here. Note that the theoretical analysis provided in this article may {\it a priori} need to be adapted to the 1D case}. For the strain energy, we propose to  consider the example Saint~Venant-Kirchhoff model (see section~\ref{sec-app-energies}), namely
\begin{equation}
\mathcal{W}(E) = \mu_L \mathrm{tr}(E^2) + \frac{\lambda_L}{2} \mathrm{tr}(E)^2,
\quad
\displaystyle \frac{\p \mathcal{W}}{\p E}(E) = 2\mu_L E + \lambda_L \mathrm{tr}(E) \I,
\label{Lame-model}
\end{equation}
where $\mu_L$ and~$\lambda_L$ are given Lam\'e coefficients. We refer to section~\ref{sec-app-energies} for further details on this model. As mentioned previously the space discretization is realized with finite elements, more specifically P1-elements. The control is distributed on the subdomain $\omega = [0.75,1.00]$. The cost function and the objective functions are chosen to be  
\begin{equation*}
c(u,\dot{u},\xi) = -\displaystyle\frac{\alpha}{2}\|\xi\|_{\L^2(\omega)}^2, \quad
\phi^{(1)}(u,\dot{u}) = \mathfrak{p}, \quad
\phi^{(2)}(u,\dot{u}) = 0.
\end{equation*}
with $\alpha >0$. Regarding the choice of the functional~$\phi^{(1)}$, unlike~\eqref{sysmainpressure} where we aim at maximizing the variations of the pressure at some time~$\tau \in (0,T)$, we rather aim at maximizing the pressure itself directly, because initially the pressure is equal to zero, in view of the choice we made for the initial conditions in Table~\ref{tab-param}. The time discretization for the state system~\eqref{sysmain} corresponds to the Crank-Nicolson method (that is the $\theta$-method with~$\theta = 0.5$), while the time discretization for the adjoint system~\eqref{sysadjoint-init} is an implicit Euler scheme. At each time step, the nonlinearity due to the strain energy terms are treated with the Newton method. The choices for the different parameters are summarized in Table~\ref{tab-param}.

\begin{table}[H]
\begin{center}
\begin{equation*}
\begin{array} {|c|c|c|c|c|c|c|c|c|c|}
	\hline
	\alpha & \kappa & \lambda_L & \mu_L & g & u_0 & \dot{u}_0 & T & \text{time step} & \text{mesh size}   \\
	\hline
	2.10^{-3} & 2.10^{-4}  & 0.05 & 0.05 & 0 & 0 & 0 & 15.0 & 0.02 & 0.01\\
	\hline
\end{array}
\end{equation*}
\vspace*{-10pt}
\caption{Values of parameters for numerical realization}
\label{tab-param}
\end{center}
\end{table}
\FloatBarrier 

The implementation has been realized in C++, using the Getfem++ Library~\cite{Getfem}.

\textcolor{black}{
\subsection{Results}}
Using Algorithm~\ref{algo-super}, we obtain that the optimal time parameter is approximately~$\tau \approx 7.9$. We provide screenshots representing the time evolution of the different variables in Figures~\ref{figControl},~\ref{figDisp} and~\ref{figVelo}. The time evolution of the pressure is represented in Figure~\ref{fig-pressure}. Note that in view of the parameters given in Table~\ref{tab-param}, with no control we would obtain the trivial states~$u=0$ and~$\dot{u} =0$ on~$(0,T)$.

\hspace*{-15pt}\begin{minipage}{\linewidth}
\centering
\begin{minipage}{0.24\linewidth}
\includegraphics[trim = 0cm 0cm 0cm 1cm, clip, scale=0.25]{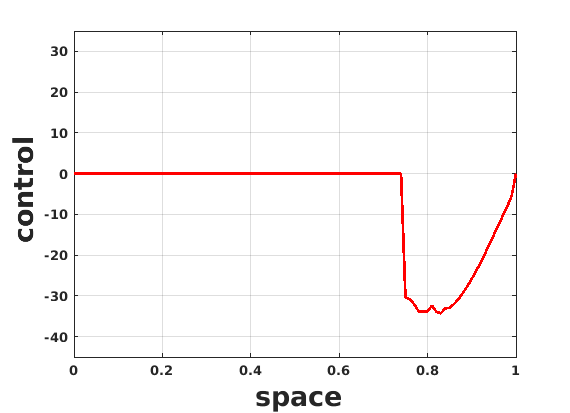}
\begin{center} $ t = 0.02 $ \end{center}
\end{minipage}
\begin{minipage}{0.24\linewidth}
\includegraphics[trim = 0cm 0cm 0cm 1cm, clip, scale=0.25]{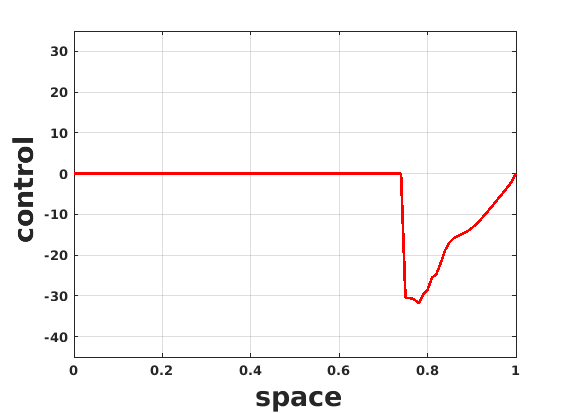}
\begin{center} $ t = 1.00 $ \end{center}
\end{minipage}
\begin{minipage}{0.24\linewidth}
\includegraphics[trim = 0cm 0cm 0cm 1cm, clip, scale=0.25]{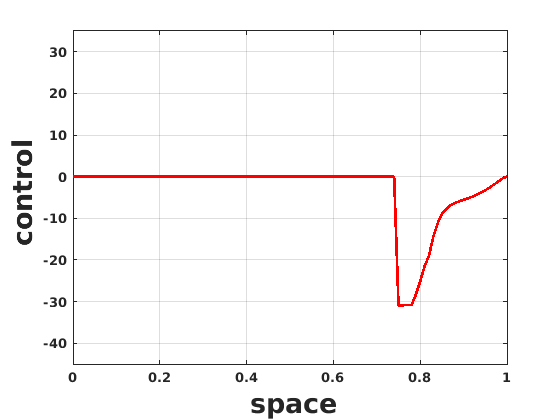}
\begin{center} $ t = 1.5 $ \end{center}
\end{minipage}
\begin{minipage}{0.24\linewidth}
\includegraphics[trim = 0cm 0cm 0cm 1cm, clip, scale=0.25]{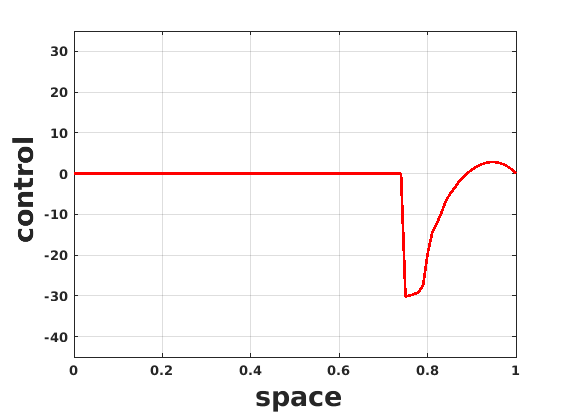}
\begin{center} $ t = 2.0 $ \end{center}
\end{minipage}
\\ 
\begin{minipage}{0.24\linewidth}
\includegraphics[trim = 0cm 0cm 0cm 1cm, clip, scale=0.25]{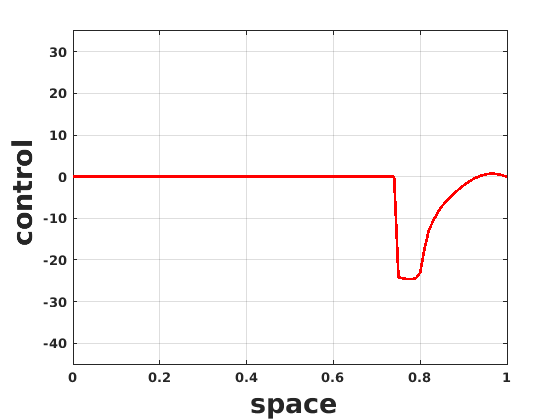}
\begin{center} $ t = 2.5 $ \end{center}
\end{minipage}
\begin{minipage}{0.24\linewidth}
\includegraphics[trim = 0cm 0cm 0cm 1cm, clip, scale=0.25]{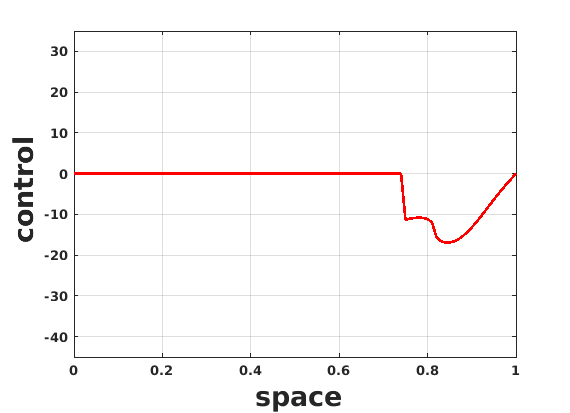}
\begin{center} $ t = 3.0 $ \end{center}
\end{minipage}
\begin{minipage}{0.24\linewidth}
\includegraphics[trim = 0cm 0cm 0cm 1cm, clip, scale=0.25]{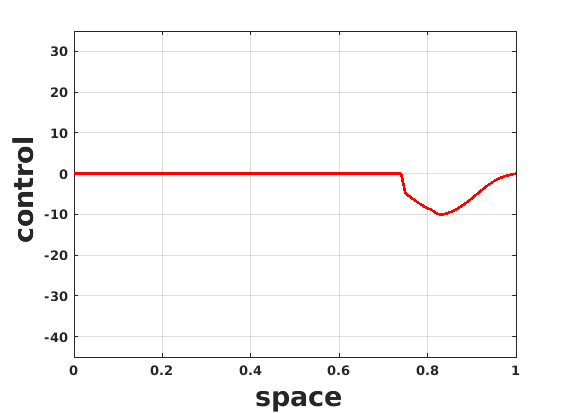}
\begin{center} $ t = 3.5 $ \end{center}
\end{minipage}
\begin{minipage}{0.24\linewidth}
\includegraphics[trim = 0cm 0cm 0cm 1cm, clip, scale=0.25]{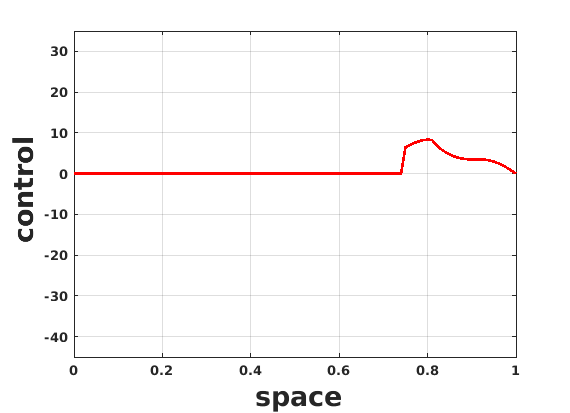}
\begin{center} $ t = 4.0 $ \end{center}
\end{minipage}
\\
\begin{minipage}{0.24\linewidth}
\includegraphics[trim = 0cm 0cm 0cm 1cm, clip, scale=0.25]{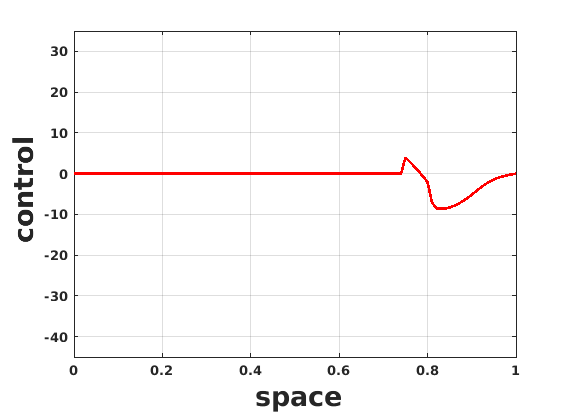}
\begin{center} $ t = 4.5 $ \end{center}
\end{minipage}
\begin{minipage}{0.24\linewidth}
\includegraphics[trim = 0cm 0cm 0cm 1cm, clip, scale=0.25]{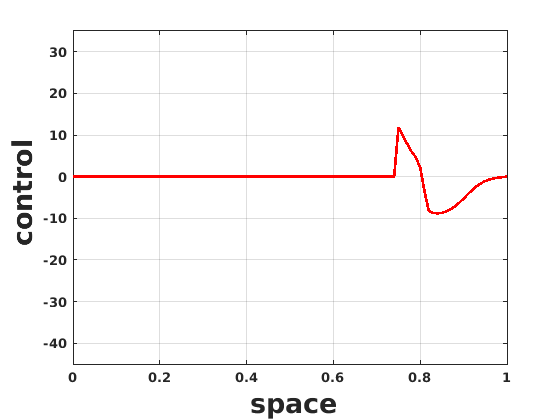}
\begin{center} $ t = 5.0 $ \end{center}
\end{minipage}
\begin{minipage}{0.24\linewidth}
\includegraphics[trim = 0cm 0cm 0cm 1cm, clip, scale=0.25]{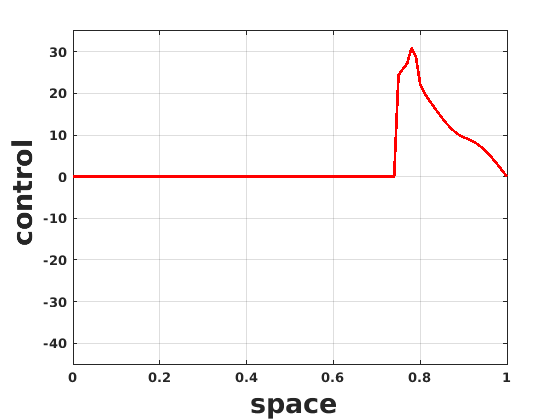}
\begin{center} $ t = 5.5 $ \end{center}
\end{minipage}
\begin{minipage}{0.24\linewidth}
\includegraphics[trim = 0cm 0cm 0cm 1cm, clip, scale=0.25]{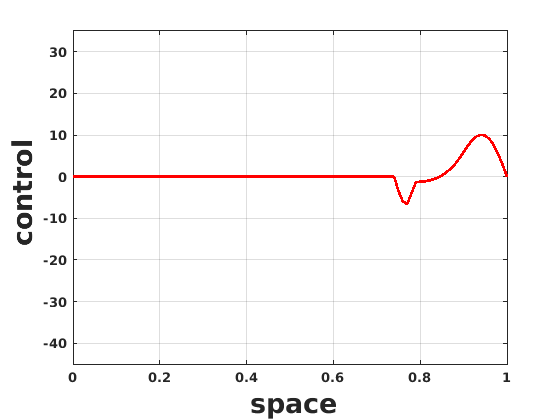}
\begin{center} $ t = 6.0 $ \end{center}
\end{minipage}
\\
\begin{minipage}{0.24\linewidth}
\includegraphics[trim = 0cm 0cm 0cm 1cm, clip, scale=0.25]{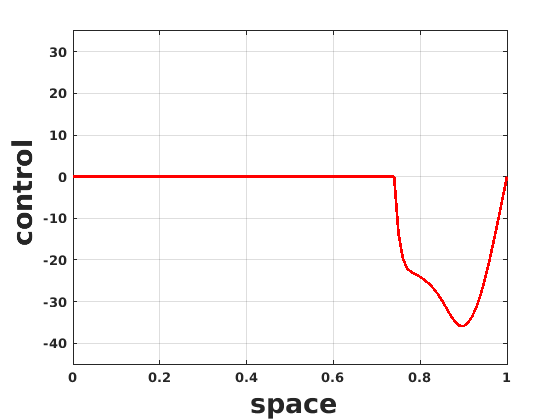}
\begin{center} $ t = 6.5 $ \end{center}
\end{minipage}
\begin{minipage}{0.24\linewidth}
\includegraphics[trim = 0cm 0cm 0cm 1cm, clip, scale=0.25]{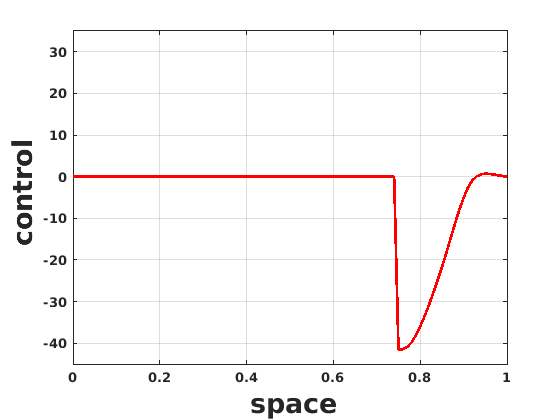}
\begin{center} $ t = 7.0 $ \end{center}
\end{minipage}
\begin{minipage}{0.24\linewidth}
\includegraphics[trim = 0cm 0cm 0cm 1cm, clip, scale=0.25]{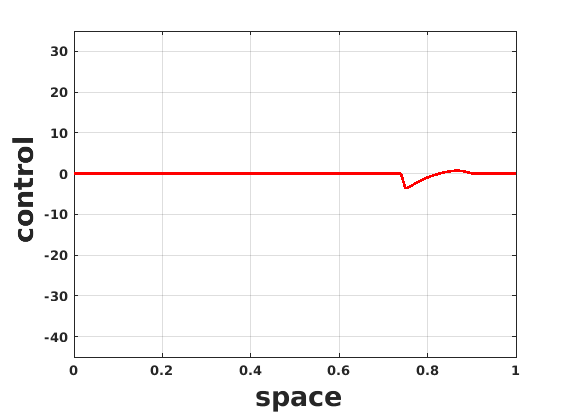}
\begin{center} $ t = 7.46 $ \end{center}
\end{minipage}
\begin{minipage}{0.24\linewidth}
\includegraphics[trim = 0cm 0cm 0cm 1cm, clip, scale=0.25]{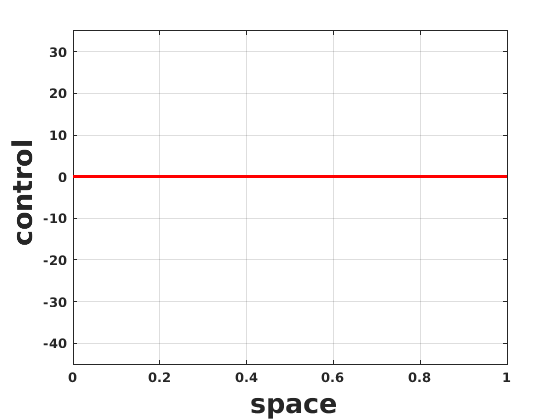}
\begin{center} $ t = 7.48 $ \end{center}
\end{minipage}
\end{minipage}
\begin{minipage}{\linewidth}
\begin{figure}[H]
\caption{Values of the control function, for \textcolor{black}{different times}. }
\label{figControl}
\end{figure}
\end{minipage}\\
\hfill \\
\FloatBarrier

In Figure~\ref{figControl} we observe that the control function, distributed on the segment~$[0.75,1.00]$, is sparse in time, in the sense that it becomes inactivated from~$t=7.48$, a short time before the optimal time parameter~$t=\tau \approx 7.9$. The fact that it is sparse in time is clearly explained by the fact that the terminal objective functional~$\phi^{(2)}$ is chosen to be equal to zero. The fact that the control function gets inactive a short time before~$\tau$ could be explained by the propagation effect that makes the control useless during a (short) time interval before~$\tau$. Note also that the sign of the control function changes rapidly, that we could explain by the necessity of to creating a wave phenomenon (see Figure~\ref{figVelo}).


\hspace*{-15pt}\begin{minipage}{\linewidth}
\centering
\begin{minipage}{0.24\linewidth}
\includegraphics[trim = 0cm 0cm 0cm 1cm, clip, scale=0.25]{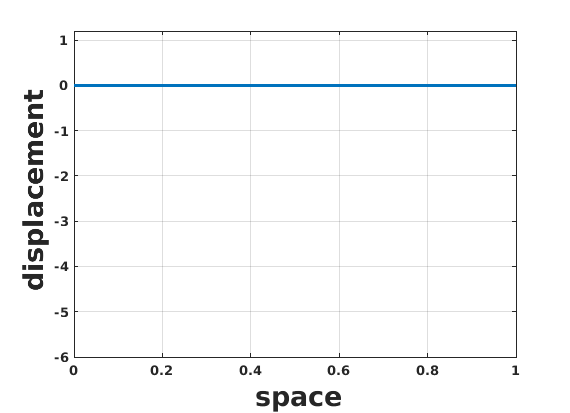}
\begin{center} $ t = 0.02 $ \end{center}
\end{minipage}
\begin{minipage}{0.24\linewidth}
\includegraphics[trim = 0cm 0cm 0cm 1cm, clip, scale=0.25]{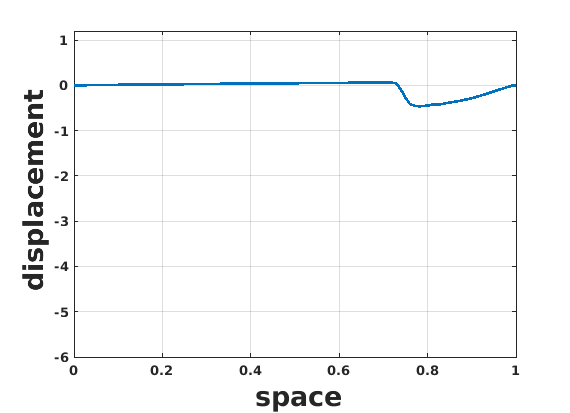}
\begin{center} $ t = 1.00 $ \end{center}
\end{minipage}
\begin{minipage}{0.24\linewidth}
\includegraphics[trim = 0cm 0cm 0cm 1cm, clip, scale=0.25]{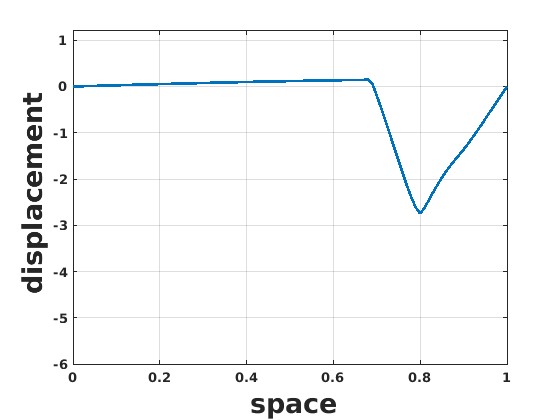}
\begin{center} $ t = 2.5 $ \end{center}
\end{minipage}
\begin{minipage}{0.24\linewidth}
\includegraphics[trim = 0cm 0cm 0cm 1cm, clip, scale=0.25]{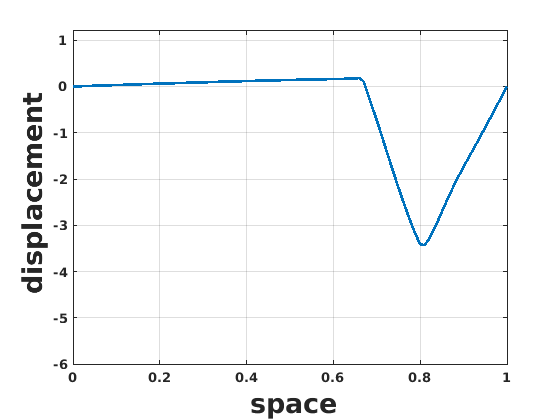}
\begin{center} $ t = 3.0 $ \end{center}
\end{minipage}
\\
\begin{minipage}{0.24\linewidth}
\includegraphics[trim = 0cm 0cm 0cm 1cm, clip, scale=0.25]{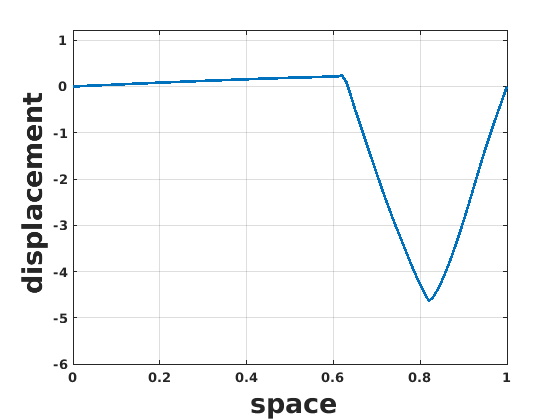}
\begin{center} $ t = 4.0 $ \end{center}
\end{minipage}
\begin{minipage}{0.24\linewidth}
\includegraphics[trim = 0cm 0cm 0cm 1cm, clip, scale=0.25]{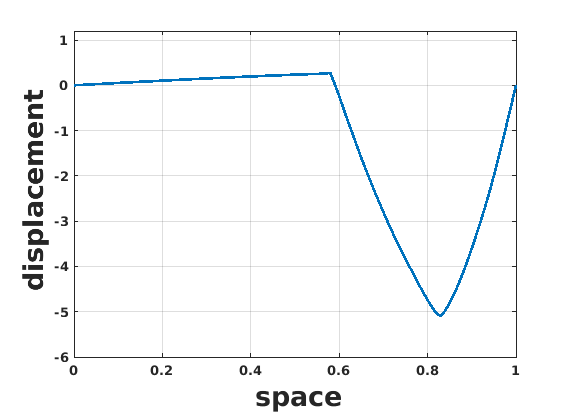}
\begin{center} $ t = 5.0 $ \end{center}
\end{minipage}
\begin{minipage}{0.24\linewidth}
\includegraphics[trim = 0cm 0cm 0cm 1cm, clip, scale=0.25]{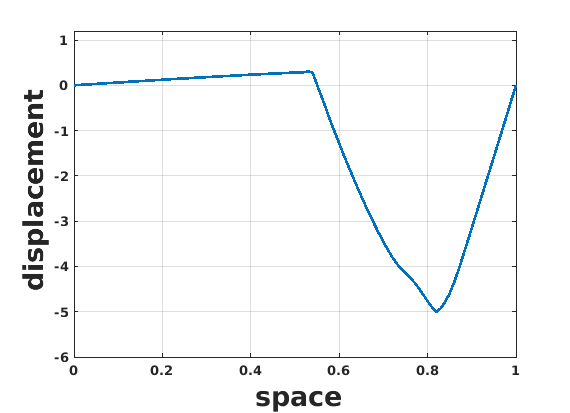}
\begin{center} $ t = 6.0 $ \end{center}
\end{minipage}
\begin{minipage}{0.24\linewidth}
\includegraphics[trim = 0cm 0cm 0cm 1cm, clip, scale=0.25]{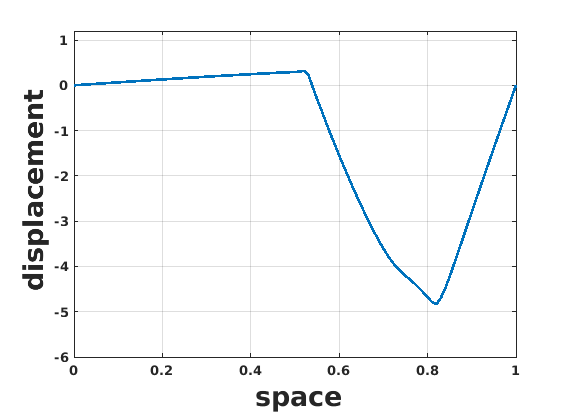}
\begin{center} $ t = 6.3 $ \end{center}
\end{minipage}
\\
\begin{minipage}{0.24\linewidth}
\includegraphics[trim = 0cm 0cm 0cm 1cm, clip, scale=0.25]{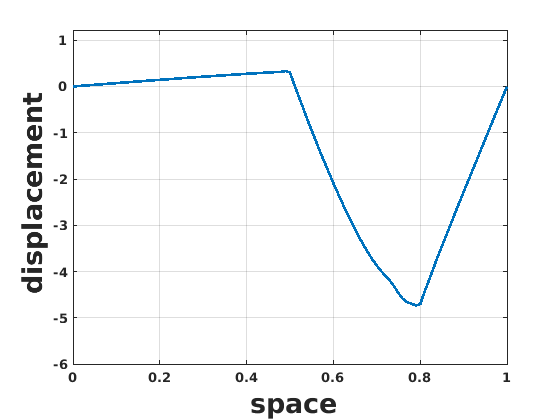}
\begin{center} $ t = 7.0 $ \end{center}
\end{minipage}
\begin{minipage}{0.24\linewidth}
\includegraphics[trim = 0cm 0cm 0cm 1cm, clip, scale=0.25]{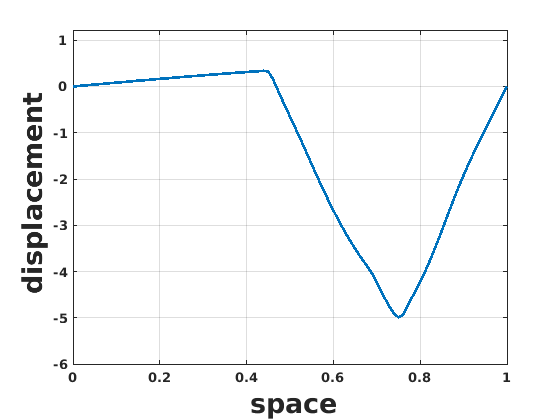}
\begin{center} $ t = 8.0 $ \end{center}
\end{minipage}
\begin{minipage}{0.24\linewidth}
\includegraphics[trim = 0cm 0cm 0cm 1cm, clip, scale=0.25]{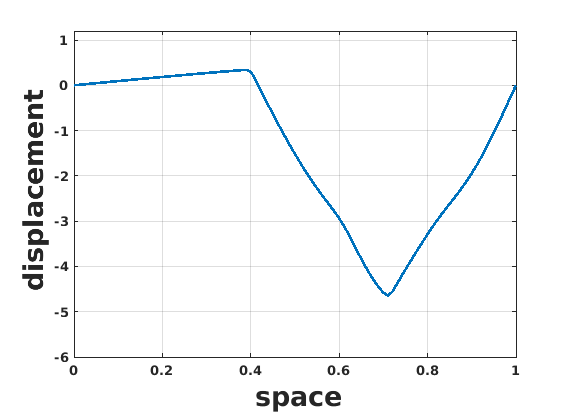}
\begin{center} $ t = 9.0 $ \end{center}
\end{minipage}
\begin{minipage}{0.24\linewidth}
\includegraphics[trim = 0cm 0cm 0cm 1cm, clip, scale=0.25]{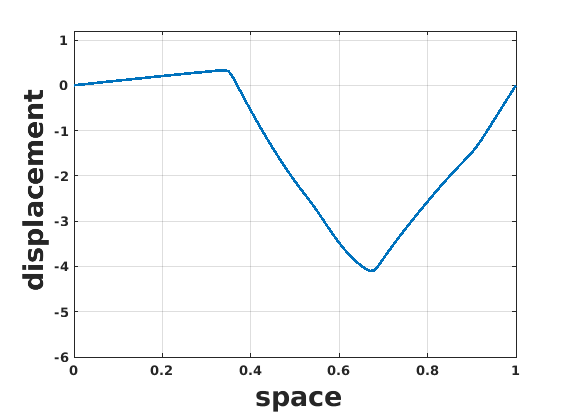}
\begin{center} $ t = 10.0 $ \end{center}
\end{minipage}
\\
\begin{minipage}{0.24\linewidth}
\includegraphics[trim = 0cm 0cm 0cm 1cm, clip, scale=0.25]{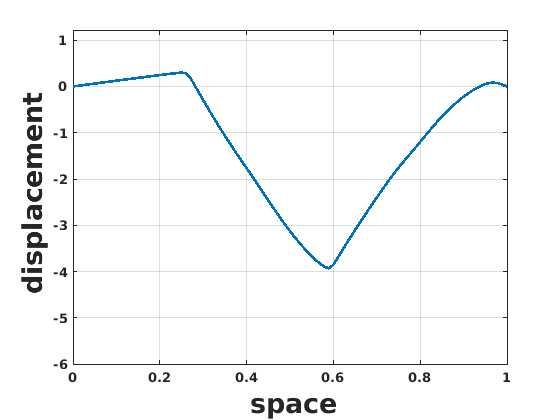}
\begin{center} $ t = 12.0 $ \end{center}
\end{minipage}
\begin{minipage}{0.24\linewidth}
\includegraphics[trim = 0cm 0cm 0cm 1cm, clip, scale=0.25]{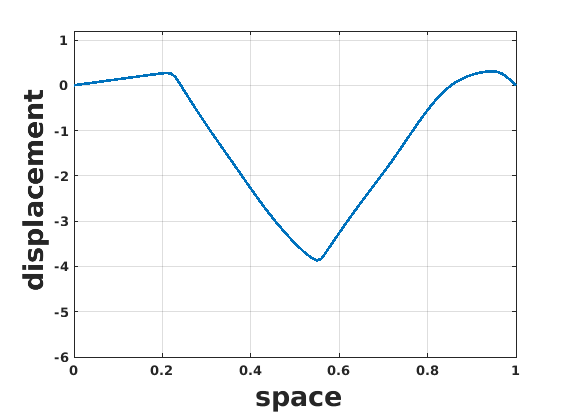}
\begin{center} $ t = 13.0 $ \end{center}
\end{minipage}
\begin{minipage}{0.24\linewidth}
\includegraphics[trim = 0cm 0cm 0cm 1cm, clip, scale=0.25]{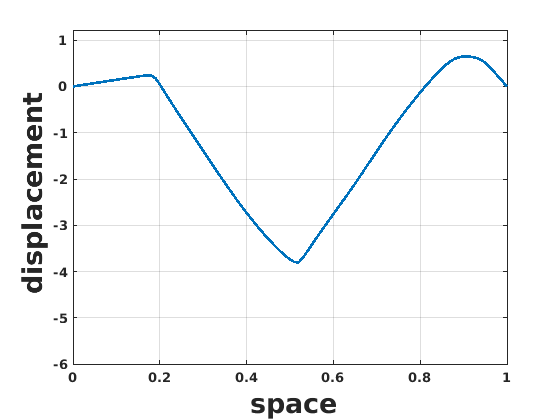}
\begin{center} $ t = 14.0 $ \end{center}
\end{minipage}
\begin{minipage}{0.24\linewidth}
\includegraphics[trim = 0cm 0cm 0cm 1cm, clip, scale=0.25]{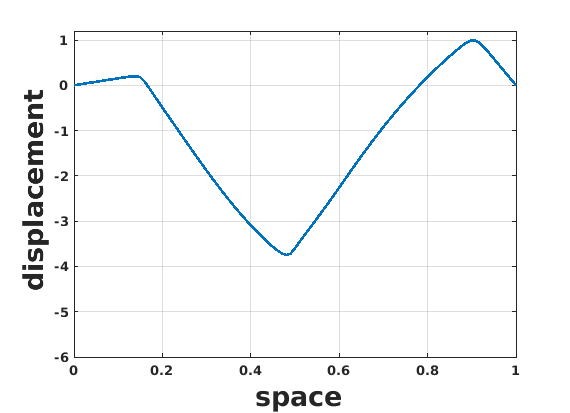}
\begin{center} $ t = 15.0 $ \end{center}
\end{minipage}
\end{minipage}
\begin{minipage}{\linewidth}
\begin{figure}[H]
\caption{Values of the displacement variable~$u$, for for \textcolor{black}{different times}. }
\label{figDisp}
\end{figure}
\end{minipage}\\
\hfill \\
\FloatBarrier

In Figure~\ref{figDisp}, where the time evolution of the displacement field is represented, we observe that the state~$u$ tends to become steady, and then define a deformed domain~$(\Id+u)(\Omega)$ of steady shape.

\newpage

\hspace*{-15pt}\begin{minipage}{\linewidth}
\centering
\begin{minipage}{0.24\linewidth}
\includegraphics[trim = 0cm 0cm 0cm 1cm, clip, scale=0.25]{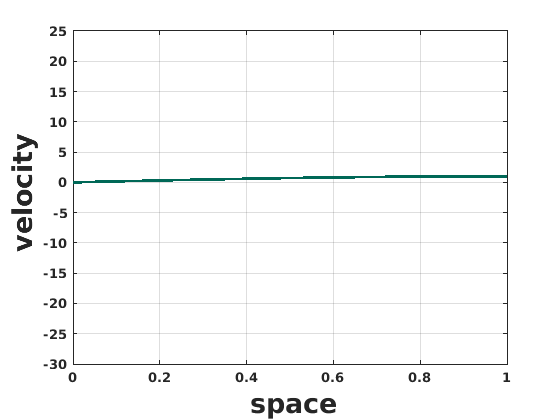}
\begin{center} $ t = 0.02 $ \end{center}
\end{minipage}
\begin{minipage}{0.24\linewidth}
\includegraphics[trim = 0cm 0cm 0cm 1cm, clip, scale=0.25]{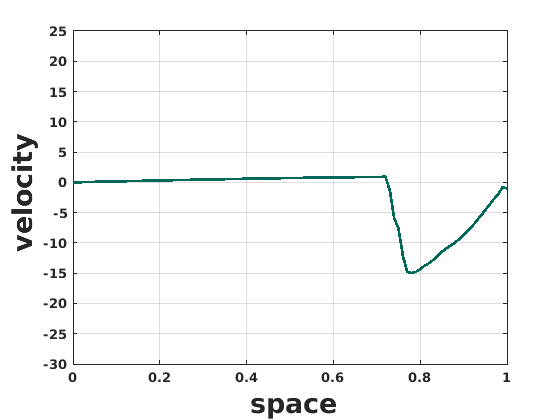}
\begin{center} $ t = 1.00 $ \end{center}
\end{minipage}
\begin{minipage}{0.24\linewidth}
\includegraphics[trim = 0cm 0cm 0cm 1cm, clip, scale=0.25]{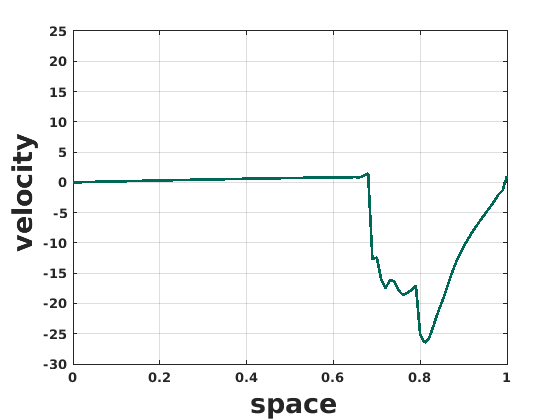}
\begin{center} $ t = 2.5 $ \end{center}
\end{minipage}
\begin{minipage}{0.24\linewidth}
\includegraphics[trim = 0cm 0cm 0cm 1cm, clip, scale=0.25]{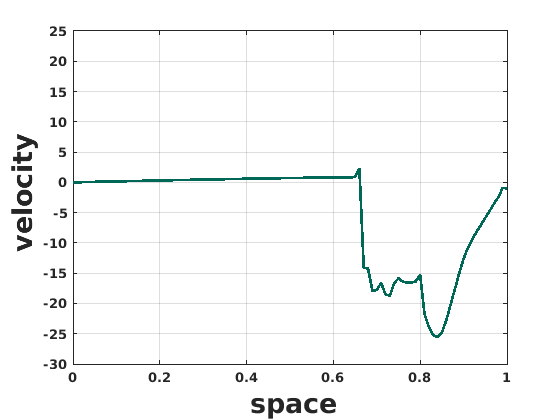}
\begin{center} $ t = 3.0 $ \end{center}
\end{minipage}
\\
\begin{minipage}{0.24\linewidth}
\includegraphics[trim = 0cm 0cm 0cm 1cm, clip, scale=0.25]{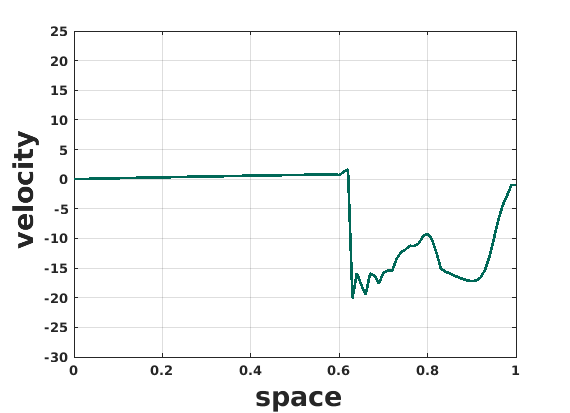}
\begin{center} $ t = 4.0 $ \end{center}
\end{minipage}
\begin{minipage}{0.24\linewidth}
\includegraphics[trim = 0cm 0cm 0cm 1cm, clip, scale=0.25]{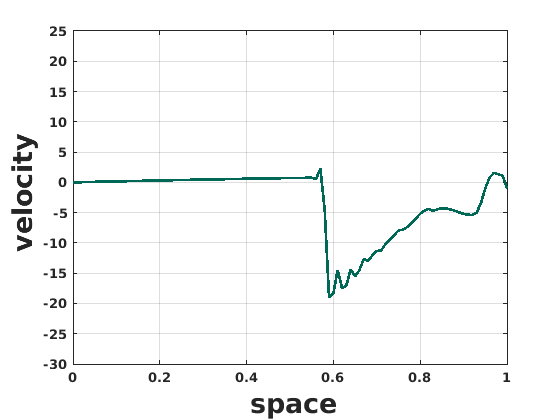}
\begin{center} $ t = 5.0 $ \end{center}
\end{minipage}
\begin{minipage}{0.24\linewidth}
\includegraphics[trim = 0cm 0cm 0cm 1cm, clip, scale=0.25]{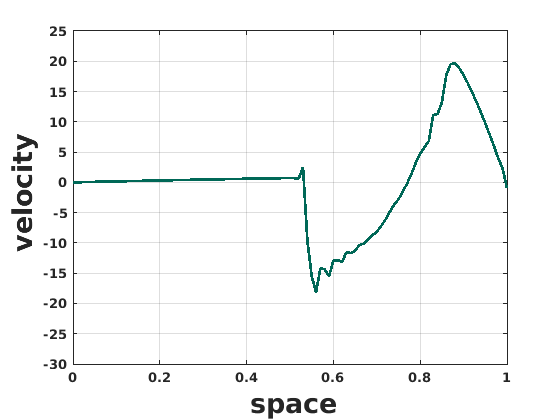}
\begin{center} $ t = 6.0 $ \end{center}
\end{minipage}
\begin{minipage}{0.24\linewidth}
\includegraphics[trim = 0cm 0cm 0cm 1cm, clip, scale=0.25]{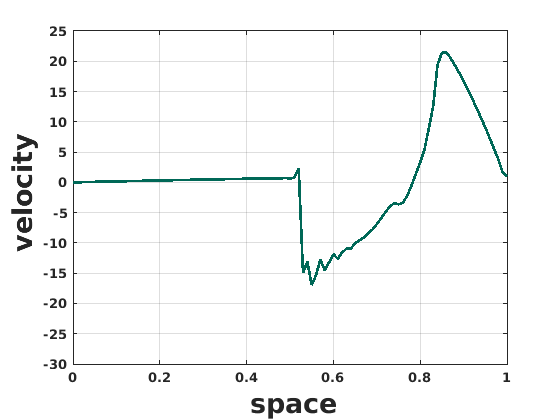}
\begin{center} $ t = 6.3 $ \end{center}
\end{minipage}
\\
\begin{minipage}{0.24\linewidth}
\includegraphics[trim = 0cm 0cm 0cm 1cm, clip, scale=0.25]{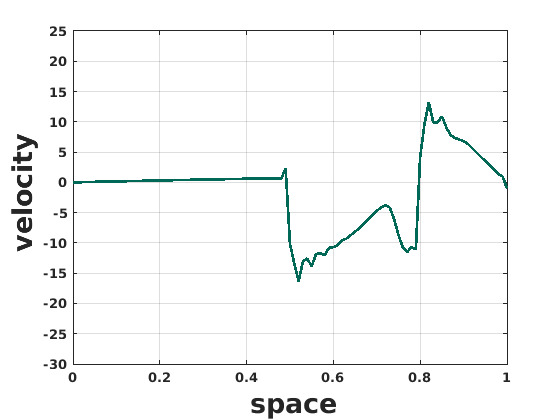}
\begin{center} $ t = 7.0 $ \end{center}
\end{minipage}
\begin{minipage}{0.24\linewidth}
\includegraphics[trim = 0cm 0cm 0cm 1cm, clip, scale=0.25]{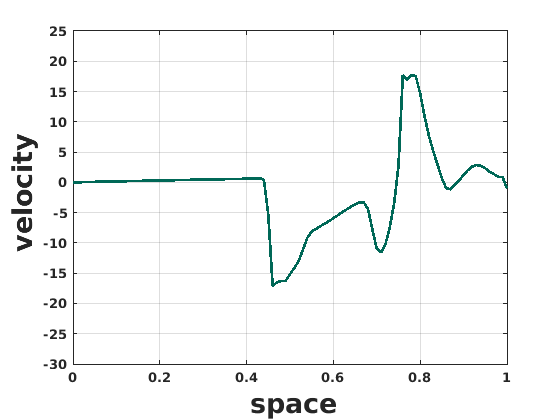}
\begin{center} $ t = 8.0 $ \end{center}
\end{minipage}
\begin{minipage}{0.24\linewidth}
\includegraphics[trim = 0cm 0cm 0cm 1cm, clip, scale=0.25]{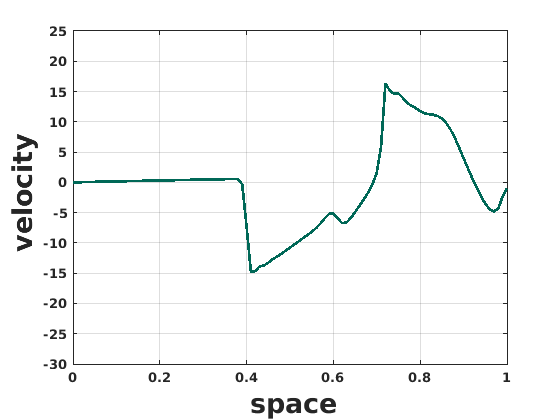}
\begin{center} $ t = 9.0 $ \end{center}
\end{minipage}
\begin{minipage}{0.24\linewidth}
\includegraphics[trim = 0cm 0cm 0cm 1cm, clip, scale=0.25]{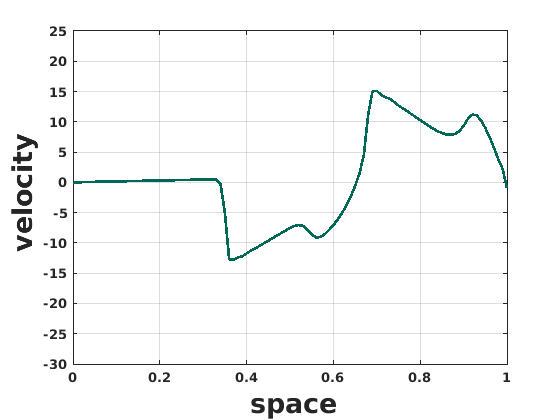}
\begin{center} $ t = 10.0 $ \end{center}
\end{minipage}
\\
\begin{minipage}{0.24\linewidth}
\includegraphics[trim = 0cm 0cm 0cm 1cm, clip, scale=0.25]{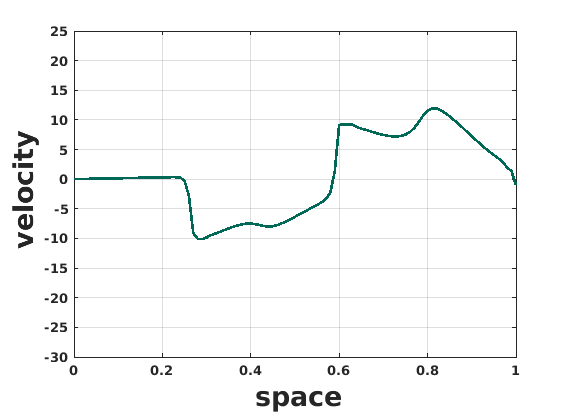}
\begin{center} $ t = 12.0 $ \end{center}
\end{minipage}
\begin{minipage}{0.24\linewidth}
\includegraphics[trim = 0cm 0cm 0cm 1cm, clip, scale=0.25]{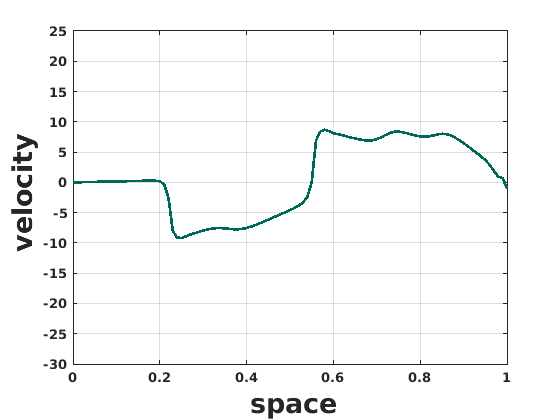}
\begin{center} $ t = 13.0 $ \end{center}
\end{minipage}
\begin{minipage}{0.24\linewidth}
\includegraphics[trim = 0cm 0cm 0cm 1cm, clip, scale=0.25]{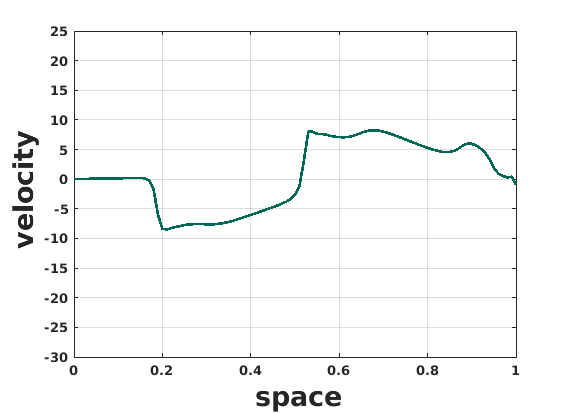}
\begin{center} $ t = 14.0 $ \end{center}
\end{minipage}
\begin{minipage}{0.24\linewidth}
\includegraphics[trim = 0cm 0cm 0cm 1cm, clip, scale=0.25]{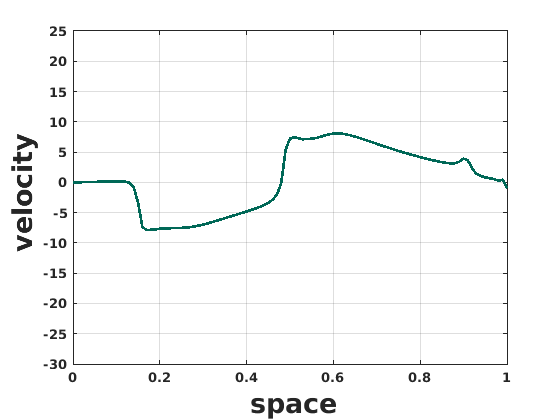}
\begin{center} $ t = 15.0 $ \end{center}
\end{minipage}
\end{minipage}
\begin{minipage}{\linewidth}
\begin{figure}[H]
\caption{Values of the displacement velocity variable~$\dot{u}$, for for \textcolor{black}{different times}. }
\label{figVelo}
\end{figure}
\end{minipage}\\
\hfill \\
\FloatBarrier

The influence of the control is more visible on the time evolution of the velocity field~$\dot{u}$ represented in Figure~\ref{figVelo}: At the beginning the state~$\dot{u}$ becomes non-positive, before changing its sign in order to create the profile of an oscillation. Next, after~$t=\tau$, it adopts a profile that is translated to the left.

\hspace*{-20pt}\begin{minipage}{\linewidth}
\begin{figure} [H]
\centering
\includegraphics[trim = 0cm 2cm 0cm 0cm, scale = 0.30]{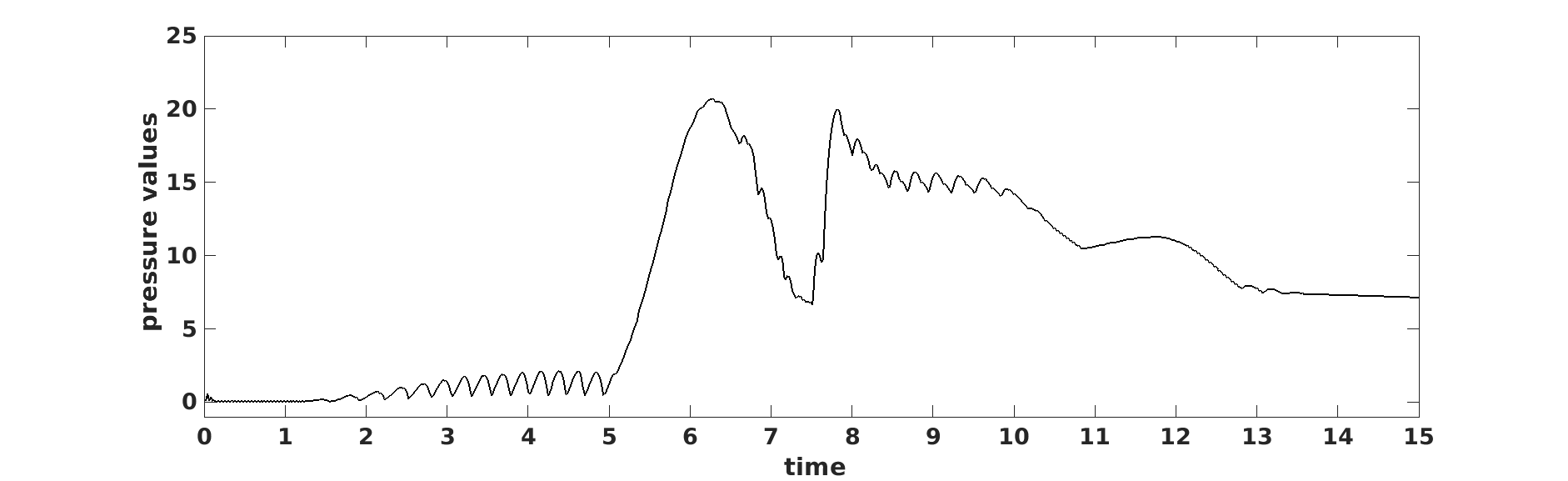}
\caption{Values of the pressure variable.}
\label{fig-pressure}
\end{figure}
\end{minipage}
\hfill \\ \hfill \\
\FloatBarrier
In Figure~\ref{fig-pressure} we observe that under the action of the control, the pressure presents small \textcolor{black}{successive oscillations, that we can maybe relate to the comments given in~\cite{Page2011}. These oscillations occur} until approximately~$t=5$ when it starts increasing steeply, for reaching it first maximum around~$t \approx 6.3$. Next the pressure decreases, before bouncing again around time~$t=5$, increasing next even more steeply for reaching another maximum around~$t=\tau \approx7.9$. Next, the pressure oscillates while decreasing, and seems to reach a pseudo-steady state from~$t=13.5$. The final value of the pressure is still larger than its initial value (that is~$0$).

\section{Conclusion} \label{sec-conclusion}
In this article we proposed a mathematical framework for the modeling of the mechanical aspects of defibrillation, based on the application of a distributed control on a part of the heart tissue. In particular, we developed an approach based on the optimal control theory, in order to enable us to maximize a class of functionals at a free time parameter, also optimized as well as the distributed control. We were able to derive rigorously first-order optimality conditions, that are possible to exploit for numerical realizations. We believe that our approach can pave the way to the development of robust and realistic numerical realizations. Further, based on the mathematical analysis we provide for the elastodynamics system with global injectivity condition, other physics-related aspects of the defibrillation problem could also be coupled to it, and addressed in the same fashion, like the electrical activity of the heart tissue for example. Let us finally mention that the complexity and the inherent technicalities seem {\it a priori} to be unavailable, as we aim at modeling phenomena with hyperelastic behavior.

\subsection*{Link to the code for numerical implementation}
The C++ code with which the numerical experiments were performed in section~\ref{sec-results} is available here:\\
\url{https://github.com/SebastienCourt/Defibrillation}

\subsection*{Acknowledgments}
The author thanks Prof.~Karl Kunisch (University of Graz \& RICAM, Linz), Prof.~Gernot Plank (Biotechmed Graz) and his research group. The discussions and ideas he had during the years he spent in Graz have \textcolor{black}{inspired} the present work. 

\appendix

\section{Appendix} \label{sec-apendix}

\subsection{Proofs of intermediate results}

This subsection is dedicated to the technical proofs of results related to the wellposedness of linear systems, namely Proposition~\ref{prop-well-auto} and Corollary~\ref{coro-well-auto-NH}. The methodology is similar to~\cite{Court2023}, but differs by its approach, as it actually deals with different systems.

\subsubsection{Proof of Proposition~\ref{prop-well-auto}} \label{sec-app-tek1}

Introduce first the following Hilbert spaces 
\begin{equation*}
\begin{array} {l}
 \mathcal{V}(\Omega) =  \left\{ v \in \HH^1(\Omega) \mid v_{|\Gamma_D} = 0, \ \displaystyle \int_{\Gamma_N} v\cdot n \, \d \Gamma_N = 0 \right\}, \\ \mathcal{V}(\Gamma_N) =  \left\{ v \in \HH^{1/2}(\Gamma_N) \mid \displaystyle \int_{\Gamma_N} v\cdot n \, \d \Gamma_N = 0 \right\}.
\end{array}
\end{equation*}
Remark that these spaces are included in those introduced in Definition~\ref{def-ass-op}, namely $\mathcal{V}_0(\Omega)$ and~$\mathcal{V}_0(\Gamma_N)$, respectively. We characterize the dual of $\mathcal{V}(\Gamma_N)$:
\begin{lemma} \label{lemma-dual}
A function~$\chi \in \mathcal{V}(\Gamma_N)'$ if and only if there exists~$\mathfrak{p} \in \R$ such that $\chi = \mathfrak{p}\, n$.
\end{lemma}

\begin{proof}
Set $\mathfrak{p} = \displaystyle \frac{1}{|\Gamma_N|} \int_{\Gamma_N} \chi \cdot n\, \d \Gamma_N$. It is easy to verify that $\chi -\mathfrak{p}\, n = 0$ in~$\mathbf{H}^{-1/2}(\Gamma_N)$, which ends the proof.
\end{proof}
\noindent Define the following operators
\begin{equation*}
\langle Av ; \varphi \rangle_{\mathcal{V}(\Omega)', \mathcal{V}(\Omega)} = \displaystyle
\kappa\int_{\Omega} \nabla v : \nabla \varphi \, \d \Omega, \quad
\langle Bw; \psi\rangle_{\mathcal{V}(\Omega)', \mathcal{V}(\Omega)} = \displaystyle
 \int_{\Omega} (\sigma_L(0).\nabla w):\nabla \psi \, \d \Omega.
\end{equation*}
Denoting $y_0=u$, $y_1 = \dot{u}$ and $y:=(y_0,y_1)^T$, the variational formulation of system~\eqref{syslih} is given for all $\varphi = (\varphi_0,\varphi_1) \in \LL^2(\Omega) \times\mathcal{V}(\Omega)$ as follows:\begin{small}
\begin{eqnarray}
\langle\dot{y}, \varphi\rangle_{\LL^2(\Omega)\times \mathcal{V}(\Omega)', \LL^2(\Omega)\times \mathcal{V}(\Omega)} & = &
\langle \dot{y}_0,\varphi_0 \rangle_{\LL^2(\Omega)} 
+ \langle \dot{y}_1,\varphi_1 \rangle_{\mathcal{V}(\Omega)',\mathcal{V}(\Omega)} \nonumber\\
& = & \langle y_1,\varphi_0 \rangle_{\LL^2(\Omega)}
- \kappa \langle \nabla y_1, \nabla \varphi_1\rangle_{\LLL^2(\Omega)}
-\langle \sigma_L(0).\nabla y_0,\nabla \varphi_1 \rangle_{\LLL^2(\Omega)} \nonumber \\
& & 
+ \langle f, \varphi_1\rangle_{\mathcal{V}(\Omega)',\mathcal{V}(\Omega)}
+ \langle g, \varphi_1 \rangle_{\mathcal{V}(\Gamma_N)',\mathcal{V}(\Gamma_N)}, \nonumber\\
\langle\dot{y}, \varphi\rangle_{\LL^2(\Omega)\times \mathcal{V}(\Omega)', \LL^2(\Omega)\times \mathcal{V}(\Omega)}
& = & \langle y_1,\varphi_0 \rangle_{\LL^2(\Omega)}
- \langle A y_1,  \varphi_1\rangle_{\mathcal{V}(\Omega)', \mathcal{V}(\Omega)}
-\langle B y_0, \varphi_1 \rangle_{\mathcal{V}(\Omega)', \mathcal{V}(\Omega)} \nonumber\\
& & 
+ \langle f, \varphi_1\rangle_{\mathcal{V}(\Omega)',\mathcal{V}(\Omega)}
+ \langle g, \varphi_1 \rangle_{\mathcal{V}(\Gamma_N)',\mathcal{V}(\Gamma_N)}.
\label{varform}
\end{eqnarray}\end{small}\noindent Using Assumption~$\mathbf{A3}$, 
we deduce that there exists $y_1 = \dot{u}\in \L^2(0,T;\mathcal{V}(\Omega)) \cap \H^1(0,T; \mathcal{V}(\Omega)')$ such that~\eqref{varform} holds for all $\varphi \in \mathcal{V}(\Omega)$, almost everywhere in $(0,T)$. After integration by parts, we obtain
\begin{equation*}
\begin{array} {l}
\ddot{u} -\kappa \Delta \dot{u} - \divg(\sigma_L(0).\nabla u)- f = 0  \quad \text{in } \mathcal{V}(\Omega)', \\
 \displaystyle \kappa \frac{\p \dot{u}}{\p n} + (\sigma_L(0).\nabla u)n- g = 0 \quad \text{in } \mathcal{V}(\Gamma_N)'.
\end{array}
\end{equation*}
From the second equation we deduce by Lemma~\ref{lemma-dual} the existence of~$\mathfrak{p}(t) 
\in \R$ such that
\begin{equation*}
\kappa \frac{\p \dot{u}}{\p n} + (\sigma_L(0).\nabla u)n + \mathfrak{p}\, n- g = 0
\quad \text{in } \HH^{-1/2}(\Gamma_N).
\end{equation*}
Note that deriving in time the third and fourth equations of system~\eqref{syslih} yields~$\langle \dot{u},n\rangle_{\HH^{-1/2}(\Gamma_N),\HH^{1/2}(\Gamma_N)} = 0$ and $\dot{u} = 0$ on~$\Gamma_D$ for almost every~$t\in[0,T]$. Taking the scalar product of the first equation of system~\eqref{syslih} by~$\dot{u}$ and using the Green formula, we obtain
\begin{equation*}
\frac{\d}{\d t} \left( \frac{1}{2}\|\dot{u}\|_{\LL^2(\Omega)}^2 + \sigma_L(0).\nabla u: \nabla u \right)
+ \kappa \| \nabla \dot{u}\|_{\LLL^2(\Omega)}^2 = 
\langle f,\dot{u} \rangle_{\LL^2(\Omega)} + \langle g,\dot{u} \rangle_{\HH^{1/2}(\Gamma_N),\HH^{-1/2}(\Gamma_N)}.
\end{equation*}
Integrating in time this equality and using Assumption~$\mathbf{A2}$, we deduce via the Young's inequality
\begin{equation*}
\nabla u \in \L^{\infty}(0,T;\LLL^2(\Omega)), \quad
\dot{u} \in \L^{\infty}(0,T;\LL^2(\Omega)), \quad
\nabla \dot{u} \in \L^2(0,T;\LLL^2(\Omega)).
\end{equation*}
Next, proceeding as previously with~$\ddot{u}$ in the role of~$\dot{u}$, we obtain similarly
\begin{equation*}
\begin{array}{l}
\displaystyle\|\ddot{u}\|_{\L^2(0,T;\LL^2(\Omega))}^2 + 
\frac{\kappa}{2} \|\nabla \dot{u}(T)\|_{\LLL^2(\Omega)}^2 +
\displaystyle\int_0^T \langle\sigma_L(0).\nabla u , \nabla \ddot{u}\rangle_{\mathcal{V}(\Omega),\mathcal{V}(\Omega)'} \d t\\ 
= \displaystyle \frac{\kappa}{2} \|\dot{u}_0\|_{\LLL^2(\Omega)}^2 
+ \langle f,\ddot{u}\rangle_{\LL^2(\Omega)}
 + \langle g,\ddot{u} \rangle_{\HH^{1/2}(\Gamma_N),\HH^{-1/2}(\Gamma_N)}.
\end{array}
\end{equation*}
Further, by integration by parts on~$(0,T)$, we deduce
\begin{equation*}
\begin{array} {rcl}
\|\ddot{u}\|_{\L^2(0,T;\LL^2(\Omega))}^2 + 
\frac{\kappa}{2} \frac{\d }{\d t} \|\nabla \dot{u}\|_{\LLL^2(\Omega)}^2 & = &
\displaystyle
\int_0^T \langle\sigma_L(0).\nabla \dot{u} , \nabla \dot{u}\rangle_{\LLL^2(\Omega)} \d t \\[10pt] & &
+ \langle f,\ddot{u}\rangle_{\LL^2(\Omega)} 
+ \langle g,\ddot{u} \rangle_{\HH^{1/2}(\Gamma_N),\HH^{-1/2}(\Gamma_N)}.
\end{array}
\end{equation*}
which shows that
\begin{equation*}
\ddot{u} \in \L^2(0,T;\LL^2(\Omega)), \quad 
\nabla \dot{u} \in \L^{\infty}(0,T;\LLL^2(\Omega)).
\end{equation*}
Further, from the first equation of~\eqref{syslih} we have $-\kappa \Delta \dot{u} = -\ddot{u} + \divg(\sigma_L(0).\nabla u) + f \in \L^2(0,T;\LL^2(\Omega))$, which yields $\dot{u} \in \L^2(0,T;\HH^2(\Omega))$ and therefore
\begin{equation*}
\dot{u} \in \L^2(0,T;\HH^2(\Omega)) \cap \H^1(0,T;\LL^2(\Omega)).
\end{equation*}
Thus system~\eqref{syslih} owns the $L^p$-maximal regularity property for $p=2$, and so for any~$p>3$, namely $\dot{u} \in \L^p(0,T;\HH^2(\Omega)) \cap \W^{1,p}(0,T;\LL^2(\Omega))$ and\begin{small}
\begin{equation*}
\begin{array} {l}
\|\dot{u}\|_{\L^p(0,T;\HH^2(\Omega)) \cap \W^{1,p}(0,T;\LL^2(\Omega))} \leq \\
\displaystyle
C \left(
\|f\|_{\L^p(0,T;\mathbf{L}^2(\Omega))}
+\|g\|_{\L^p(0,T;\mathbf{H}^{1/2}(\Gamma_N))\cap\W^{1/(2p'),p}(0,T;\mathbf{L}^2(\Gamma_N))} 
+ \|(u_0,\dot{u}_0)\|_{\mathbf{H}^2(\Omega)\times \mathbf{H}^{2/(p')}(\Omega)}
\right).
\end{array}
\end{equation*}\end{small}\noindent Further, estimate~\eqref{est-trace-emb3} yields also\begin{small}
\begin{equation}
\begin{array} {l}
\displaystyle
\left\|\frac{\p \dot{u}}{\p n} \right\|_{\W^{1/(2p'),p}(0,T;\mathbf{H}^{1/2-1/p}(\Gamma_N)')}  \leq \\
C \left(\|f\|_{\L^p(0,T;\mathbf{L}^2(\Omega))}
+\|g\|_{\L^p(0,T;\mathbf{H}^{1/2}(\Gamma_N))\cap\W^{1/(2p'),p}(0,T;\mathbf{L}^2(\Gamma_N))}
+ \|(u_0,\dot{u}_0)\|_{\mathbf{H}^2(\Omega)\times \mathbf{H}^{2/(p')}(\Omega)}
\right).
\end{array}
\label{est-norm-deriv}
\end{equation}\end{small}\noindent Since the coefficients of~$\sigma_L(0)$ are in~$\L^{\infty}(\Omega)$, we also have\begin{small}
\begin{equation}
\begin{array} {l}
\displaystyle
\|(\sigma_L(0).\nabla u)n \|_{\W^{1/(2p'),p}(0,T;\mathbf{H}^{1/2-1/p}(\Gamma_N)')} \leq \\
C \left(
\|f\|_{\L^p(0,T;\mathbf{L}^2(\Omega))}
+\|g\|_{\L^p(0,T;\mathbf{H}^{1/2}(\Gamma_N))\cap\W^{1/(2p'),p}(0,T;\mathbf{L}^2(\Gamma_N))}
+ \|(u_0,\dot{u}_0)\|_{\mathbf{H}^2(\Omega)\times \mathbf{H}^{2/(p')}(\Omega)}
\right).
\label{est-norm-deriv2}
\end{array}
\end{equation}\end{small}\noindent Next, using the second equation of system~\eqref{syslih}, we obtain the following expression for the pressure
\begin{equation*}
\mathfrak{p} = \displaystyle
 \frac{1}{|\Gamma_N|}\left(
\int_{\Gamma_N} g \cdot n \, \d \Gamma_N 
-\kappa\left\langle\frac{\p \dot{u}}{\p n}+ (\sigma_L(0).\nabla u)n; n\right\rangle_{\HH^{1/2-1/p}(\Gamma_N)',\HH^{1/2-1/p}(\Gamma_N)}
\right).
\end{equation*}
Since the normal vector of~$\Gamma_N$ satisfies
\begin{equation*}
n\in \mathbf{W}^{2-1/p,p}(\Gamma_N) \hookrightarrow \mathbf{H}^{2-1/p}(\Gamma_N) \hookrightarrow \mathbf{H}^{1/2-1/p}(\Gamma_N),
\end{equation*}
we deduce that the pressure satisfies $\mathfrak{p} \in \W^{1/(2p'),p}(0,T;\R)$, and\begin{small}
\begin{equation*}
\begin{array} {rcl}
\displaystyle
\| \mathfrak{p}\|_{\W^{1/(2p'),p}(0,T;\R)} & \leq & C\left(
\|g\|_{\W^{1/(2p'),p}(0,T;\LL^2(\Gamma_N))} \right. \\ & & + 
\left. \displaystyle\left\|\frac{\p \dot{u}}{\p n} + (\sigma_L(0).\nabla u)n
\right\|_{\W^{1/(2p'),p}(0,T;\HH^{1/2-1/p}(\Gamma_N)')}
\right).
\end{array}
\end{equation*}\end{small}\noindent Combined with estimates~\eqref{est-norm-deriv}-\eqref{est-norm-deriv2}, this yields
\begin{equation}
\begin{array} {rcl}
\|\mathfrak{p} \|_{\mathcal{P}_{p,T}} & \leq & C \left(
\|f\|_{\L^p(0,T;\mathbf{L}^2(\Omega))}
+\|g\|_{\L^p(0,T;\mathbf{H}^{1/2}(\Gamma_N))\cap\W^{1/(2p'),p}(0,T;\mathbf{L}^2(\Gamma_N))}\right. \\
& & \left.+ \|(u_0,\dot{u}_0)\|_{\mathbf{H}^2(\Omega)\times \mathbf{H}^{2/(p')}(\Omega)}
\right). \label{est-p-prop32}
\end{array}
\end{equation}
Therefore from the second equation of system~\eqref{syslih} we have $\displaystyle \kappa\frac{\p \dot{u}}{\p n} + (\sigma_L(0).\nabla u)n = g -\mathfrak{p}\, n \in \mathcal{G}_{p,T}(\Gamma_N)$, and then we deduce the regularity of~$(u,\dot{u})$ in $\mathcal{U}_{p,T}(\Omega) \times \dot{\mathcal{U}}_{p,T}(\Omega)$ from Proposition~\ref{prop-well-auto}. Further, estimate~\eqref{est-prop32} yields
\begin{equation*}
\begin{array} {l}
\|u\|_{\mathcal{U}_{p,T}(\Omega)} + 
\|\mathfrak{p}\|_{\mathcal{P}_{p,T}}\\[10pt]  \leq 
C\left(\|(u_0,\udotini)\|_{\mathcal{U}^{(0,1)}_p(\Omega)} +
\|f\|_{\mathcal{F}_{p,T}(\Omega)} + 
\|g\|_{\mathcal{G}_{p,T}(\Gamma_N)} +
\| \mathfrak{p}\, n\|_{\mathcal{G}_{p,T}(\Gamma_N)}
\right) \\
\\ \leq C\left(\|(u_0,\udotini)\|_{\mathcal{U}^{(0,1)}_p(\Omega)} +
\|f\|_{\mathcal{F}_{p,T}(\Omega)} + 
\|g\|_{\mathcal{G}_{p,T}(\Gamma_N)} +
\| \mathfrak{p}\|_{\mathcal{P}_{p,T}}
\right).
\end{array}
\end{equation*}
Combined with~\eqref{est-p-prop32}, this yields the announced estimate and completes the proof.

\subsubsection{Proof of Corollary~\ref{coro-well-auto-NH}} \label{sec-app-tek2}

Notice that the constraint of system~\eqref{siszeroNH} also writes
\begin{equation*}
\int_{\Gamma_N} \left(u-\frac{1}{|\Gamma_N|}hn\right)
\cdot n\, \d \Gamma_N = 0 \quad \quad \text{in $(0,T)$.}
\end{equation*}
We then proceed with a lifting method. We need to define an extension of~$\displaystyle \frac{1}{|\Gamma_N|}hn$ in~$\Omega$. Let us first define extensions of $\displaystyle \frac{1}{|\Gamma_N|}h(0)n$ and $\displaystyle \frac{1}{|\Gamma_N|}\dot{h}(0)n$. Recall that $h(0) \in \R$ and $\dot{h}(0)\in \R$ do not depend on the space variable. Since by assumption $n \in \WW^{2-1/p,p}(\Gamma_N) \hookrightarrow\WW^{2/{p'}-1/p,p}(\Gamma_N) $, there exists~$H_0 \in \WW^{2,p}(\Omega)$ and $\dot{H}_0 \in \WW^{2/{p'},p}(\Omega)$ respective extensions of $\displaystyle  \frac{1}{|\Gamma_N|} h(0)n$ and $\displaystyle  \frac{1}{|\Gamma_N|} \dot{h}(0)n$ such that
\begin{subequations} \label{lifteq}
\begin{eqnarray}
{H_0}_{| \Gamma_N} = \displaystyle \frac{1}{|\Gamma_N|} h(0)n, \qquad
\| H_0 \|_{\WW^{2,p}(\Omega)} 
\leq C\|h(0) \|_{\R},\label{lifteq1}\\[5pt]
{\dot{H}}_{0| \Gamma_N} = \displaystyle \frac{1}{|\Gamma_N|} \dot{h}(0)n, \qquad
\| \dot{H}_0 \|_{\WW^{2/{p'},p}(\Omega)} 
\leq C\|\dot{h}(0) \|_{\R}.\label{lifteq2}
\end{eqnarray}
\end{subequations}
We now define an extension $H$ of~$\displaystyle \frac{1}{|\Gamma_N|}hn$ by solving first the following heat equation with mixed boundary conditions, dealing with~$\dot{H}$ as unknown:
\begin{equation*}
\begin{array} {l}
\ddot{H}- \kappa\Delta \dot{H} = 0  \text{ in } \Omega \times (0,T), \quad
\dot{H}_{|\Gamma_N} = \displaystyle \frac{1}{|\Gamma_N|}\dot{h} n  
\  \text{on } \Gamma_N \times (0,T), \\[10pt]
\dot{H} = 0  \ \text{on } \Gamma_D \times (0,T), \quad
\dot{H}(0) = \dot{H}_0 \  \text{in } \Omega.
\end{array}
\end{equation*}
Since~$n \in \WW^{2-1/p,p}(\Gamma_N)$, we have $\displaystyle \frac{1}{|\Gamma_N|}\dot{h}n \in \mathcal{H}_{p,T}(\Gamma_N)$. From~\cite{Pruess2002} the solution of this equation satisfies
\begin{eqnarray}
\| \dot{H}\|_{\dot{\mathcal{U}}_{p,T}(\Omega)} & \leq & 
C\left( \| \dot{h}n\|_{\mathcal{H}_{p,T}(\Gamma_N)} + \| \dot{H}_0\|_{\WW^{2/{p'},p}(\Omega)}\right)\nonumber \\
& \leq &  C\left(
\|\dot{h} \|_{\W^{1-1/{2p},p}(0,T;\R)} + \|\dot{h}(0)\|_{\R}\right)
 = C\|\dot{h}\|_{\mathcal{H}_{p,T}},
\label{estliftx}
\end{eqnarray}
where we have used~\eqref{lifteq2}. Using~\eqref{est-trace-emb2}, we deduce in particular
\begin{equation}
\|\dot{H}(0)\|_{\mathbf{W}^{2/{p'},p}(\Omega)}
+\kappa\left\|\frac{\p \dot{H}}{\p n}\right\|_{\mathcal{G}_{p,T}(\Gamma_N)} \leq
C \|\dot{h} \|_{\mathcal{H}_{p,T}}.
\label{est-dirichlet-hbar}
\end{equation}
Further, we set $H(\cdot,t) = H_0 + \displaystyle \int_0^t \dot{H}(\cdot,s)\d s$, which implies
\begin{equation}
\|H\|_{\mathcal{U}_{p,T}(\Omega)} \leq C\left(
\|H_0\|_{\WW^{2,p}(\Omega)} + \|\dot{H}\|_{\dot{\mathcal{U}}_{p,T}(\Omega)}
\right) \leq C
\left(
\|h(0)\|_{\R} + \|\dot{h}\|_{\mathcal{H}_{p,T}}
\right), \label{est-HU}
\end{equation}
where we have used~\eqref{lifteq1} and~\eqref{estliftx}, and thus we obtain in particular
\begin{equation}
\|H\|_{\L^p(0,T;\WW^{2,p}(\Omega))}
+\|(\sigma_L(0).\nabla H)n\|_{\mathcal{G}_{p,T}(\Gamma_N)} \displaystyle
 \leq  C\left( 
\|h(0)\|_{\R} + \| \dot{h}\|_{\mathcal{H}_{p,T}}
\right).
\label{est-dirichlet-hbar2}
\end{equation}
Now define $\bar{u} := u - H$. We rewrite system~\eqref{siszeroNH} as 
\begin{equation}
\begin{array} {rcl} 
\ddot{\bar{u}} - \kappa \Delta \dot{\bar{u}} -\divg(\sigma_L(0).\nabla \bar{u}) = 
f + \divg(\sigma_L(0).\nabla H) & & 
\text{in $\Omega \times (0,T)$},\\
\kappa \displaystyle \frac{\p \dot{\bar{u}}}{\p n} + (\sigma_L(0).\nabla \bar{u}) n
+\mathfrak{p}\, n = 
g - \kappa\frac{\p \dot{H}}{\p n}- (\sigma_L(0).\nabla H)n & &  
\text{on $\Gamma_N \times (0,T)$},\\
\displaystyle
\int_{\Gamma_N} \bar{u}\cdot n\, \d \Gamma_N = 0 & & \text{in $(0,T)$} ,\\
\bar{u} = 0 & & \text{on $\Gamma_D \times (0,T)$},\\
u(\cdot,0) = u_0-h_0, \quad \dot{\bar{u}}(\cdot,0) = \udotini - \dot{h}_0 & &
\text{in $\Omega$},
\end{array}
\end{equation}
We recognize system~\eqref{syslih} satisfied by~$(\bar{u},\mathfrak{p})$. Then from Proposition~\ref{prop-well-auto} there exists a unique~$(\bar{u},\mathfrak{p}) \in \mathcal{U}_{p,T}(\Omega) \times \mathcal{P}_{p,T}$, and it satisfies
\begin{equation*}
\begin{array} {rcl}
\|\bar{u}\|_{\mathcal{U}_{p,T}(\Omega)} + \| \mathfrak{p} \|_{\mathcal{P}_{p,T}}
& \leq & \displaystyle
C \left(\displaystyle
\|f\|_{\mathcal{F}_{p,T}} + \|H\|_{\L^p(0,T;\WW^{2,p}(\Omega))} \right. \\
& & \displaystyle \left. + \| g\|_{\mathcal{G}_{p,T}(\Gamma_N)}
+\|(\sigma_L(0).\nabla H)n\|_{\mathcal{G}_{p,T}(\Gamma_N)}
+ \kappa\left\|\frac{\p \dot{H}}{\p n}\right\|_{\mathcal{G}_{p,T}(\Gamma_N)}
 \right. \\
& & \displaystyle \left.  + 
\|(u_0,\dot{u}_0)\|_{\mathcal{U}_{p}^{(0,1)}(\Omega)} +
\|(H_0,\dot{H}_0)\|_{\mathcal{U}_{p}^{(0,1)}(\Omega)}
\right), \\
& \leq & C\left( \displaystyle
\|f\|_{\mathcal{F}_{p,T}} + \| g\|_{\mathcal{G}_{p,T}(\Gamma_N)} +
\|h(0)\|_{\R} + \| \dot{h}\|_{\mathcal{H}_{p,T}}
\right. \\
& & \left. +
\|(u_0,\dot{u}_0)\|_{\mathcal{U}_{p}^{(0,1)}(\Omega)} \right)
\end{array}
\end{equation*}
where we have used~\eqref{est-dirichlet-hbar}-\eqref{est-dirichlet-hbar2} and~\eqref{lifteq1}-\eqref{lifteq2}. Further, the couple $(u,\mathfrak{p}) = (\bar{u}+ H,\mathfrak{p})$ satisfies
\begin{equation*}
\|u\|_{\mathcal{U}_{p,T}(\Omega)} + \|\mathfrak{p}\|_{\mathcal{P}_{p,T}}
\leq 
\|\bar{u}\|_{\mathcal{U}_{p,T}(\Omega)}
 + \|\mathfrak{p}\|_{\mathcal{P}_{p,T}}
+ \|H\|_{\mathcal{U}_{p,T}(\Omega)},
\end{equation*}
which, combined with the previous estimate and~\eqref{est-HU}, yields
\begin{equation*}
\begin{array}{rcl}
\|u\|_{\mathcal{U}_{p,T}(\Omega)} + \|\mathfrak{p}\|_{\mathcal{P}_{p,T}} 
& \leq & 
C\left( \displaystyle
\|f\|_{\mathcal{F}_{p,T}} + \| g\|_{\mathcal{G}_{p,T}(\Gamma_N)} +
 \|h\|_{\mathcal{H}_{p,T}}
+\|(u_0,\dot{u}_0)\|_{\mathcal{U}_{p}^{(0,1)}(\Omega)}
\right).
\end{array}
\end{equation*}
This provides the existence of solution~$(u,\mathfrak{p})$. Uniqueness is due to the linearity of system~\eqref{siszeroNH}: Considering two solutions of~\eqref{siszeroNH}, their difference satisfies system~\eqref{syslih} with zero right-hand-sides, and therefore is equal to zero, because solutions of~\eqref{syslih} are unique, which concludes the proof.

\subsection{Examples of strain energies} \label{sec-app-energies}
Let us give a set of examples of classical strain energies from the literature, and show that they satisfy the assumptions~$\mathbf{A1}$--$\mathbf{A3}$.

\subsubsection{The Saint Venant-Kirchhoff's model.} It corresponds to the following strain energy
\begin{equation*}
\mathcal{W}_1(E) = \mu_L \trace\left( E^2 \right) + \frac{\lambda_L}{2} \trace(E)^2,
\end{equation*}
where $\mu_L >0$ and $\lambda_L \geq 0$ are the so-called Lam\'e coefficients. The energy is clearly twice differentiable, its first- and second-derivatives of $\mathcal{W}_1$ are given respectively by
\begin{equation*}
	\check{\Sigma}_1(E) = 2\mu_L E + \lambda_L \trace(E)\I, \qquad
\frac{\p^2\mathcal{W}_1}{\p E^2}(E) = \frac{\p \check{\Sigma}_1}{\p E}(E)  =  2\mu_L \III + \lambda_L \I \otimes \I,
\end{equation*}
where $\III \in \R^{d\times d \times d \times d}$ denotes the identity tensor of order~$4$. In particular, we see that if the matrix~$E$ is symmetric, then~$\check{\Sigma}_1$ defines a symmetric matrix. Therefore Assumptions~$\mathbf{A1}$-$\mathbf{A2}$ are satisfied by this strain energy. Further, regarding Assumption~$\mathbf{A3}$, we see that
\begin{equation*}
\sigma_L(0).\nabla v = 
\nabla v \check{\Sigma}_1(0) + \displaystyle \frac{1}{2}
\frac{\p^2\mathcal{W}_1}{\p E^2}(0).(\nabla v + \nabla v^T) = 
2\mu_L\epsilon(v) + \lambda_L \mathrm{trace}(\epsilon(v)) \I.
\end{equation*}
where we have introduced the notation~$\epsilon(v) = \frac{1}{2}(\nabla v + \nabla v^T)$ for the symmetric part of~$\nabla v$. The operator~$\sigma_L(0)$ then corresponds to the well-known linearized Lam\'e operator, which defines the coercive operator~$-\divg(\sigma_L(0).\nabla v)$ under the condition~$v_{|\Gamma_D} = 0$, in virtue of the Petree-Tartar lemma~\cite[Lemma~A.38 page~469]{Ern}. We can also refer to the Korn's inequality for this claim. Therefore we can claim that Assumption~$\mathrm{A3}$ is satisfied for this example.

\subsubsection{The Fung's model.} It corresponds to the following strain energy
\begin{equation*}
	\mathcal{W}_2(E) = \mathcal{W}_2(0) + \beta \left(\exp\left(\gamma \ \trace(E^2)\right) - 1\right),
\end{equation*}
where $\mathcal{W}_2(0) \geq 0$, $\beta >0$ and $\gamma > 0$ are given coefficients. The space $\WWW^{1,p}(\Omega)$ is invariant under composition of the exponential function when $p>d$ (see~\cite{BB1974}, Lemma~A.2. page 359). The first- and second-derivatives of $\mathcal{W}_2$ are given respectively by
\begin{equation*}
	\check{\Sigma}_2(E)  =  2\gamma\beta \exp\left(\gamma \ \trace(E^2) \right) E, \qquad
	\frac{\p \check{\Sigma}_2}{\p E}(E)  =  \beta \exp\left(\gamma \ \trace(E^2) \right)
	\left( 2\gamma \III + (2\gamma)^2 E \otimes E\right).
\end{equation*}
Again, if~$E$ is symmetric, then~$\check{\Sigma}_2$ is symmetric. Assumptions~$\mathbf{A1}$-$\mathbf{A2}$ are then satisfied. For Assumption~$\mathbf{A3}$, we need to evaluate
\begin{equation*}
\sigma_L(0).\nabla v = 
\nabla v \check{\Sigma}_2(0) + \displaystyle 
\frac{\p^2\mathcal{W}_2}{\p E^2}(0).(\epsilon(v)) = 
2\beta \gamma \epsilon(v).
\end{equation*}
Like in the previous example, the operator~$-\divg(\sigma_L(0).\nabla v)$ is coercive, and thus Assumption~$\mathbf{A3}$ is satisfied.

\subsubsection{The Ogden's model.} The family of strain energies corresponding to this model are linear combinations of energies of the following form
\begin{equation*}
	\mathcal{W}_3(E)  =   \trace\left((2E+\I)^{\gamma} - \I \right),
\end{equation*}
where $\gamma \in \R$. Since the tensor $2E+\I$ is real and symmetric, the expression $(2E+\I)^{\beta} $ makes sense for any $\beta \in \R$ by diagonalizing~$2E+\I$, and the energy $\mathcal{W}_3(E)$ can be expressed in terms of the eigenvalues of $2E+\I$. Since $2E(u)+\I = (\I+\nabla u)^T(\I+\nabla u)$, if $(\lambda_i)_{1\leq i \leq d}$ denote the singular values of $\I+\nabla u$, and $(\mu_i)_{1\leq i \leq d}$ denote those of $E(u)$, we have
\begin{equation*}
	\mathcal{W}_3(E)  =  \sum_{i=1}^d \left(\lambda_i^{2\gamma}  -1\right)
	= \sum_{i=1}^d \left((1+2\mu_i)^\gamma -1\right), \qquad
	\check{\Sigma}_3(E) =  2\gamma (2E+\I)^{\gamma-1}.
\end{equation*}
Denoting by $(v_l)_{1\leq l\leq d}$ the normalized orthogonal eigenvectors of $E$, we further write
\begin{equation*}
\check{\Sigma}_3(E)  = \sum_{i=1}^d 2\gamma(2\mu_i+1)^{\gamma-1}v_i \otimes v_i.
\end{equation*}
Note that the operator $v_i \otimes v_i$ is the projection on $\mathrm{Span}(v_i)$, and Assumption~$\mathbf{A2}$ is satisfied by~$\check{\Sigma}_3$. The sensitivity of the eigenvalues and eigenvectors with respect to the matrix can be derived for example from~\cite{Stewart} (see Theorem~IV.2.3 page~183, and Remark~2.9 page 239). Thus after calculations we get
\begin{equation*}
\begin{array} {rcl}
\displaystyle \frac{\p \check{\Sigma}_3}{\p E}(E) & = & \displaystyle 2\gamma
\sum_{i=1}^d 2(\gamma-1)(2\mu_i+1)^{\gamma-2} (v_i\otimes v_i) \otimes (v_i\otimes v_i) \displaystyle \\
& & + \displaystyle  2\gamma\sum_{i=1}^d (2\mu_i+1)^{\gamma-1} 
\sum_{j\neq i} \frac{1}{\mu_i-\mu_j}(v_j\otimes v_i) \otimes (v_j\otimes v_i) .
\end{array}
\end{equation*}
This expression shows that the strain energy $\mathcal{W}_3$ fulfills also Assumption~$\mathbf{A1}$. Finally, considering the vectors~$(v_i)_{1\leq i\leq d}$ of the canonical basis as normalized orthogonal eigenvectors of matrix $E=0$ (with eigenvalues $\mu_i = 0$), we evaluate\begin{small}
\begin{equation*}
\begin{array} {l}
\sigma_L(0).\nabla u = 
\nabla u \check{\Sigma}_3(0) + \displaystyle 
\frac{\p\check{\Sigma}_3}{\p E}(0).(\epsilon(u)) = 
2\gamma \nabla u + 
\displaystyle 4\gamma(\gamma-1)
\sum_{i=1}^d \Big((v_i\otimes v_i):\epsilon(u) \Big) (v_i\otimes v_i), \\
(\sigma_L(0).\nabla u) : \nabla u  = 
2\gamma |\nabla u|^2_{\R^{d\times d}}
+4\gamma(\gamma-1) \displaystyle 
\sum_{i=1}^d \Big((v_i\otimes v_i):\epsilon(u) \Big)^2.
\end{array}
\end{equation*}\end{small}\noindent We then obtain a coercive operator, provided that for example~$\gamma\geq 1$, and deduce that Assumption~$\mathbf{A3}$ is also satisfied for this model.

\subsection{A Lagrangian mechanics perspective} \label{sec-app-Lag}

Introduce formally the Lagrangian functional associated with Problem~\eqref{mainpbtilde2} as follows:
\begin{small}
\begin{equation*}
\begin{array} {rcl}
\tilde{\mathcal{L}}(\tilde{y}_0,\tilde{y}_1,\tilde{\mathfrak{p}}, \tilde{\xi}, \tau,\tilde{\zeta}_0, \tilde{\zeta}_1, \tilde{\pi}) & = &
\displaystyle 
\int_0^2 \dot{\mu}c(\tilde{y}_0,\tilde{y}_1,\tilde{\xi})\d s +
 \phi^{(1)}(\tilde{y}_0,\tilde{y}_1)(1) + \phi^{(2)}(\tilde{y}_0,\tilde{y}_1)(2) \\[10pt]
& & +\displaystyle 
\int_0^2 
\langle \dot{\tilde{y}}_1 - \dot{\mu} f(\tilde{\xi}) , \tilde{\zeta}_1\rangle_{\LL^p(\Omega),\LL^{p'}(\Omega)} \d s \\[10pt]
& & \displaystyle +\int_0^2 \dot{\mu}\langle \kappa \nabla \tilde{y}_1 + 
\sigma(\nabla \tilde{y}_0), 
\nabla \tilde{\zeta}_1\rangle_{\WW^{1,p}(\Omega),\WW^{1,p}(\Omega)'}
\d s\\[10pt]
& & \displaystyle
+ \langle \dot{\tilde{y}}_0-\dot{\mu}\tilde{y}_1,\tilde{\zeta}_0 \rangle_{\dot{\mathcal{U}}_{p,2}(\Omega), \dot{\mathcal{U}}_{p,2}(\Omega)'} 
- \langle \tilde{g},\tilde{\zeta}_1\rangle_
{\mathcal{G}_{p,2}(\Gamma_N),\mathcal{G}_{p,2}(\Gamma_N)'}\\[10pt]
& & \displaystyle + \int_0^2 \dot{\mu}\tilde{\pi} \int_{\Omega} \det \Phi(\tilde{y}_0) \, 
\d \Omega\, \d s \\[10pt]
& & +\displaystyle \int_0^2 \dot{\mu}\tilde{\mathfrak{p}}\int_{\Gamma_N}  
\tilde{\zeta}_1 \cdot \cof(\Phi(\tilde{y}_0))n \, \d \Gamma_N\, \d s.
\end{array}
\end{equation*}
\end{small}\noindent Recall that the dependence of~$\tilde{\mathcal{L}}$ with respect to~$\tau$ is represented by the change of variables $\dot{\mu} = \dot{\mu}(\cdot,\tau)$. Coming back to the original variables, namely 
\begin{equation*}
\begin{array} {lll}
\tilde{y}_0(\cdot,s) = y_0(\cdot,\mu(s,\tau)), &
\tilde{y}_1(\cdot,s) = y_1(\cdot,\mu(s,\tau)), & 
\tilde{\mathfrak{p}}(s) = \mathfrak{p}(\mu(s,\tau)),\\[5pt]
\tilde{\xi}(\cdot,s) = \xi(\cdot,\mu(s,\tau)), & & \\[5pt]
\tilde{\zeta}_0(\cdot,s) = \zeta_0(\cdot,\mu(s,\tau)), &
\tilde{\zeta}_1(\cdot,s) = \zeta_1(\cdot,\mu(s,\tau)), &
\tilde{\pi}(s) = \pi(\mu(s,\tau)), \\[5pt]
\tilde{f}(\cdot,s) = f(\cdot,\mu(s,\tau)), &
\tilde{g}(\cdot,s) = g(\cdot,\mu(s,\tau)), & 
\end{array}
\end{equation*}
we have
\begin{equation*}
\tilde{\mathcal{L}}(\tilde{y}_0,\tilde{y}_1,\tilde{\mathfrak{p}}, \tilde{\xi}, \tau,\tilde{\zeta}_0, \tilde{\zeta}_1, \tilde{\pi}) = 
\mathcal{L}(y_0,y_1,\mathfrak{p}, \xi, \tau,\zeta_0, \zeta_1, \pi),
\end{equation*}
where
\begin{small}
\begin{equation*}
\begin{array} {rcl}
\mathcal{L}(y_0,y_1,\mathfrak{p}, \xi,\tau,\zeta_0, \zeta_1, \pi) & = &
\displaystyle 
\int_0^T c(y_0,y_1,\xi)\d t + \phi^{(1)}(y_0,y_1)(\tau) + \phi^{(2)}(y_0,y_1)(T) \\[10pt]
& & + \displaystyle 
\int_0^T \langle \dot{y}_1 - f(\xi), \zeta_1\rangle_{\LL^p(\Omega),\LL^{p'}(\Omega)} \d t\\[10pt]
& & + \displaystyle \int_0^T\langle \kappa \nabla y_1 + \sigma(\nabla y_0), \nabla \zeta_1
 \rangle_{\WW^{1,p}(\Omega),\WW^{1,p}(\Omega)'}
\d t\\[10pt]
& & \displaystyle
+ \langle \dot{y}_0-y_1,\zeta_0 \rangle_{\dot{\mathcal{U}}_{p,T}(\Omega), \dot{\mathcal{U}}_{p,T}(\Omega)'} 
- \langle g,\zeta_1\rangle_{\mathcal{G}_{p,T}(\Gamma_N),\mathcal{G}_{p,T}(\Gamma_N)'}\\[10pt]
& & \displaystyle + \int_0^T \pi \int_{\Omega} \det \Phi(y_0) \, \d \Omega\, \d t
+ \int_0^T \mathfrak{p}\int_{\Gamma_N}  \zeta_1 \cdot \cof(\Phi(y_0))n \, \d \Gamma_N\, \d t.
\end{array}
\end{equation*}
\end{small}\noindent Differentiating~$\tilde{\mathcal{L}}$ with respect to~$(\tilde{\zeta}_0,\tilde{\zeta}_1,\tilde{\pi})$ yields system~\eqref{sysmaintilde2x}. Differentiating~$\tilde{\mathcal{L}}$ with respect to~$(\tilde{y}_0,\tilde{y}_1,\tilde{\mathfrak{p}})$ yields system~\eqref{sysadjoint}. And differentiating~$\tilde{\mathcal{L}}$ with respect to~$(\tilde{\xi},\tau)$ yields~\eqref{id-deriv}. Therefore a critical point of functional~$\tilde{\mathcal{L}}$ satisfies the optimality conditions stated in Theorem~\ref{th-optcond}. Actually, following the approach adopted in~\cite{Maxmax1}, we could show that an optimal solution of Problem~\eqref{mainpbtilde} is necessarily a critical point of function~$\tilde{\mathcal{L}}$. Further, the optimality conditions stated as in Corollary~\ref{cor-optcond} correspond to a critical point of mapping~$\mathcal{L}$.

\begin{remark} \label{remark-sym}
Let us detail the term of~$\mathcal{L}$ derived from the strain energy, namely
\begin{equation*}
\begin{array} {rcl}
\langle  \sigma(\nabla y_0), \nabla \zeta_1
 \rangle_{\WW^{1,p}(\Omega),\WW^{1,p}(\Omega)'} & = & 
\langle(\I+\nabla y_0) \check{\Sigma}(E(y_0)), \nabla  \zeta_1 \rangle_{\WW^{1,p}(\Omega),\WW^{1,p}(\Omega)'} \\
& = & 
\langle \check{\Sigma}(E(y_0)), (\I+\nabla y_0)^T\nabla  \zeta_1 \rangle_{\WW^{1,p}(\Omega),\WW^{1,p}(\Omega)'} \\
& = & \langle \check{\Sigma}(E(y_0)), E'(y_0).\zeta_1 \rangle_{\WW^{1,p}(\Omega),\WW^{1,p}(\Omega)'}
\end{array}
\end{equation*}
with $E'(y_0).\zeta_1 = \frac{1}{2}((\I+\nabla y_0)^T\nabla \zeta_1 + \nabla \zeta_1^T(\I+\nabla y_0)$, because the tensor~$\check{\Sigma}$ is assumed to be symmetric in Assumption~$\mathbf{A2}$. Further, recall that $\check{\Sigma}(E(y_0)) = \displaystyle \frac{\p \mathcal{W}}{\p E}(E(y_0))$, and thus
\begin{equation*}
\langle  \sigma(\nabla y_0), \nabla \zeta_1
 \rangle_{\WW^{1,p}(\Omega),\WW^{1,p}(\Omega)'} = 
\displaystyle \frac{\p }{\p y_0}(\mathcal{W}( E(y_0))).\zeta_1.
\end{equation*}
\end{remark}

\subsection{The control operator in the context of cardiac electrophysiology} 
\label{sec-app-control}
The control is realized through a distributed right-hand-side $f$ in equation~\eqref{sysmain1}. In practice this function is expressed in terms of the {\it fiber direction}, denoted by $\hat{\mathfrak{f}}$, namely a vector tangent to the tissue, depending on the geometry and considered as a part of the data. More precisely, $f$ is chosen in the form
\begin{equation*}
	f  =  \divg(s_a \hat{\mathfrak{f}} \otimes \hat{\mathfrak{f}} ),
\end{equation*}
where $s_a$ is a scalar function, depending on space and time, that we choose as being the command, denoted formally by $\xi$ throughout the paper. The tensor $s_a \hat{\mathfrak{f}} \otimes \hat{\mathfrak{f}}$ is the so-called {\it active stress tensor}. Since on $\p \Omega$ the vector $\hat{\mathfrak{f}}$ is tangent, by the Green formula the following inner product by any test function $\varphi$ writes simply as
\begin{equation*}
\langle f(\xi) ; \varphi \rangle_{\LL^2(\Omega)} = 	
-\int_{\Omega} \divg(\xi \hat{\mathfrak{f}} \otimes \hat{\mathfrak{f}})\cdot \varphi \, \d \Omega =
	\int_{\Omega} \xi( \hat{\mathfrak{f}}  \otimes \hat{\mathfrak{f}}  : \nabla \varphi)\, \d \Omega.
\end{equation*}
Denoting by $\omega \subset \Omega$ the control domain, an example of control space for the distributed control function $\xi$ on~$\omega$ is the following
\begin{equation*}
\mathcal{X}_{p,T}(\omega) = \L^p(0,T;\W^{1,p}(\omega)).
\end{equation*}
In this example the control function $\xi$ is only scalar, but since the quantity to maximize, namely the variations of the pressure $\mathfrak{p}$, is also scalar, and moreover depends only on time, we expect that the set of controls is rich enough for our purpose.

\end{document}